\documentclass{article}
\usepackage{latexsym}
\usepackage{graphics}
\usepackage{amsmath}
\usepackage{amsfonts}
\begin{document}
\bibliographystyle{plain}
\title{On the decomposition of Global Conformal Invariants I}
\author{Spyros Alexakis}
\maketitle
\newtheorem{proposition}{Proposition}
\newtheorem{theorem}{Theorem}
\newtheorem{lemma}{Lemma}
\newtheorem{conjecture}{Conjecture}
\newtheorem{observation}{Observation}
\newtheorem{formulation}{Formulation}
\newtheorem{definition}{Definition}
\newtheorem{corollary}{Corollary}

\begin{abstract}
This is the first of two papers where we address and partially
confirm a conjecture of Deser and Schwimmer, originally postulated
in high energy physics.
 The objects of study are
 scalar Riemannian quantities constructed out of the curvature
 and its covariant derivatives, whose
 integrals over compact manifolds are invariant under
 conformal changes of the underlying metric. Our main
 conclusion is that each such quantity that locally depends
 only on the curvature tensor (without covariant derivatives)
can be written as a linear combination of the
Chern-Gauss-Bonnet integrand and a scalar conformal
 invariant.
\end{abstract}

\section{Introduction}
\subsection{Outline of the problem.}
\par Consider any Riemannian manifold
 $(M^n,g^n)$. The basic local objects that describe the geometry
 of the metric $g^n$ are the curvature tensor $R_{ijkl}$
and the Levi-Civita connection $\nabla_{g^n}$. We are interested
in {\it intrinsic scalar quantities} $P(g^n)$. These scalar
quantities, as defined by Weyl (see also \cite{e:ntrm} and
\cite{bfg:spdcn}), are polynomials in the components of the
tensors $R_{ijkl},\dots , \nabla^m_{r_1\dots r_m}R_{ijkl},\dots$
and $g^{ij}$, with two basic features: The values of these
polynomials must be invariant under changes of the coordinate
system (or isometries), and there must also be a number $K$ so
that under the re-scaling $g^n\longrightarrow t^2g^n$
($t\in\mathbb{R}_{+}$), we have $P(t^2 {g}^n)=t^K P(g^n)$. We then
say that $P(g^n)$ is a {\it scalar Riemannian invariant} of weight
$K$.
\newline

\par  It is a classical result, implied in Weyl's work
\cite{w:cg},  that any such Riemannian invariant $P(g^n)$ of
weight $K$ can be written as
 a linear combination

\begin{equation}
\label{lincomb}
P(g^n)=\Sigma_{l\in L} a_l C^l(g^n)
\end{equation}
 of complete contractions $C^l(g^n)$ in
 the form:

\begin{equation}
\label{contraction} contr(\nabla^{m_1}_{r_1\dots r_{m_1}}
R_{i_1j_1k_1l_1}\otimes\dots
\otimes \nabla^{m_r}_{t_1\dots t_{m_r}}R_{i_rj_rk_rl_r})
\end{equation}
for which $C^l(t^2g^n)=t^KC^l(g^n)$.

\par This notion of {\it intrinsic} extends to vector fields. We
define an intrinsic vector field $T_a(g^n)$ ($a$ is the free
index) of weight $K$ to be a polynomial in the components of
 the
tensors $R_{ijkl},\dots ,\nabla^m_{r_1\dots r_m}R_{ijkl},
\dots$ and $g^{ij}$,
with the property that under changes of coordinates (isometries)
that map the coordinate functions $x^1,\dots ,x^n$ to the
coordinate functions $y^1,\dots ,y^n$, $T_a(g^n)$ must satisfy the transformation law:

$$T'_{\alpha}(g^n)=
T_i(g^n)
\frac{\partial x^i}{\partial y^\alpha}$$
 where $T'_\alpha$ stands for
 the vector field expressed in the new coordinate system.
 Moreover, we say that $T_a(g^n)$ has
 weight $K$ if $T_a(t^2g^n)=t^{K+1}T_a(g^n)$.

\par By Weyl's work, we have that an intrinsic vector field of
weight $K$ must be a linear combination of partial contractions,
with one free index, in the form:

\begin{equation}
\label{pcontraction} pcontr(
\nabla^{m_1}_{r_1\dots r_{m_1}}R_{ijkl}\otimes\dots \otimes
\nabla^{m_r}_{t_1\dots t_{m_r}}R_{i_rj_rk_rl_r})
\end{equation}

 \par We recall that under general conformal re-scalings
 $\hat{g}^n=e^{2\phi(x)}g^n$ the volume form
 re-scales by the formula $dV_{\hat{g}^n}=e^{n\phi(x)}
dV_{g^n}$, in particular for any constant $t$ we have
$dV_{t^2g^n}=t^ndV_{g^n}$. Thus, for any scalar Riemannian invariant
$P(g^n)$ of weight $-n$ we have that
${\int}_{M^n}P(g^n)dV_{g^n}$ is scale-invariant for all
 compact and orientable  manifolds $M^n$.
\newline

\par The problem we address in this paper and in
\cite{a:dgciII} is to find all the Riemannian invariants $P(g^n)$
of weight $-n$ for which the integral:

\begin{equation}
\label{defprop}
  \int_{M^n}P(g^n)dV_{g^n}
\end{equation}
 is invariant under conformal re-scalings of the metric $g^n$
 on any $M^n$ compact without boundary.
\newline

\par In other words, we are requiring that  for every real-valued function
$\phi(x)\in C^{\infty}(M^n)$ we must have that  for
$\hat{g}^n(x)= e^{2\phi(x)}g^n$:

\begin{equation}
\label{expl1}
\int_{M^n}P(\hat{g}^n)dV_{\hat{g}^n}
=\int_{M^n}P(g^n)  dV_{g^n}
\end{equation}

\par This question was originally raised by Deser and
 Schwimmer in \cite{ds:gccaad}
(see also \cite{hs:hwa} and \cite{be:clwid8}) in the context of understanding ``conformal anomalies''.
 On the other hand, an answer to this question would also
 shed light on the structure of $Q$-curvature in high
 dimensions. The problem, as posed in
\cite{ds:gccaad}, is the following:

\begin{conjecture}[Deser-Schwimmer]
\label{conj} Consider a Riemannian scalar $P(g^n)$ of weight $-n$,
for some even $n$. Suppose that for any compact manifold
$(M^n,g^n)$ the quantity

\begin{equation}
\label{donne} {\int}_{M^n}P(g^n)dV_{g^n}
\end{equation}
is invariant under any conformal change of metric
$\hat{g}^n(x)=
 e^{2\phi (x)} g^n(x)$. Then $P(g^n)$ must be a linear combination of
 three``obvious candidates'', namely:

\begin{equation}
\label{post} P(g^n)=W(g^n)+div_iT_i(g^n)+c\cdot
\operatorname{Pfaff}(R_{ijkl})
\end{equation}

\begin{enumerate}

\item{$W(g^n)$ is a scalar conformal invariant of weight
$-n$, ie
it satisfies $W(e^{2\phi(x)}g^n)(x)=e^{-n\phi(x)}W(g^n)(x)$
for every $\phi\in C^\infty(M^n)$ and every $x\in M^n$.}

\item{$T_i(g^n)$ is a Riemannian vector field of weight
$-n+1$.
(Since for any compact $M^n$ we have $\int_{M^n}div_i
T_i(g^n)dV_{g^n}=0$.)}

\item{$\operatorname{Pfaff}(R_{ijkl})$ stands for the Pfaffian of
the curvature $R_{ijkl}$. (Since for any compact Riemannian
$(M^n,g^n)$ $\int_{M^n}
\operatorname{Pfaff}(R_{ijkl})dV_{g^n}=
\frac{2^n\pi^{\frac{n}{2}}(\frac{n}{2}-1)!}{2(n-1)!}\chi(M^n)$.)}
\end{enumerate}
\end{conjecture}

\par In this paper and in \cite{a:dgciII} we show:

\begin{theorem}
\label{weakt2} Conjecture \ref{conj} is true, in the following
restricted version:

Let us suppose that (\ref{donne}) holds, and additionally that
$P(g^n)$ locally depends only on the curvature tensor
$R_{ijkl}$ and not its covariant derivatives
$\nabla^mR_{ijkl}$
(meaning that $P(g^n)$ is a linear combination of contractions in
the form (\ref{contraction})
 with $m_1=\dots =m_r=0$).
Then, there exists a a scalar conformal invariant $W(g^n)$ of
 weight $-n$ that locally depends only on the Weyl tensor,
and also a constant $c$ so that:

\begin{equation}
\label{post2} S(g^n)=W(g^n)+c\cdot \operatorname{Pfaff}
(R_{ijkl})
\end{equation}
where $\operatorname{Pfaff}(R_{ijkl})$ stands for the Pfaffian of
the curvature $R_{ijkl}$.
\end{theorem}

\par The proof of the above will shed light both on global conformal invariants
that locally depend only on the curvature tensor (and not its covariant
derivatives), but also on the structure of the Pfaffian of the
 curvature tensor.
\newline

\subsection{Geometric Applications of the Deser-Schwimmer Conjecture:
$Q$-curvature and re-normalized volume.}

\par $Q$-curvature is a Riemannian scalar quantity introduced by Branson for
each even dimension $n$ (see \cite{b:fd}). In dimension 2,
$Q^2(g^2)=R(g^2)$ (the scalar curvature), while in dimension $4$
its structure is well-understood and has been extensively studied.
Its fundamental property is that $Q^n(g^n)$ has weight $-n$ in
dimension $n$ and that the integral $\int_{M^n}Q^n(g^n)dV_{g^n}$
over compact
 manifolds $M^n$ is invariant under conformal charges of the underlying
metric $g^n$. Thus, if one proves Conjecture \ref{conj} in full
 strength, one would obtain that $Q^n(g^n)$ can be decomposed
 as in the right hand side of (\ref{post}), in fact with
 $c\ne 0$.

\par This fact is all the more interesting due to the
{\it nice} transformation law of $Q$-curvature under conformal changes $\hat{g}^n=e^{2\phi(x)}g^n$. One then has
 that $e^{n\phi(x)}Q^n(\hat{g}^n)=Q^n(g^n)+P^n_{g^n}(\phi)$,
where $P^n_{g^n}(\phi)$ is a {\it conformally invariant differential operator}, originally constructed in
\cite{gjms:cipl}. Thus, prescribing the $Q$-curvature
can be informally interpreted as prescribing a modified
 version of the Chern-Gauss-Bonnet integrand
$\operatorname{Pfaff}(R_{ijkl})$. This {\it modified
 Pfaffian} enjoys  a nice transformation law under conformal re-scalings,
 rather than the messy transformation that governs
$\operatorname{Pfaff}(R_{ijkl})$.

\par This understanding of the structure of
 $Q$-curvature in any even dimension raises the
 question whether the strong results of Chang, Yang, Gursky, Qing et
 al in dimension $4$ (see for example \cite{cgy:ematcg},
\cite{cqy:tccem},
 \cite{g:pecidowasep}), have analogues in higher dimensions.
  Moreover, a proof of Conjecture \ref{conj} in full strength
 will
  lead to a better understanding of the notion of
{\it re-normalized volume}
  for conformally compact Einstein manifolds.
\newline

\par Conformally compact Einstein manifolds
 have been the focus of much
research in recent years, see ~\cite{cqy:tccem}, \cite{gz:smcg},
\cite{q:rccem}, \cite{wy:cbacc}, to name just a few. What follows
is a very brief discussion of the objects of study, largely
 reproduced
 from \cite{gz:smcg}.

\par We consider manifolds with boundary,
$(X^{n+1}, g^{n+1})$,
$\partial X^{n+1}=M^n$, where the boundary $M^n$ carries a
conformal structure $[h^n]$. We consider a defining function $x$
for $\partial X^{n+1}$ in $X$:

$$x|_{\mathaccent23{X}}>0, \text{    } x|_{\partial X}=0, \text{   }
dx|_{\partial X}\ne 0$$
\par We then say that $g^{n+1}$ is a {\it conformally compact} metric on
$X^{n+1}$ with {\it conformal infinity} $[h^n]$, if we can find a
smooth
 metric $\overline{g}^{n+1}$ on $\overline{X}^{n+1}$ so that in
$\mathaccent23{X}^{n+1}$:

$$g^{n+1}=\frac{\overline{g}^{n+1}}{x^2}, \text{  }
\overline{g}^{n+1}|_{\partial X^{n+1}}\in [h^n]$$ A conformally
compact metric is {\it asymptotically
 hyperbolic}, in the sense that its sectional curvatures approach $-1$ as $x$ approaches 0.  We notice
that since we can pick different defining functions, the metric
$g^{n+1}$ in the interior $\mathaccent23{X}^{n+1}$ defines a whole
conformal class on the boundary. In the rest of this discussion,
we will be considering conformally compact manifolds $(X^{n+1},
g^{n+1})$ which in addition are Einstein.
\newline

\par Conformally compact Einstein manifolds are studied as
models for the Anti-de-Sitter/Conformal Field Theory
(AdS-CFT) correspondence in string theory. In order to
compute the partition function for the conformal field
 theory in the supergravity approximation, one must evaluate
 the Einstein action for the
metric $g^{n+1}$, which in the case at hand is proportional
to the volume of $(X^{n+1},g^{n+1})$. Since this volume is
 clearly infinite ($g^{n+1}$ is asymptotically hyperbolic)
one regularizes it through re-normalization, thus introducing
 the {\it re-normalized volume}. We briefly
 discuss this re-normalization procedure and its relation to
 $Q$-curvature below. For a more detailed discussion we refer
 the reader to \cite{g:varccem}, \cite{gw:casoacc},
\cite{w:adssh} and the references therein.

 It is known that each choice of metric $h\in [h^n]$ on the boundary $M^n$
 uniquely determines a defining function $x$ in a collar neighborhood of $\partial X^{n+1}$ in $X^{n+1}$, say
 $\partial X^{n+1}\times [0,\epsilon]$, so that $g^{n+1}$ takes the form:

\begin{equation} \label{aggouria} g^{n+1}=x^{-2}(dx^2
+h_x),\text{ } h_0=h,
\end{equation}
 where $h_x$ is a 1-parameter
family of metrics on $\partial X^{n+1}$. We then consider the
volume of the region $R_\epsilon=\{ x>\epsilon\}$ in $(X^{n+1},
g^{n+1})$, expanded out in powers of $\epsilon$, and let
$\epsilon\rightarrow 0$.
 Given that $g^{n+1}$ is Einstein, we have that if $n$ is odd, then:

\begin{equation}
\label{jp}
vol_{g^{n+1}}(\{ x>\epsilon\})=c_0\epsilon^{-n}+c_2\epsilon^{-n+2}+
\dots +c_{n-1}\epsilon^{-1}+V+o(1)
\end{equation}
whereas if $n$ is even, then:

\begin{equation}
\label{jp2}
vol_{g^{n+1}}(\{ x>\epsilon\})=c_0\epsilon^{-n}+c_2\epsilon^{-n+2}+
\dots +c_{n-1}\epsilon^{-2}+Llog(\frac{1}{\epsilon})+V+o(1)
\end{equation}

\par Moreover, if $n$ is odd and since $g^{n+1}$ is Einstein,
Graham and Zworski showed in \cite{gz:smcg}
that $V$ is independent of the choice of metric $h^n$ in the conformal
  class $[h^n]$. (Recall that this choice was used in order to write out
$g^{n+1}$ in the form (\ref{aggouria}), and hence also in defining
the region $R_\epsilon$, therefore $V$ depends apriori on the
choice $h^n\in [h^n]$).
 For $n$ odd, $V$ is called the {\it re-normalized volume}
of $(X^{n+1}, g^{n+1})$.

\par For $n$ even, one has that $V$ is {\it not} independent
 of the
choice of metric $h^n$ in the conformal class $[h^n]$. In this
case it is the quantity $L$ that demonstrates this invariance.
This quantity $L$ represents the failure of defining the
re-normalized volume independently of the defining function $x$.
It is therefore called the ``conformal anomaly'' in the physics
literature. Moreover, Graham-Zworski  have shown that
$L=c_n\cdot\int_{M^n}Q(h^n)dV_{h^n}$, where $h^n$ is an
 arbitrary metric in
the conformal class $[h^n]$. Hence, a proof of Conjecture
\ref{conj}
would immediately imply that $L$ can be written out as:

\begin{equation}
\label{indoi}
L=\int_{M^n} W(h^n)dV_{h^n}+(Const)\cdot \chi (M^n)
\end{equation}
where $W(h^n)$ is a scalar conformal invariant of weight $-n$ and
$M^n=\partial X^{n+1}$,while $\chi(M^n)$ stands for the Euler
 characteristic of  $M^n$ and $(Const)\ne 0$.
\newline

\par Another significant result has recently been obtained by
Chang, Qing and Yang, \cite{cqy:pc}, relating the
re-normalized
volume $V$ with the $Q$-curvature of $g^{n+1}$ and
hence with the Euler characteristic of the  manifold $X^{n+1}$.
They show
 that if Conjecture 1 is true, then for $n$ odd one
can express the re-normalized volume of
$(X^{n+1},g^{n+1})$ via the $Q$-curvature:

\begin{equation}
\label{chandqiya}
\begin{split}
&R.V.[(X^{n+1},g^{n+1})]=(Const)_{n+1}\cdot \int_{X^{n+1}}
Q^{n+1}(g^{n+1})dV_{g^{n+1}}=
\\&\int_{X^{n+1}}W(g^{n+1})
dV_{g^{n+1}}+(const)_{n+1}\cdot \chi (X^{n+1})
\end{split}
\end{equation}
where $(Const)_{n+1}, (const)_{n+1}$ are nonzero dimensional
constants and $W(g^{n+1})$ is a scalar conformal
invariant of weight $-n-1$. Here the left hand side
 stands for the re-normalized volume of the manifold
$(X^{n+1},g^{n+1})$. Hence, it follows that the re-normalized
volume explicitly depends on the topology of $X^{n+1}$, via its
Euler characteristic. A result related to (\ref{chandqiya}) has
been independently established (by an entirely different
 method) by Albin in \cite{a:rcipem}.

\par This identity raises the question of whether one can adapt
 the powerful techniques developed for the study of
$Q$-curvature to the study of conformally compact Einstein
manifolds. Strong results have already been obtained in
dimension $4$, see \cite{cqy:tccem}. For higher dimensions one
 might try to extend the work of Brendle \cite{b:gechofcg}
 to this setting. Another question would be whether one can
obtain expressions analogous to (\ref{indoi}) and
(\ref{chandqiya}) for the re-normalized areas and conformal
anomalies of submanifolds, defined by Graham and
 Witten in \cite{gw:casoacc}.

\subsection{Outline of the paper.}
\label{outpaper}

\par Our theorem is a structure result for $P(g^n)$. We use the
``global'' conformal invariance under integration of $P(g^n)$ to
derive information on its algebraic expression.

In this paper we introduce the main tool that will show Theorem
\ref{weakt2}, the so-called super divergence formula. For each
$P(g^n)$ that satisfies (\ref{expl1}), we define an operator
$I_{g^n}(\phi)$ that measures the ``non-conformally invariant
part'' of $P(g^n)$ (see (\ref{defI}) below). We then use the
property (\ref{explI}) of $I_{g^n}(\phi)$ to derive the super
divergence formula for $I_{g^n}(\phi)$ (in fact, for a polarized
version of $I_{g^n}(\phi)$). This formula, which in our opinion is
also of independent interest, provides us with an understanding of
the structure of $I_{g^n}(\phi)$. In the sequel to this paper,
\cite{a:dgciII}, we will use the super
 divergence formula to
derive information on the structure of $P(g^n)$ and prove Theorem
\ref{weakt2}.
\newline

\par The super divergence formula is proven in a number of steps.
A more primitive version is the ``simple divergence formula''
in section \ref{simplediv}.
This is then refined three times in section \ref{threerefs}
and we obtain the super
divergence formula in subsection \ref{superdiv}. The only background material needed for all
this work is a slight extension of Theorem B.4 in
\cite{beg:itccg},
which itself is a generalization of a classical theorem of Weyl
in \cite{w:cg}. This extension is discussed in section \ref{transdim}.
Roughly, Theorem B.4 in \cite{beg:itccg} and our
Theorem \ref{UMxi} below assert that a linear identity involving
complete contractions which holds for all values we can give to the
tensors in those contractions, must then also hold {\it formally}.

\section{Background material.}
\label{backmat}

\subsection{Definitions and Identities.}
\label{def&id}

We note that whenever we refer to a manifold $M^n$, we will be
assuming it to be {\it compact} and {\it orientable}. Moreover,
$n$ will be a {\it fixed, even} dimension throughout this paper. We begin by
recalling a few definitions and formulas.

\begin{definition}
\label{defcontr}
In this paper, we will be dealing with complete contractions
 and their linear combinations. Any complete contraction:

$$C=contr((A^1)_{i_1\dots i_s}\otimes\dots \otimes
(A^t)_{j_1\dots j_q})$$
will be seen as a formal expression.
Each factor $(A^l)_{i_1\dots i_s}$ is an ordered set of
 slots. Given the factors $(A^1)_{i_1\dots i_s},\dots ,
(A^t)_{j_1\dots j_q}$, a complete contraction is then
 seen as a set of pairs of slots $(a_1,b_1),\dots ,(a_w,b_w)$,
with the following properties: if $k\ne l$, $\{a_l,b_l\}
\bigcap \{a_k,b_k\}=\emptyset$, $a_k\ne b_k$,
${\bigcup}_{i=1}^w \{a_i,b_i\}=\{i_1,\dots ,j_q\}$.
Each pair corresponds to a particular contraction.

\par Two complete contractions:
$$contr((A^1)_{i_1\dots i_s}\otimes\dots \otimes
(A^t)_{j_1\dots j_w})$$

and
$$contr((B^1)_{f_1\dots f_q}\otimes\dots \otimes
(B^{t'})_{v_1\dots v_z})$$

will be identical if $t=t'$, $(A^l)=(B^l)$ and if the
${\mu}^{th}$ index in $A^l$ contracts against the ${\nu}^{th}$
 index in $A^r$, then the ${\mu}^{th}$ index in $B^l$ contracts  against the ${\nu}^{th}$ index in $B^r$. For any
 complete contraction, we define its length to stand for the
 number of its factors.

\par We can also consider linear combinations of complete contractions:

$${\Sigma}_{l\in L} a_l (C_1)^l$$
and

$${\Sigma}_{r\in R} b_r (C_2)^r$$

Two linear combinations as above are considered identical if
 $R=L$ and $a_l=b_l$ and $(C_1)^l=(C_2)^l$.
A linear combination of complete contractions as above is identically
 zero if for every $l\in L$ we have that $a_l=0$.

\par For any complete contraction $C$, we will say a
 factor $(A)_{r_1\dots r_{s_l}}$ has
 an internal contraction if two indices in
$(A)_{r_1\dots r_{s_l}}$
are contracting between themselves.

\par All the above definitions extend to partial contractions
and their linear combinations.
\end{definition}

\par We also introduce two
 language conventions: For any linear combination of complete
 contractions $\Sigma_{l\in L} a_l C^l$, when we speak of a
 {\it sublinear combination}, we will mean some linear combination
 $\Sigma_{l\in L'} a_l C^l$ where $L'\subset L$. Also, when we say
 that an identity between linear combinations of complete
 contractions:

\begin{equation}
\label{chuda} \Sigma_{r\in R} a_r C^r=\Sigma_{t\in T} a_t C^t
\end{equation}
 {\it holds modulo complete contractions of length $\ge
\lambda$}, we will mean that we have an identity:
\begin{equation}
\label{chuda2} \Sigma_{r\in R} a_r C^r=\Sigma_{t\in T} a_t
C^t+\Sigma_{u\in U} a_u C^u
\end{equation}
where each $C^u$ has at least $\lambda$ factors.
\newline

\begin{definition}
\label{symmetrization}
 Now, for each tensor $T_{ab\dots d}$ and each subset
$\{d,e,\dots f\}\subset \{a,b,\dots ,d\}$, we define the symmetrization
 of the tensor $T_{ab\dots d}$ over the slots $d,e,\dots ,f$:

\par Let $\Pi$ stand for the set of permutations of the ordered set
$\{d,e,\dots ,f\}$. For each $\pi\in \Pi$, we define $\pi
T_{ab\dots f}$ to stand for the tensor that arises from
$T_{ab\dots f}$ by permuting the slots $d,e,\dots ,f$ according to
the permutation $\pi$.  We then define the symmetrization of the
tensor $T_{ab\dots d}$ over the slots $d,e,\dots ,f$ to be:

$$\Sigma_{\pi\in \Pi}\frac{1}{|\Pi|} \cdot \pi T_{ab\dots d}$$

\par If $\{d,e,\dots f\}= \{a,b,\dots ,d\}$, we will denote that
symmetrization by $T_{(ab\dots d)}$.
\end{definition}

\par We recall a few basic facts from Riemannian geometry.
 Consider any Riemannian manifold $(M^n,g^n)$ and any
$x_0\in M^n$. We pick any coordinate system $x^1,\dots ,x^n$
 and denote by $X_i$
the coordinate vector fields, ie the vector fields
$\frac{\partial}{\partial x^i}$. We will write $\nabla_i$
instead of $\nabla_{X_i}$.

\par The curvature tensor $R_{ijkl}$ of $g^n$ is given
 by the formula:

\begin{equation}
\label{curvature}
[\nabla_i\nabla_j-\nabla_j\nabla_i]X_k=R_{ijkl}X^l
\end{equation}

\par In a coordinate system, we can also express it in terms
 of the Christoffel symbols:

\begin{equation}
\label{kormaki}
R^l_{ijk}={\partial}_j{\Gamma}^l_{ik}-{\partial}_k
{\Gamma}^l_{ij}+
{\Sigma}_m({\Gamma}^m_{ik}{\Gamma}^l_{mj}-
{\Gamma}^m_{ik}{\Gamma}^l_{mk})
\end{equation}

\par Moreover, the Ricci tensor $Ric_{ik}$ is then:

\begin{equation}
\label{ricci}
Ric_{ik}=R_{ijkl}g^{jl}
\end{equation}

We recall  the two Bianchi identities:

\begin{equation}
\label{1bianchi} R_{ABCD}+R_{CABD}+R_{BCAD}=0
\end{equation}

\begin{equation}
\label{2bianchi} {\nabla}_A R_{BCDE}+{\nabla}_C
R_{ABDE}+{\nabla}_B R_{CADE}=0
\end{equation}

\par We also recall how the basic geometric objects transform under the conformal
change $\hat{g}^n(x)= e^{2\phi(x)}g^n(x)$. These formulas come
from \cite{e:ncg}.

\begin{equation}
\label{curvtrans}
\begin{split}
&R_{ijkl}^{\hat{g}^n}=e^{2\phi (x)}[R^{g^n}_{ijkl}+
\nabla_{il}{\phi}
g_{jk}+\nabla_{jk}{\phi}g_{il}-\nabla_{ik}{\phi}g_{jl}-\nabla_{jl}{\phi}g_{ik}
+\nabla_i{\phi}\nabla_k{\phi}g_{jl}+\nabla_j{\phi}\nabla_l{\phi}g_{ik}
\\&-\nabla_i{\phi}
\nabla_l{\phi} g_{jk} -\nabla_j{\phi}\nabla_k{\phi}g_{il}
+|\nabla\phi|^2g_{il}g_{jk}- |\nabla\phi|^2g_{ik}g_{lj}]
\end{split}
\end{equation}

\begin{equation}
\label{ricci}
{Ric}^{\hat{g}^n}_{\alpha\beta}={Ric}^{g^n}_{\alpha\beta}+
(2-n){\nabla}^2_{\alpha\beta}\phi
- {\Delta}\phi
g^n_{\alpha\beta}+(n-2)({\phi}_{\alpha}{\phi}_{\beta}-{\phi}^k{\phi}_k
g^n_{\alpha\beta})
\end{equation}

While the transformation law for the Levi-Civita
connection is:

\begin{equation}
\label{levicivita} {\nabla}^{\hat{g}^n}_k {\eta}_l=
\nabla_k^{g^n}{\eta}_l -\nabla_k{\phi} {\eta}_l -\nabla_l{\phi}
{\eta}_k +\nabla^s{\phi} {\eta}_s g^n_{kl}
\end{equation}

\par We now focus on complete contractions $C(g^n)$ in the
 form (\ref{curvature}). We still think of these objects as formal expressions, but also as functions of the metric
$g^n$. Thus, for complete contractions in the
 form (\ref{contraction}), contracting two lower indices
${}_a,{}_b$ will mean that we multiply by $g^{ab}$ and then sum over $a,b$. We have that under the rescaling
$\hat{g}^n=t^2g^n$
the tensors $\nabla^mR_{ijkl}$ and $(g^n)^{ij}$ transform by
$\nabla^m_{r_1\dots r_m}R_{ijkl}(t^2g^n)=t^2
\nabla^m_{r_1\dots r_m}R_{ijkl}(g^n)$,
$(g^n)^{ij}(t^2 g^n)=t^{-2}(g^n)^{ij}(g^n)$. (We will sometimes write $\nabla^mR_{ijkl}$ instead of
$\nabla^m_{r_1\dots r_m}R_{ijkl}$, for brevity). Thus, for
 each $C(g^n)$ in the form (\ref{contraction}), if we define
$K=-\Sigma_{i=1}^r (m_i+2)$, we will have that $C(t^2g^n)=
t^K C(g^n)$. We define $K$ to be the {\it weight} of this complete contraction.
\newline

\par For future reference, we will consider more general
 complete contractions defined on manifolds $(M^n,g^n)$
and define their weight.

\begin{definition}
\label{weight}
 We consider any complete contraction
$C_{g^n}(V^1,\dots ,V^x)$ in the form:

\begin{equation}
\label{triangul}
contr(\nabla^{m_1}R_{ijkl}\otimes\dots\otimes\nabla^{m_r}
R_{ijkl}\otimes V^1_{a_1\dots a_{f_1}}\otimes\dots\otimes
V^x_{b_1\dots b_{f_x}})
\end{equation}
defined for any $x_0\in M^n$.
Here the tensors $V^y_{a_1\dots a_{f_y}}$ are auxiliary
 tensors (all of whose indices are lowered) that have a
 scaling property under re-scalings of the metric:
$V^y_{a_1\dots a_{f_y}}(t^2g^n)=t^{C_y}
V^y_{a_1\dots a_{f_y}}(t^2g^n)$.
(An example for a tensor $V^y_{a_1\dots a_{f_y}}$ would be the
 $y^{th}$ iterated covariant
 derivative of a function $\psi$, in which case $C_y=0$).
Moreover, all the tensors here are over
$TM^n|_{x_0}$. The particular contractions of any two lower
indices will be with respect to the quadratic form
$(g^n)^{ij}(x_0)$.

We then define the weight of such a complete contraction to be
$W=-\Sigma_{i=1}^r(m_i+2)-\Sigma_{i=1}^x(f_i-C_y)$.
As for the previous case, we then have that:
$C_{t^2g^n}(V^1,\dots ,V^x)=t^WC_{g^n}(V^1,\dots ,V^x)$.
\end{definition}

In this whole paper, when we write a complete contractions
 and include the metric $g^n$ in the notation, we will imply
 that the contraction is defined on manifolds (and possibly
also depending on additional auxiliary objects, for example
  scalar functions) and will
 have a weight, as defined above. Unless otherwise stated, all
 complete contractions will have weight $-n$.

\subsection{The operator $I_{g^n}(\phi)$ and its polarizations.}
\label{I&P}

\par For this paper and in \cite{a:dgciII}, we will be
 considering $P(g^n)$ as
a linear combination in the form:

\begin{equation}
\label{pointofref}P(g^n)=\Sigma_{l\in L} a_l C^l(g^n)
\end{equation}
where each $C^l(g^n)$ is in the form (\ref{contraction}) and
 has weight $-n$. We assume that $P(g^n)$ satisfies
(\ref{expl1}).

\par We define a differential operator, which will depend
both on the metric $g^n$ and and auxiliary
$\phi\in C^\infty (M^n)$:

\begin{equation}
\label{defI}
I_{g^n}(\phi)(x)=e^{n\phi(x)}P(e^{2\phi(x)}g^n)(x)-P(g^n)(x)
\end{equation}

  We then have by (\ref{expl1}) that:

\begin{equation}
\label{explI} \int_{M^n}I_{g^n}(\phi)dVg^n=0
\end{equation}
 for every compact manifold $(M^n,g^n)$ and any function
$\phi\in C^{\infty}(M^n)$.
 Then, using the transformation laws
(\ref{curvtrans}) and (\ref{levicivita}) we see that
$I_{g^n}(\phi)$ is a differential operator acting on
the function $\phi$. In particular, we can pick any $A>0$
functions ${\psi}_1(x),\dots ,{\psi}_A(x)$ , and choose:
$$\phi (x)= {\Sigma}_{l=1}^A {\psi}_l(x)$$
 Hence, we have a differential operator $I_{g^n}({\psi}_1,\dots ,
{\psi}_A)(x)$, so that, by (\ref{explI}):
$$ {\int}_{M^n}I_{g^n}({\psi}_1,\dots ,{\psi}_A)dV_{g^n}=0$$
 for any $(M^n,g^n)$, $M^n$ compact and any
functions ${\psi}_1(x),\dots ,{\psi}_A(x)\in C^{\infty}(M^n)$.

  Now, for any given functions ${\psi}_1(x),\dots ,{\psi}_A(x)$, we
can consider re-scalings: $$\lambda_1 \psi_1 (x),\dots ,\lambda_A
\psi_{A}(x)$$

  Hence, as above we will have the equation:
\begin{equation}
\label{ilambda}
{\int}_{M^n}I_{g^n}(\lambda_1\psi_1,\dots
,\lambda_A\psi_A)dV_{g^n}=0
\end{equation}

  We can
then see ${\int}_{M^n}I_{g^n}({\lambda}_1{\psi}_1, \dots ,
{{\lambda}_A}{\psi}_A)dV_{g^n}$ as a
polynomial in the factors ${\lambda}_1,\dots,{\lambda}_A$.
 Call this polynomial ${\Pi}({\lambda}_1,\dots ,{\lambda}_A)$.

 \par But then relation (\ref{ilambda}) gives us that this polynomial $\Pi$ is
identically zero. Hence, each coefficient of each monomial in the variables
  ${\lambda}_1 ,\dots ,{\lambda}_A$ must be zero. We want to
pick out a particular such monomial. Pick out any integer $1\le Z
\le A$.

 \par Then in $I_{g^n}({\lambda}_1{\psi}_1, \dots ,
{\lambda}_A{\psi}_A)$ (seen as a multi-variable polynomial in
${\lambda}_1 ,\dots ,{\lambda}_A$) consider the coefficient of the
monomial ${\lambda}_1 \cdot\dots \cdot{\lambda}_Z$. This will be a
differential operator in the functions ${\psi}_1,\dots ,
{\psi}_Z$, which we will denote by $I^Z_{g^n} ({\psi}_1, \dots
,{\psi}_Z)$. By elementary properties of polynomials and by the
definition of $I_{g^n}(\phi)$ in (\ref{defI}) we have:

\begin{equation}
\begin{split}
\label{derlami} &I^Z_{g^n} ({\psi}_1, \dots ,{\psi}_Z)=\\&
{\partial}_{{\lambda}_1}\dots {\partial}_{{\lambda}_Z}
[e^{n({\lambda}_1 {\psi}_1+\dots +{\lambda}_Z
{\psi}_Z)}P(e^{2({\lambda}_1 {\psi}_1+\dots +{\lambda}_Z
{\psi}_Z)}g^n)]|_{{\lambda}_1=0,\dots ,{\lambda}_Z=0}
\end{split}
\end{equation}

\par The precise form of $I^Z_{g^n} ({\psi}_1, \dots ,{\psi}_Z)$, given
$P(g^n)$, can be calculated using the transformation laws in the
previous section. We will be doing this in \cite{a:dgciII}.
 For the time being, just note that by (\ref{ilambda}) we
 have the equation:

\begin{lemma}
\begin{equation}
\label{lI} \int_{M^n}I^Z_{g^n} ({\psi}_1, \dots
,{\psi}_Z)dV_{g^n} =0
\end{equation}
for every compact $(M^n,g^n)$ and any ${\psi}_1,\dots ,
{\psi}_Z\in C^{\infty}(M^n)$.
\end{lemma}
 {\it Proof:} Straightforwardly from relation (\ref{ilambda}) and the
 equation (\ref{derlami}). $\Box$
\newline

 \par From all the above, it is straightforward to see that
$I^Z_{g^n} ({\psi}_1, \dots ,{\psi}_Z)$ is a linear combination of
complete contractions of weight $-n$ in the form:

\begin{equation}
\label{linicontrpsi}
\begin{split}
&contr({\nabla}_{r_1\dots r_{m_1}}^{m_1}R_{i_1j_1k_1l_1}\otimes
\dots \otimes {\nabla}_{v_1\dots
v_{m_s}}^{m_s}R_{i_sj_sk_sl_s}\otimes
\\&{\nabla}^{{\nu}_1}_ {{\chi}_1\dots
{\chi}_{{\nu}_1}}{\psi}_1 \otimes\dots \otimes
{\nabla}^{{\nu}_Z}_{{\omega}_1\dots {\omega}_{{\nu}_Z}} {\psi}_Z)
\end{split}
\end{equation}

\par For the rest of this paper, we will only be using the
 fact that $I^Z_{g^n}(\psi_1,\dots ,\psi_Z)$ satisfies
(\ref{lI}) and that it is a linear combination of complete
 contractions in the form (\ref{linicontrpsi}).

\section{The Trans-Dimensional Isomorphisms.}
\label{transdim}

\par The aim of this section is to show that there is a natural isomorphism
of linear combinations of  weight $-n$ complete contractions
 in the form
(\ref{linicontrpsi}), between dimensions $N$ and
$n$, if $N\ge n$. In order to make this statement precise and
to provide a proof, we will recall some terminology and facts
from the appendices in \cite{beg:itccg}. The main
``known fact'' that we will be using is Theorem \ref{UMxi} in
 the next subsection. This theorem is a slight generalization
 of Theorem B.4 in \cite{beg:itccg}, and it can be
 proven using the same ideas. The appendices in
\cite{beg:itccg} generalize
 classical theorems that can be found in \cite{w:cg}.

\subsection{Known facts}
\label{knowfact}

\par The appendices of \cite{beg:itccg} deal with identities
involving linear combinations of complete contractions. The main
assertion there is that under certain hypotheses, when a linear
identity involving complete contractions holds ``by
substitution'', it must then also hold ``formally''. We will be
explaining these notions in this subsection. For more details, we refer the reader to \cite{beg:itccg}.
\newline

\par We introduce the ``building blocks'' of our complete
contractions. Firstly, we consider  symmetric tensors. Let us
consider a family of sets of symmetric tensors
$\{T^\alpha=\{ T^\alpha_0,
T^\alpha_i,\dots ,T^\alpha_{i_1\dots i_s},\dots
\}\}_{\alpha\in A}$
($T^\alpha_0$
is just a scalar, ie a tensor of rank zero), defined over the
vector space $\mathbb{R}^n$. Here each $\alpha\in A$ is not a free
index of the tensors $T^\alpha_{i_1\dots i_s}$. It just serves to
distinguish the tensors $T^{\alpha_1}_{i_1\dots i_s}$ and
$T^{\alpha_2}_{i_1\dots i_s}$ when $\alpha_1\ne \alpha_2$.
\newline

\par Our second building block will be a list
 of tensors that resemble the covariant derivatives of the curvature tensor:

\begin{definition}
\label{lincurv} A set of linearized curvature tensors is
 defined to be a list  of tensors
$R=\{R_{ijkl},\dots ,R_{f_1\dots f_s,ijkl},\dots \}$ defined
 over $\mathbb{R}^n$, so that each $R_{x_1\dots x_s,ijkl}$
satisfies the following identities:

\begin{enumerate}
\item{$R_{x_1\dots x_s,ijkl}$ \text{is symmetric in }$x_1,\dots
,x_s$} \item{$R_{x_1\dots [x_s,ij]kl}=0$} \item{$R_{x_1\dots
x_s,[ijk]l}=0$} \item{$R_{x_1\dots x_s,ijkl}=-R_{x_1\dots
x_s,jikl}$, $R_{x_1\dots x_s,ijkl}=-R_{x_1\dots x_s,ijlk}$}
\end{enumerate}

where in general,
$T_{r_1\dots r_m[i_1 i_2 i_3]f_1\dots f_d}$ will stand for
 the sum over all the cyclic permutations of the indices
$i_1,i_2 ,i_3$ (in the case where two of the indices
$i_1, i_2,i_3$ are antisymmetric).
\end{definition}

\par Our third building block is the following set:

\begin{definition}
\label{spec}
 Let us consider a set of tensors
$\Xi=\{ {\Xi}^{k_1}_i,\dots {\Xi}^{k_s}_{i_1\dots i_s},
\dots \}$, where the free indices are $i_1,\dots ,i_s,k_s$.
We assume that each tensor ${\Xi}^{k_s}_{i_1\dots i_s}$
 is symmetric in the indices
 $i_1,\dots ,i_s$.
We will call any such tensor a special tensor.
Any such set $\Xi$ will be called a set of special tensors.
\end{definition}

 We can then form complete contractions of tensors that belong
 to the sets $\bigcup_{\alpha\in
A}\{T^\alpha\},R,\Xi$. They will be in the form:

\begin{equation}
\label{fefmyapp2}
contr(u^{l_1}\otimes \dots \otimes u^{l_Z}\otimes
R^{r_1}\otimes \dots \otimes R^{r_m}\otimes {\Xi}^{z_1}
\otimes\dots \otimes {\Xi}^{z_x})
\end{equation}

where each tensor $u^{l_i}$ belongs to the set
$\bigcup_{\alpha\in
A}\{T^\alpha\}$, each tensor
$R_{r_j}$ belongs to the set $R=\{R_{ijkl},\dots ,
R^s_{f_1\dots f_s,ijkl},\dots \}$ and each tensor ${\Xi}^z$
belongs to the set
$\Xi=\{ {\Xi}^k_i,\dots {\Xi}^k_{i_1\dots i_s},\dots \}$.
A particular contraction of two lower indices will be with
 respect to the Kronecker $\delta^{ij}$, while for an upper
 and lower index we will be using the Einstein summation
 convention.
 We can consider linear
combinations of such complete contractions: $\Lambda
(\bigcup_{\alpha\in A}\{T^\alpha\},R, \Xi)={\Sigma}_{l\in L} a_l
C^l(\bigcup_{\alpha\in A}\{T^\alpha\},R,\Xi)$.

\par For each complete contraction $C(\bigcup_{\alpha\in A}\{T^\alpha\},R,\Xi)$
that contains a factor $t=R_{i_1\dots i_s,ijkl}$, we will say that
we apply the third identity in Definition \ref{lincurv}
to the indices $i,j,k$ (or that
we permute indices according to the third identity)
  if we substitute the
complete contraction $C(\bigcup_{\alpha\in A}\{T^\alpha\},R,\Xi)$,
by $-C_1(\bigcup_{\alpha\in
A}\{T^\alpha\},R,\Xi)-C_2(\bigcup_{\alpha\in A}\{T^\alpha\},R,\Xi)$,
where $C_1(\bigcup_{\alpha\in A}\{T^\alpha\},R,\Xi)$ is obtained from
\\ $C(\bigcup_{\alpha\in A}\{T^\alpha\},R,\Xi)$ by substituting
$t$ by $R_{i_1\dots i_s,kijl}$ and $C_2(\bigcup_{\alpha\in
A}\{T^\alpha\},R,\Xi)$ is obtained
 from $C(\bigcup_{\alpha\in A}\{T^\alpha\},R,\Xi)$ by substituting
 $t$ by $R_{i_1\dots i_s,jkil}$.
We similarly define what it means to apply the second identity in
Definition \ref{lincurv}. It is clear what is meant by applying
the first and fourth identities (or by permuting indices according
to the first and
 fourth identities).

\begin{definition}
\label{defvan}
We will say that such a
linear combination of complete contractions vanishes formally
if we can can make the linear combination zero
using the following list of operations:

\par By permuting factors in the complete contractions, by
permuting indices in the factors in $\bigcup_{\alpha\in
A}\{T^\alpha\}$, by using the
 identities of the factors in $R$,
by permuting the indices $i_1,\dots ,i_s$ in the factors
 ${\Xi}^{k_s}_{i_1\dots i_s}$
 and by applying the rule $a\cdot C^l(\bigcup_{\alpha\in
A}\{T^\alpha\},R,\Xi)+ b\cdot C^l(\bigcup_{\alpha\in
A}\{T^\alpha\},R,\Xi)= (a+b)\cdot C^l(\bigcup_{\alpha\in
A}\{T^\alpha\},R,\Xi)$.

\par Also, we will say that
the linear combination $\Lambda(\bigcup_{\alpha\in
A}\{T^\alpha\},R,\Xi)$ {\it vanishes upon
substitution} if for each set of tensors $\bigcup_{\alpha\in
A}\{T^\alpha\}$, $R$ and $\Xi$
 that have the above properties, the value of
$\Lambda(\bigcup_{\alpha\in
A}\{T^\alpha\},R,\Xi)$ is zero.
\end{definition}

\par The following theorem is then an extension of Theorem
B.4 in \cite{beg:itccg} and it follows by the same ideas.

\begin{theorem}
\label{UMxi}
 Let us consider a linear combination of complete contractions
\\ $\Lambda (\bigcup_{\alpha\in
A}\{T^\alpha\},R,\Xi)={\Sigma}_{l\in L} a_l
C^l(\bigcup_{\alpha\in
A}\{T^\alpha\},R,\Xi)$ as above. For each complete contraction $C^l$,
 we denote by $Z^\sharp_l$ the number of symmetric tensors of rank $\ge 1$.
  We also recall that $m_l$ is the number of linearized curvature tensors and
$x_l$ the number of special tensors. We assume that for each $C^l$
the sum $Z^\sharp_l+2m_l+2x_l$  is less than or equal to
 $n$.
\newline

We then have that if $\Lambda (\bigcup_{\alpha\in
A}\{T^\alpha\},R,\Xi)$ vanishes upon substitution in dimension
$n$, it must also vanish formally.
\end{theorem}

\par We note that the theorem above also applies when there
are no factors from the set $\Xi$ in our linear combination.

\subsection{Corollaries of Theorem \ref{UMxi}.}

\par We derive two corollaries of Theorem
 \ref{UMxi}. We will now be considering complete contractions on
 manifolds.

\par Consider an auxiliary list of symmetric tensors
$\Omega=\{\Omega_{i_1},\dots ,\Omega_{i_1\dots i_s},\dots \}$.
We impose the condition that these tensors must remain invariant under re-scalings of the metric $g^n$, ie
$\Omega_{i_1\dots i_s}(t^2g^n)=\Omega_{i_1\dots i_s}(g^n)$.
We then focus our attention on complete contractions
$C^l_{g^n}({\psi}_1,\dots ,{\psi}_Z,\Omega)$ of the form:

\begin{equation}
\label{partlinisymomega}
\begin{split}
&contr({\nabla}^{m_1}_{r_1\dots r_{m_1}}R_{ijkl}\otimes\dots \otimes {\nabla}^{m_s}_{t_1\dots t_{m_s}}R_{ijkl}\otimes
\\& {\nabla}^{p_1}_{a_1\dots a_{p_1}}{\psi}_1\otimes
\dots \otimes {\nabla}^{p_Z}_{b_1\dots b_{p_Z}}{\psi}_Z
\otimes\Omega_{i_1\dots i_{h_1}}\otimes\dots\otimes \Omega_{u_1\dots u_{h_y}})
\end{split}
\end{equation}
 We assume that $y\ge 0$ (in other words, there
may also be no factors $\Omega_{i_1\dots i_s}$).
 If we
 write $C^l_{g^r}({\psi}_1,\dots ,{\psi}_Z,\Omega)$
 (replacing $g^n$ by $g^r$),
we will be referring to a complete contraction as above, but
 defined on an $r$-dimensional manifold.
We will call this the {\it re-writing} of the
complete contraction $C^l_{g^n}({\psi}_1,\dots ,{\psi}_Z,
\Omega)$ in dimension $r$. Also, when
 we speak of the value of $C^l_{g^r}({\psi}_1,\dots ,{\psi}_Z,
\Omega)(x_0)$, we will mean the value of the
 above complete contraction at a point $x_0$ on a manifold
$(M^r,g^r)$, for functions ${\psi}_1,\dots ,{\psi}_Z$ defined
 around $x_0\in M^r$ and for symmetric tensors
$\Omega_{i_1\dots i_s}$
defined at $x_0$.  This terminology extends to linear combinations.

\par Finally, a note about the weight of the complete
 contractions: By our definition of weight, if
$C^l_{g^r}({\psi}_1,\dots ,{\psi}_Z,\Omega)$ has weight $-n$,
then in the notation of (\ref{partlinisymomega}):

\begin{equation}
\label{proxwrw}
\Sigma_{i=1}^s (m_i+2)+\Sigma_{i=1}^Z p_i+\Sigma_{i=1}^y h_i=n
\end{equation}
Thus, if we have $Z^\sharp$ factors $\nabla^{p_i}\psi_i$ with
 $p_i\ge 1$, the above implies that:

\begin{equation}
\label{proxwrw2}
Z^\sharp+2s+y\le n
\end{equation}

\begin{definition}
\label{realform}
 We will say that a relation between complete contractions in the
 form (\ref{partlinisymomega}):

$${\Sigma}_{l\in L} a_l C^l_{g^n}(\psi_1,\dots ,\psi_Z,\Omega)=0$$
 will hold {\it formally} if we can make the above sum
 identically zero by performing the following operations:
We may permute factors in any complete contraction
$C^l_{g^n}(\psi_1,\dots ,\psi_Z)$ and also permute indices within
the factors $\Omega_{i_1\dots i_s}$. Furthermore, for each factor
$\nabla^p_{r_1\dots r_p}\psi_h$, with $p=2$ we may permute
$r_1,r_2$, while for $p>2$, we may apply the identity:

\begin{equation}
\label{number} [{\nabla}_A{\nabla}_B - {\nabla}_B{\nabla}_A] X_C=R_{ABCD}X^D
\end{equation}
 and  for each factor
${\nabla}^mR_{ijkl}$, we may apply the identities:

\begin{enumerate}

\item{ ${\nabla}^m_{r_1\dots r_m}R_{ijkl}=-{\nabla}^m_{r_1\dots
r_m}R_{jikl}=-{\nabla}^m_{r_1\dots r_m}R_{ijlk}$.}

\item{ ${\nabla}^m_{r_1\dots [r_m}R_{ij]kl}=0$}

\item{ ${\nabla}^m_{r_1\dots r_m}R_{[ijk]l}=0$}

\item{The identity (\ref{number}) above.}
\end{enumerate}

\par The application of the second and third identities above
 has been defined. To apply the fourth identity to a factor
 $\nabla^p\psi_h$ or $\nabla^mR_{ijkl}$ means that for
 each complete contraction $C_{g^n}({\psi}_1,\dots ,{\psi}_Z,\Omega)$
 of the form (\ref{fefmyapp2}), for each factor
$\nabla^p_{r_1\dots r_p}\psi_h$ or
 ${\nabla}^m_{r_1\dots r_m}R_{ijkl}$ in
$C_{g^n}({\psi}_1,\dots ,{\psi}_Z,\Omega)$ and each pair of consecutive
derivative indices $r_{s-1},r_s$ we may write:
$$C_{g^n}({\psi}_1,\dots ,{\psi}_Z,\Omega)=C'_{g^n}({\psi}_1,\dots
,{\psi}_Z,\Omega) + {\Sigma}_{h\in H} a_h C^h_{g^n}({\psi}_1,\dots ,
\psi_Z,\Omega)$$

  where $C'_{g^n}({\psi}_1,\dots ,{\psi}_Z,\Omega)$ is obtained from
$C_{g^n}({\psi}_1,\dots ,{\psi}_Z,\Omega)$
 by substituting the factor $\nabla^p_{r_1\dots r_p}\psi_h$ or ${\nabla}^m_{r_1\dots r_m}
R_{ijkl}$ by $\nabla^p_{r_1\dots r_sr_{s-1}\dots r_p}\psi_h$ or
${\nabla}^m_{r_1\dots r_s r_{s-1}\dots  r_m} R_{ijkl}$,
respectively, and ${\Sigma}_{h\in H} a_h C^h_{g^n}({\psi}_1,\dots
,{\psi}_Z,\Omega)$ is obtained from $C_{g^n}({\psi}_1,\dots
,{\psi}_Z,\Omega)$ by substituting the factor
$\nabla^p_{r_1\dots r_p}\psi_h$
${\nabla}^m_{r_1\dots r_m}R_{ijkl}$ by one of the
summands in the following expressions, respectively, and then
summing again:

\begin{equation}
\label{prin9atobrw}
\begin{split}
&\Sigma_{\{a_1,\dots a_x\},\{b_1,\dots b_{s-2-x}\}
\subset\{r_1,\dots r_{s-2}\}, \{a_1,\dots a_x\}\bigcap\{b_1,\dots
b_{s-1-x}\}=\emptyset}
\\&(\nabla^x_{a_1\dots a_x}{R_{r_{s-1}r_sr_{s+1}}}^d)
(\nabla^{s-1-x}_{b_1\dots b_{s-1-x}}\nabla^{m-s-1}_{dr_{s+2}\dots
r_p} \psi_h+\dots
\\& +(\nabla^x_{a_1\dots a_x}
{R_{r_{s-1}r_sr_p}}^d)(\nabla^{s-1-x}_{b_1\dots b_{s-1-x}}
\nabla^{m-s-1}_{r_{s+1}\dots d}\psi_h)
\end{split}
\end{equation}

\begin{equation}
\label{9atobrw}
\begin{split}
&\Sigma_{\{a_1,\dots a_x\},\{b_1,\dots b_{s-2-x}\}
\subset\{r_1,\dots r_{s-2}\}, \{a_1,\dots a_x\}\bigcap\{b_1,\dots
b_{s-2-x}\}=\emptyset}
\\&(\nabla^x_{a_1\dots a_x}{R_{r_{s-1}r_{s}r_{s+1}}}^d)
(\nabla^{s-1-x}_{b_1\dots b_{s-2-x}}\nabla^{m-s-1}_{dr_{s+2}\dots
r_m} R_{ijkl})+\dots
\\& +(\nabla^x_{a_1\dots a_x}
{R_{r_{s-1}r_{s}l}}^d)(\nabla^{s-1-x}_{b_1\dots b_{s-1-x}}
\nabla^{m-s-1}_{r_{s+1}\dots r_m}R_{ijkd})
\end{split}
\end{equation}
\end{definition}

\par Now, our first corollary of Theorem \ref{UMxi}:

\begin{lemma}
\label{corUM}
Consider complete contractions
$C^l_{g^n}({\psi}_1,\dots ,{\psi}_Z,
\Omega)$, each in the form (\ref{partlinisymomega}) and with weight
$-n$, so that the following identity:

\begin{equation}
\label{subtofor1}
F_{g^n}(\psi_1,\dots ,\psi_Z,\Omega)=
\Sigma_{l\in L} a_l C^l_{g^n}(\psi_1,\dots ,
\psi_Z,\Omega)=0
\end{equation}
 holds at any point $x_0$ for any metric $g^n$
and any functions ${\psi}_1,\dots ,{\psi}_Z$ defined around
 $x_0$ and any symmetric tensors $\Omega_{i_1\dots i_s}$ defined over
$TM^n|_{x_0}$. We then have that the above identity must hold
 formally.
\end{lemma}

{\it Proof:} Let us consider the minimum length $\tau$, among
all the complete contractions in (\ref{subtofor1}).
 Let us index the complete
contractions $C^l_{g^n}({\psi}_1,\dots ,{\psi}_Z,\Omega)$ of length
 $\tau$ in the set $L^{\tau}\subset L$. Suppose we can show
 that, applying the above operations, we can make
${\Sigma}_{l\in L^{\tau}} a_l C^l_{g^n}({\psi}_1,\dots ,
{\psi}_Z,\Omega)$ formally
equal to a linear combination ${\Sigma}_{r\in R} a_r C^r_{g^n}(\psi_1,\dots ,\psi_Z,\Omega)$, where each complete contraction
$C^r_{g^n}(\psi_1,\dots ,\psi_Z,\Omega)$ has length
$\ge \tau +1$.

\par If we can prove the above claim then using a finite
 number of
 iterations we will have proven our Lemma. This is true since
 there is obviously a number $T$, so that all the complete
 contractions that arise by iteratively applying the
 identities of Definition \ref{realform} to the complete
 contractions $C^l_{g^n}(\psi_1,\dots ,\psi_Z,\Omega)$,
$l\in L$ must have length $\le T$. This follows just by the
finiteness of the index set $L$.
 The rest of this proof will focus on showing that
 claim.

\par In order to accomplish this, let us begin with a
 definition. For any complete contraction
$C^l_{g^n}({\psi}_1,\dots ,{\psi}_Z,\Omega)$, let $lin
C^l(R,{\Psi}_1,\dots ,{\Psi}_Z,\Omega)$ stand for the complete
contraction between linearized curvature tensors and symmetric
tensors that is obtained from $ C^l_{g^n}({\psi}_1,\dots
,{\psi}_Z,\Omega)$ by substituting each factor
${\nabla}^m_{t_1\dots t_m}R_{ijkl}$ by a linearized curvature
tensor $R_{t_1\dots t_m,ijkl}$, and each factor
$\nabla^p_{r_1\dots r_p}\psi_h$ by a symmetric $p$-tensor
$\Psi^h_{r_1\dots r_p}$. We will prove a fact which will be
used many times in the
 future.

\begin{lemma}
\label{realtolin}
In the above notation, given (\ref{subtofor1}),
we will have that:
$${\Sigma}_{l\in L^{\tau}} a_l linC^l(R,{\Psi}_1,\dots ,
{\Psi}_Z,\Omega)=0$$ formally.
\end{lemma}

{\it Proof:} We recall the following fact, which follows from the
proof of Theorem 2.6 in \cite{e:ntrm}: Given any set $R$ of
linearized curvature tensors $R_{t_1\dots t_m,ijkl}(x_0)$, there
is a Riemannian metric defined around $x_0$ so that for any
$m$:

\begin{equation}
\label{lintoreal}
({\nabla}^m_{t_1\dots t_m}R_{ijkl})^{g^n}(x_0)=
R_{t_1\dots t_m,ijkl}(x_0) +C(R)_{t_1\dots t_m,ijkl}
\end{equation}
where $C(R)_{t_1\dots t_m,ijkl}$ stands for
a polynomial in the components of the linearized curvature
 tensors. We have that this polynomial depends only on $m$ and
the indices $t_1,\dots ,t_m,i,j,k,l$. Furthermore, we have that
each monomial in that polynomial will have degree at
 least 2.

\par For any set $R$ of linearized curvature tensors,
we will call the metric $g^n$ for which (\ref{lintoreal})
holds the {\it associated} metric.
 Now, for any choice of symmetric tensors
$\{ T^1_0,T^1_i,\dots T^1_{i_1\dots i_s},\dots\}$,
 $\dots$,
$\{T^Z_0,T^Z_i,\dots , T^Z_{i_1\dots i_s},\dots \}$,
we have that there are functions ${\psi}_1,\dots ,{\psi}_Z$
 defined around $x_0$ so that:
${\nabla}^s_{i_1\dots i_s}{\psi}_l(x_0)= T^l_{i_1\dots i_s}$ (for
some arbitrary ordering of the indices $i_1,\dots,i_s$ on the left
hand side), and also for each permutation $\pi(i_1\dots i_s)$ of
the indices $i_1,\dots ,i_s$:

\begin{equation}
\label{lintoreal2} {\nabla}^p_{\pi(i_1\dots i_s)}\psi_h(x_0)=
{\nabla}^p_{i_1\dots i_s}\psi_h(x_0)  +C(R,\psi_h)_{i_1\dots i_s}
\end{equation}
where $C(R,T^h)_{i_1\dots i_s}$ stands for a polynomial in the
components of the linearized curvature
 tensors and of one component of a tensor from the set $T^h$
(of rank $\ge 1$).
  We have that this polynomial depends only on $p$ and
the indices $i_1,\dots ,i_s$. Furthermore, we have that each
monomial in that polynomial will have degree at
 least 2.

 For any choice of symmetric tensors $T^l_{i_1\dots i_s}$, we
define the functions ${\psi}_l$ to be their {\it associated}
functions.

\par Now, we pick any set $R$ of linearized curvature tensors
 and any set $T$ of symmetric tensors and consider the value
 of $F_{g^n}({\psi}_1,\dots ,{\psi}_Z,\Omega)$ for the associated
 metric $g^n$ and the associated functions ${\psi}_l$.
 By virtue of our remarks, we see that there is a
 fixed polynomial $\Pi(T,R,\Omega)$ in the vector space of
components of the sets $T$ and $R$, so that for any given set
 $R$ of linearized curvature tensors and any set $T$ of
 symmetric tensors at $x_0$,

$$\Pi(T,R,\Omega)=F_{g^n}({\psi}_1,\dots ,{\psi}_Z,\Omega)=0$$

\par Furthermore, we observe from (\ref{lintoreal}) that
 each monomial in $\Pi(T,R,\Omega)$
has degree at least $\tau$. Finally, if $\Pi(T,R,\Omega)|_{\tau}$
stands for the sublinear combination of monomials of degree
$\tau$
in $\Pi(T,R,\Omega)$, we have that:

 $$\Pi(T,R,\Omega)|_{\tau}=0$$

for every set $R$ of linearized curvature tensors and every
 sets $T$, $\Omega$ of symmetric tensors. But given equations
(\ref{lintoreal}) and (\ref{lintoreal2}) we see that:

\begin{equation}
\label{click} \Pi(T,R,\Omega)|_{\tau}={\Sigma}_{l\in L^{\tau}} a_l
linC^l(R,{\psi}_1,\dots ,{\psi}_Z,\Omega)=0
\end{equation}

\par Hence, in view of Theorem \ref{UMxi}, we have that
 (\ref{click}) must hold formally. $\Box$
\newline

So, for each
$linC^l_{g^n}({\psi}_1,\dots ,{\psi}_Z,\Omega)$
 there is a sequence of permutations for the factors
$\Psi^l_{t_1\dots t_a}$, $,\Omega_{i_1\dots i_s}$
 and of applications
 of the identities of a linearized curvature tensor to the factors
$R_{t_1\dots t_m,ijkl}(x_0)$ so that (\ref{click}) will hold by
virtue of the identity $a\cdot
C(\bigcup_{i=1}^Z\{T^i\},R,\Omega)+b\cdot
C(\bigcup_{i=1}^Z\{T^i\},R, \Omega)=(a+b)\cdot
C(\bigcup_{i=1}^Z\{T^i\},R,\Omega)$.

\par We then repeat these operations to the sublinear combination
${\Sigma}_{l\in L^{\tau}} a_l \\ C^l_{g^n}({\psi}_1,\dots ,
{\psi}_Z,\Omega)$.
 The only difference is that the indices $t_1,\dots ,t_m$
in each factor ${\nabla}^m_{t_1\dots t_m}R_{ijkl}(x_0)$ and the
indices $i_1,\dots ,i_p$ in each factor $\nabla^s_{i_1\dots
i_s}\psi_h$ are not symmetric. Nonetheless, we may permute the
indices $i_1,\dots ,i_s$ in each factor $\nabla^s_{i_1\dots
i_s}\psi_h$ and  the indices $t_1,\dots ,t_m$ in each factor
${\nabla}^m_{t_1\dots t_m}R_{ijkl}$ and introduce correction
 terms, which are complete contractions
in the form (\ref{partlinisym}) of length $\ge \tau +1$.
 Hence, repeating the permutations
which made (\ref{click}) identically zero, we have our
claim. $\Box$
\newline

\par We now make a note about the notation we will be
 using: We have considered complete contractions $C^l_{g^n}(\psi_1,
\dots ,\psi_Z,\Omega)$ in the general form
(\ref{partlinisymomega}), and we have explained that there may
also be no factors $\Omega_{i_1\dots i_s}$. We make the extra
convention that if we refer to a complete contraction
$C^l_{g^n}(\psi_1,\dots ,\psi_Z)$, we will imply that it is in the
form (\ref{partlinisymomega}) {\it and} has no factors
$\Omega_{i_1\dots i_s}$. Therefore, it will be in the form:

\begin{equation}
\label{partlinisym}
\begin{split}
&contr({\nabla}^{m_1}_{r_1\dots r_{m_1}}R_{ijkl}\otimes\dots \otimes
{\nabla}^{m_s}_{t_1\dots t_{m_s}}R_{ijkl}\otimes
\\& {\nabla}^{\nu_1}_{a_1\dots a_{\nu_1}}{\psi}_1\otimes
\dots \otimes {\nabla}^{\nu_Z}_{b_1\dots b_{\nu_Z}}{\psi}_Z)
\end{split}
\end{equation}

\par Our next Lemma is another corollary of Theorem
\ref{UMxi}. We must again introduce a definition.

\par We focus on complete contractions $C^l_{g^n}({\psi}_1,
\dots ,{\psi}_Z,\Xi)$ of the form:

\begin{equation}
\label{partlinisymxi}
\begin{split}
&contr({\nabla}^{m_1}_{r_1\dots r_{m_1}}R_{ijkl}\otimes\dots \otimes
{\nabla}^{m_s}_{t_1\dots t_{m_s}}R_{ijkl}\otimes
\\& {\nabla}^{p_1}_{a_1\dots a_{p_1}}{\psi}_1\otimes
\dots \otimes {\nabla}^{p_Z}_{b_1\dots b_{p_Z}}{\psi}_Z \otimes {\Xi}^{k_1}_{i_1\dots i_s} \otimes\dots\otimes
{\Xi}^{k_f}_{j_1\dots j_t})
\end{split}
\end{equation}

\par In the manifold context, we impose the re-scaling
condition ${\Xi}^{k_1}_{i_1\dots i_s}(t^2g^n)=
{\Xi}^{k_1}_{i_1\dots i_s}(g^n)$ on the special tensors.
When we wish to apply the theorem to a particular case of
 special tensors, we will easily see that this condition
 holds.

\begin{definition}
\label{realformxi}
 We will say that a relation between complete contractions
 in the form (\ref{partlinisymxi}):

$${\Sigma}_{l\in L} a_l C^l_{g^n}({\psi}_1,\dots ,
{\psi}_Z,\Xi)=0$$
 will hold {\it formally} if we can make the above sum
 identically zero by performing the following operations:
We may interchange factors in any complete contraction
$C^l_{g^n}({\psi}_1,\dots ,{\psi}_Z)$ and  also permute the
 indices $i_1,\dots ,i_s$ among each factor
${\Xi}^k_{i_1\dots i_s}$. Furthermore, for
 each factor ${\nabla}^mR_{ijkl}$, we may apply the
identities:

\begin{enumerate}

\item{ ${\nabla}^m_{r_1\dots r_m}R_{ijkl}=-\nabla^m_{r_1\dots
r_m}R_{jikl}=-\nabla^m_{r_1\dots r_m}R_{ijlk}$.}

\item{ ${\nabla}^m_{r_1\dots [r_m}R_{ij]kl}=0$}

\item{ ${\nabla}^m_{r_1\dots r_m}R_{[ijk]l}=0$}

\item{ We may permute the indices $r_1,\dots ,r_m$ by applying of the formula:
$[{\nabla}_A{\nabla}_B -
{\nabla}_B{\nabla}_A] X_C= R_{ABCD}X^D$, as defined in the
 previous definition.}
\end{enumerate}
and for any factor $\nabla^p_{i_1\dots i_p}\psi_h$ we may permute
the factors $i_1,i_2$ if $p=2$ and apply the identity
$[{\nabla}_A{\nabla}_B -
{\nabla}_B{\nabla}_A] X_C= R_{ABCD}X^D$, as defined in the previous definition if
$p>2$.
\end{definition}

\par We then have:

\begin{lemma}
\label{corUMxi} Consider complete contractions
$C^l_{g^n}({\psi}_1,\dots , {\psi}_Z, \Xi)$, each in the form
(\ref{partlinisymxi}) and with weight $-n$, so that the identity:

\begin{equation}
\label{subtofor}
{\Sigma}_{l\in L} a_l C^l_{g^n}({\psi}_1,\dots ,
{\psi}_Z,\Xi)=0
\end{equation}

 holds at any point $x_0$, for any metric $g^n$, any
  functions ${\psi}_1,\dots ,{\psi}_Z$  defined around
 $x_0$ and any special tensors
${\Xi}^k_{i_1\dots i_s}(x_0)$ defined at $x_0$. Assume also that each special tensor in each $C^l$ has rank at least 4.
 We then have that the above identity must hold formally.
\end{lemma}

{\it Proof:} We prove this corollary by using Theorem
\ref{UMxi}, in the same way that we proved Lemma \ref{corUM}
 using Theorem \ref{UMxi}.

\par We only need to observe that for each complete
 contraction in the form (\ref{partlinisymxi}) with weight
$-n$,
if $r_i$ stands for the rank of the $i^{th}$ special tensor
 then:

\begin{equation}
\label{proxwrw3}
\Sigma_{i=1}^s (m_i+2)+\Sigma_{i=1}^Z p_i+\Sigma_{i=1}^f
(r_i-2)=n
\end{equation}
For each $C^l_{g^n}({\psi}_1,\dots ,{\psi}_Z,
\Xi)$, we again denote by $Z^\sharp$ the number of factors
$\nabla^{p_h}\psi_h$ for which $p_h\ne0$.
Thus, since we are assuming that each special factor has rank
 at least 4, we deduce that for each complete contraction
$C^l_{g^n}({\psi}_1,\dots ,{\psi}_Z,\Xi)$:

\begin{equation}
\label{proxwrw4}
Z^\sharp +2s+2f \le n
\end{equation}

\par Let $\tau$ be the minimum length among all the  contractions $C^l_{g^n}({\psi}_1,\dots ,{\psi}_Z,
\Xi)$, $l\in L$. We define the subset $L^{\tau}\subset L$
to be the index set of all complete contractions
$C^l_{g^n}({\psi}_1,\dots ,{\psi}_Z, \Xi)$ with length
$\tau$. As before, we define the linear combination of
 complete contractions involving linearized curvature tensors
 rather than ``genuine'' covariant derivatives of the curvature
 tensor, and also symmetric tensors $\Psi^h$ rather than ``genuine''
 factors $\nabla^p\psi_h$:

$${\Sigma}_{l\in L^{\tau}} a_l linC^l(R, {\Psi}_1,\dots ,
{\Psi}_Z,\Xi)$$

 and we show that

$${\Sigma}_{l\in L^{\tau}} a_l linC^l(R, {\Psi}_1,\dots ,
{\Psi}_Z,\Xi)=0$$ formally. We then deduce that an equation:

\begin{equation}
\label{fulltolin2}
{\Sigma}_{l\in L^{\tau}} a_l C^l_{g^n}({\psi}_1,\dots ,
{\psi}_Z,\Xi)= {\Sigma}_{r\in R} a_r
C^r_{g^n}({\psi}_1,\dots ,{\psi}_Z,\Xi)
\end{equation}

where each $C^r_{g^n}({\psi}_1,\dots ,{\psi}_Z,\Xi)$
has length $\ge \tau +1$, will hold formally. By inductive
 repetition of this argument, we have our Lemma.
$\Box$
\newline

\par These Lemmas will prove useful in the future. For now, we
note that there are many definitions of an identity holding
formally. However,
 there will be no confusion, since in each of the above cases
the complete contractions involve
tensors that belong to different categories.
Furthermore, in spite of the equivalence that the above theorems
and their corollaries imply, whenever we mention an identity in this
paper, we will mean (unless we explicitly
state otherwise) that it holds at any point and for every
metric and set of functions (and maybe special tensors
$\Xi$ or symmetric tensors $\Omega$).

\subsection{The Isomorphism.}

\par We now conclude that:

\begin{proposition}
\label{propdim}
Suppose that
$\{ C^r_{g^N}({\psi}_1,\dots ,{\psi}_Z)\}_{r\in R}$
are complete contractions in the form
(\ref{partlinisym}) of weight $-n$.
 Suppose $N\ge n$. We then have that

$${\Sigma}_{r\in R} a_r C^r_{g^N}({\psi}_1,\dots ,{\psi}_Z)
(x_0)=0$$
for every $(M^n,g^n)$, every $x_0\in M^n$ and any
functions ${\psi}_l$ defined around $x_0$  if and only if:

$${\Sigma}_{r\in R} a_r C^r_{g^n}({\psi}_1,\dots ,{\psi}_Z)
(x_0)=0$$
for every $(M^n,g^n)$, every $x_0\in M^n$ and any
functions ${\psi}_l$ defined around $x_0$.
\end{proposition}

{\it Proof:} The above follows by virtue of Lemma
\ref{corUMxi}. $\Box$

\section{The silly divergence formula.}

\par Our aim here is to obtain a formula that expresses
$I^Z_{g^n}(\psi_1,\dots ,\psi_Z)$ as a divergence of a Riemannian
vector field. This first, rather easy, divergence formula is not
useful in itself. It will be used, however, in the derivation of
the much more subtle {\it simple divergence
 formula} in the next section. For now, we claim:

\begin{proposition}
\label{iisdiv}
Consider any $I^Z_{g^n}(\psi_1,\dots ,\psi_s)$, a
linear combination of contractions in the form
(\ref{linicontrpsi})
for which $\int_{M^n}I^Z_{g^n}(\psi_1,\dots ,\psi_s)
dV_{g^n}=0$ for
 every compact $(M^n,g^n)$ and any $\psi_1,\dots ,\psi_s\in
C^\infty (M^n)$. Note that $I^Z_{g^n}(\psi_1,\dots ,\psi_s)$ defined in
 (\ref{derlami}) satisfies this property.

\par We then have that $I^Z_{g^n} (\psi_1,\dots ,\psi_Z)(x)$ is
formally equal to the divergence of a Riemannian vector-valued
differential operator of weight $-n+1$
in $\psi_1(x),\dots , \psi_Z(x)$.
\end{proposition}

{\it Proof:} In view of Lemma \ref{corUM} in the previous
subsection, it suffices to show that there is a vector field
$T^i_{g^n}(\psi_1,\dots ,\psi_Z)$ of weight $-n+1$ so that:

$$I^Z_{g^n}(\psi_1,\dots
 ,\psi_Z)(x_0)= div_i T^i_{g^n}(\psi_1,\dots ,\psi_Z)(x_0)$$

for any metric $g^n$ and for any functions ${\psi}_1,\dots
,{\psi}_Z$ around $x_0$. In order to show
 this we do the following:

 Suppose that

 $$I^Z_{g^n}({\psi}_1,\dots  ,{\psi}_Z)= {\Sigma}_{j\in J} a_j
 C^j_{g^n}({\psi}_1,\dots  ,{\psi}_Z)$$

   where each of the complete contractions $C^j_{g^n}({\psi}_1,\dots
   ,{\psi}_Z)$ is in the form (\ref{linicontrpsi}). Let us sort
   out the different values of ${\nu}_1$ that can appear among the
   different complete contractions $C^j_{g^n}({\psi}_1,\dots
   ,{\psi}_Z)$. Suppose the set of those different values is the
   set $L=\{{\lambda}_1,\dots ,{\lambda}_K\}$ where
   $0\le{\lambda}_1< \dots < {\lambda}_K$.

   \par Let $J_K\subset J$ be the set of the complete contractions
   $C^j_{g^n}({\psi}_1,\dots
   ,{\psi}_Z)$ with ${\nu}_1={\lambda}_K$. We then consider the
   linear combination :

   $$F_{g^n}({\psi}_1,\dots  ,{\psi}_Z)= {\Sigma}_{j\in J_K} a_j
 C^j_{g^n}({\psi}_1,\dots  ,{\psi}_Z)$$

where each complete contraction $C^j_{g^n}(\psi_1,\dots
,{\psi}_Z)$ is in the form (\ref{linicontrpsi}) with the same
number ${\lambda}_K$ of derivatives on ${\psi}_1$. Out of
$F_{g^n}({\psi}_1,\dots  ,{\psi}_Z)$, we
construct the following vector-valued differential operator:

$$F^i_{g^n}({\psi}_1,\dots  ,{\psi}_Z)= {\Sigma}_{j\in J_K} a_j
 (C^j)^i_{g^n}({\psi}_1,\dots  ,{\psi}_Z)$$

where $(C^j)^i_{g^n}({\psi}_1,\dots  ,{\psi}_Z)$ is made out of
 $C^j_{g^n}({\psi}_1,\dots  ,{\psi}_Z)$ by erasing the index
 ${\chi}_1$ in (\ref{linicontrpsi}) and making the index that
 contracted against it in (\ref{linicontrpsi}) into a free
 index.

\par Let us then observe the following:

\begin{lemma}
\label{crude} We have that the differential operator
$$\tilde{F}_{g^n}({\psi}_1,\dots  ,{\psi}_Z)=
F_{g^n}({\psi}_1,\dots
,{\psi}_Z)- div_i F^i_{g^n}({\psi}_1,\dots  ,{\psi}_Z)$$
will be
formally equal to a linear combination of complete contractions in
 the form (\ref{linicontrpsi}) (of weight $-n$), each of
 which has ${\lambda}_K-1$
 derivatives on the function ${\psi}_1$.
\end{lemma}

 {\it Proof:} This is straightforward
 to observe by the construction of the vector-valued operators
 $(C^j)^i_{g^n}({\psi}_1,\dots
 ,{\psi}_Z)$: Let the derivative ${\nabla}_i$ in the divergence of
 each $(C^j)^i_{g^n}({\psi}_1,\dots  ,{\psi}_Z)$
 hit the factor ${\nabla}^{{\lambda}_K-1}{\psi}_1$. That summand
 in the divergence will cancel out the complete contraction $C^j_{g^n}({\psi}_1,\dots
 ,{\psi}_Z)$. Every other complete contraction in
$div_i F^i_{g^n}({\psi}_1,\dots  ,{\psi}_Z)$
will have ${\lambda}_K-1$ derivatives on ${\psi}_1$.
 This gives our desired conclusion. $\Box$

\par But then repeated application of Lemma \ref{crude} will give
us the following:
\newline

\par We can subtract a divergence of a vector field $L^i_{g^n}
({\psi}_1,\dots  ,{\psi}_Z)$ of weight $-n+1$
from $I^Z_{g^n}({\psi}_1,\dots  ,{\psi}_Z)$, so that

$$R_{g^n}({\psi}_1,\dots  ,{\psi}_Z)= I^Z_{g^n}({\psi}_1,\dots  ,{\psi}_Z)-
div_iL^i_{g^n}({\psi}_1,\dots  ,{\psi}_Z)$$

 is a linear combination of complete contractions in the form
 (\ref{linicontrpsi}), each of which has ${\nu}_1=0$.

 \par We then observe that:

\begin{lemma}
\label{crude2}
 In the above notation, $R_{g^n}({\psi}_1,\dots
 ,{\psi}_Z)$ must vanish formally.
 \end{lemma}

{\it Proof:} First observe that for any Riemannian manifold
$(M^n,g^n)$ we will have:

$${\int}_{M^n} R_{g^n}({\psi}_1,\dots
 ,{\psi}_Z) dV_{g^n}=0$$

  This is straightforward to observe, because of Lemma 5 and the
  definition of $R_{g^n}({\psi}_1,\dots
 ,{\psi}_Z)$-it is obtained from $I^Z_{g^n}({\psi}_1,\dots
 ,{\psi}_Z)$ by subtracting a divergence.

 \par Now, let us write $R_{g^n}({\psi}_1,\dots
 ,{\psi}_Z)$ as follows:

 $$R_{g^n}({\psi}_1,\dots
 ,{\psi}_Z)={\Sigma}_{l\in L} a_l C^l_{g^n}({\psi}_2,\dots
 ,{\psi}_Z)\cdot{\psi}_1$$

  We will then have that the equation:

  \begin{equation}
  \label{onlypsi}
{\int}_{M^n} R_{g^n}({\psi}_1,\dots
 ,{\psi}_Z) dV_{g^n}= {\int}_{M^n} {\psi}_1\cdot[{\Sigma}_{l\in L} a_l
 C^l_{g^n}({\psi}_2,\dots
 ,{\psi}_Z)]dV_{g^n}
  \end{equation}

  holds for any function ${\psi}_1$, and also the sum ${\Sigma}_{l\in L}
  a_l C^l_{g^n}({\psi}_2,\dots
 ,{\psi}_Z)$ is independent of the function ${\psi}_1$. But that
 shows us that $R_{g^n}({\psi}_1,\dots  ,{\psi}_Z)$ must vanish by
 substitution. Hence, by Theorem \ref{UMxi} , it must vanish formally. $\Box$

\section{The simple divergence formula.}
\label{simplediv}
\subsection{The transformation law for $I^Z_{g^N}$ and
 definitions.}

\par Let $I^Z_{g^{n}} ({\psi}_1,\dots ,{\psi}_Z)$ be as in Proposition \ref{iisdiv}.
 We then have that $I^Z_{g^{n}} ({\psi}_1,\dots ,{\psi}_Z)$
 is a divergence
of a vector-valued differential operator in ${\psi}_1(x),\dots
,{\psi}_Z(x)$.
 This is useful in itself, but we cannot extract information directly from
this fact about $P(g^n)$. Nevertheless, it is useful in the
following respect:

\par We have a relation:
\begin{equation}
\label{idivn} I^Z_{g^n}({\psi}_1,\dots  ,{\psi}_Z)=
div_iL^i_{g^n}({\psi}_1,\dots  ,{\psi}_Z)
\end{equation}

 which holds formally.
But then Proposition \ref{propdim} tells us the following:

\begin{lemma}
\label{freedim}
 Relation (\ref{idivn}) holds for any dimension $N\ge n$. That is,
considering the complete contractions and the Riemannian vector
fields in (\ref{idivn}) in any dimension $N\ge n$, we have the
formula:

\begin{equation}
\label{idivN} I^Z_{g^N}({\psi}_1,\dots  ,{\psi}_Z)=
div_iL^i_{g^N}({\psi}_1,\dots  ,{\psi}_Z)
\end{equation}
\end{lemma}

{\it Proof:} Straightforward from Propositions \ref{propdim}
and \ref{iisdiv}. $\Box$
\newline

\par Therefore, we will have that for any
$(M^N,g^N)$ and any ${\psi}_1,\dots ,{\psi}_Z\in
C^{\infty}(M^N)$:

\begin{equation}
\label{nbig} {\int}_{M^N}I^Z_{g^N} ({\psi}_1,\dots
,{\psi}_Z)(x)dV_{g^N}=0
\end{equation}

\par Now, equation (\ref{nbig}) is {\it not} scale-invariant.
This can be used to our advantage in the following way:
\newline

\par Pick out any point $x_0\in M^N$. Pick
out a small geodesic ball around $x_0$, of radius $\epsilon$.
From now on, we will assume the functions
${\psi}_1,\dots ,{\psi}_Z$ to be compactly
supported in $B(x_0,\epsilon)$. Then we can pick any coordinate
system around $x_0$ and write out $I^Z_{g^N} ({\psi}_1,\dots
,{\psi}_Z)(x)$ in that coordinate system. In that coordinate
system we will have that:

\begin{equation}
\label{bignint}
{\int}_ {\mathbb{R}^N} I^Z_{g^N}({\psi}_1,\dots
 ,{\psi}_Z)(x)dV_{g^n} =0
\end{equation}

\par Now, let our coordinate system around $x_0$ be $\{ x_1,\dots
,x_N\}$. For that coordinate system, we will denote each point in
$B(x_0,\epsilon)$ by $\vec{x}$. Let also $\vec{\xi}$ be an
arbitrary vector in $\mathbb{R}^N$.
We then consider the following
conformal change of metric in $B(x_0,\epsilon)$:
 ${\hat{g}}^N(x)=
e^{2\vec{\xi}\cdot\vec{x}}g^N(x)$. We have that (\ref{bignint})
must also hold for this metric. The volume form will re-scale as
follows:

$$dV_{\hat{g}^N}(x)=e^{N\vec{\xi}\cdot\vec{x}}
dV_{g^n}(x)$$

\par Now, we have that $I^Z_{g^N}({\psi}_1,\dots
 ,{\psi}_Z)(x)$ is a linear combination of complete
 contractions in the form (\ref{linicontrpsi}). So, in order to
 find how any given complete contraction in the form
 (\ref{linicontrpsi}) transforms under the above conformal change,
 it suffices to find how each of its factors will transform. In
 order to do that, we can employ the identities of the first
 section.

\par The transformation law of Ricci curvature, for this special conformal transformation,
is given by equation (\ref{ricci}), substituting $\phi$ by
$\vec{x}\cdot\vec{\xi}$. Recall that $\nabla_i(\vec{\xi}\cdot
\vec{x})=\vec{\xi}_i$.

\begin{equation}
\label{riccixi} {Ric}^{\hat{g}^N}_{\alpha\beta}(x)=
{Ric}^{g^N}_{\alpha\beta}(x)+
(2-N){\nabla}^2_{\alpha\beta}(\vec{\xi}\cdot\vec{x}) -
{\Delta}_{g^N} (\vec{\xi}\cdot\vec{x})g^N_{\alpha\beta}
+(N-2)({\vec{\xi}}_{\alpha}{\vec{\xi}}_{\beta}-
{\vec{\xi}}^k{\vec{\xi}}_kg^N_{\alpha\beta})
\end{equation}

 \par The scalar curvature will transform as:

\begin{equation}
\label{scalarxi} R^{\hat{g}^N}(x)= e^{-2\vec{\xi}\cdot
\vec{x}}[R^{g^N}+ 2(1-N){\Delta}_{g^N}(\vec{\xi}\cdot\vec{x})
-(N-1)(N-2){\vec{\xi}}^k{\vec{\xi}}_k]
\end{equation}

\par And the full curvature tensor:

\begin{equation}\begin{split}
\label{curvxi} & R^{\hat{g}^N}_{ijkl}(x)=
e^{2\vec{\xi}\cdot\vec{x}}\{ R^g_{ijkl}(x) +
[{\vec{\xi}}_i{\vec{\xi}}_k g^N_{jl} -{\vec{\xi}}_i{\vec{\xi}}_l
g^N_{jk} +{\vec{\xi}}_j{\vec{\xi}}_l g^N_{ik}
-{\vec{\xi}}_j{\vec{\xi}}_k g^N_{il}]
\\& -{\nabla}^2_{ik}(\vec{\xi}\cdot\vec{x})g^N_{jl}
-{\nabla}^2_{jl}(\vec{\xi}\cdot\vec{x})g^N_{ik} +
{\nabla}^2_{jk}(\vec{\xi}\cdot\vec{x})g_{il}^N +
{\nabla}^2_{il}(\vec{\xi}\cdot\vec{x})g^N_{jk}
\\&+|\vec{\xi}|^2g^N_{il}g^N_{jk}- |\vec{\xi}|^2g^N_{ik}
g^N_{lj}]
\}
\end{split}
\end{equation}

\par Hence, in order to find the transformation laws for the covariant
derivatives of the full curvature tensor, the Ricci curvature
tensor and of the factors $\nabla^p{\psi}_h$, we will need the
transformation law for the Levi-Civita
 connection in the case at hand:

\begin{equation}
\label{levicivitaxi} ({\nabla}_k {\eta}_l)^{\hat{g}^N}(x)=
({\nabla}_k {\eta}_l)^{g^N}- {\vec{\xi}}_k {\eta}_l -{\vec{\xi}}_l
{\eta}_k + {\vec{\xi}}^s {\eta}_s g^N_{kl}
\end{equation}

 \par These relations show that in (\ref{bignint}), under the re-scaling
$g^N(x) \longrightarrow \hat{g}^N(x)= e^{2\vec{\xi} \cdot\vec{x}}
g^N(x)$, the integrand $I^Z_{g^{N}}({\psi}_1,\dots
 ,{\psi}_Z)(x)$ transforms as follows:

\begin{equation}
\label{itransxi}\begin{split} &I^Z_{\hat{g}^N}({\psi}_1,\dots
 ,{\psi}_Z)(x)= e^{-n\vec{\xi}\cdot \vec{x}} [I^Z_{g^{N}}({\psi}_1,\dots
 ,{\psi}_Z)(x)+\\& S^Z_{g^{N}}({\psi}_1,\dots
 ,{\psi}_Z,\vec{\xi})(x)]
\end{split}
\end{equation}

\par  Where $S^Z_{g^{N}}({\psi}_1,\dots
 ,{\psi}_Z,\vec{\xi})(x)$ is obtained by applying the transformation laws
described above to each factor in every complete contraction in
$I^Z_{g^{N}}({\psi}_1,\dots
 ,{\psi}_Z)(x)$. $S^Z_{g^{N}}({\psi}_1,\dots
 ,{\psi}_Z,\vec{\xi})(x)$ will be a linear combination of complete contractions, each of which
 depends on $\vec{\xi}$. Hence equation (\ref{bignint}) will give, for the
metric $(g^N)^{\vec{\xi}}(x)$:

\begin{equation}
\label{bignintxi} {\int}_{\mathbb{R}^N}e^{(N-n)\vec{\xi}\cdot
\vec{x}} [I^Z_{g^N}({\psi}_1,\dots
 ,{\psi}_Z)+ S^Z_{g^N}({\psi}_1,\dots
 ,{\psi}_Z,\vec{\xi})] dV_{g^n}=0
\end{equation}

\par Roughly speaking, our goal for this subsection will be to
perform integrations by parts for the complete contractions in
$S^Z_{g^{N}}({\psi}_1,\dots
 ,{\psi}_Z,\vec{\xi})(x)$ in order to reduce equation
 (\ref{bignint}) to the form:

  $${\int}_{\mathbb{R}^N}e^{(N-n)\vec{\xi}\cdot\vec{x}}[I^Z_{g^{N}}({\psi}_1,\dots
 ,{\psi}_Z)+ L^Z_{g^{N}}({\psi}_1,\dots
 ,{\psi}_Z)]dV_{g^n}=0$$

 where $L^Z_{g^{N}}({\psi}_1,\dots
 ,{\psi}_Z)(x)$ will be independent of $\vec{\xi}$. This will be
 done and explained in rigor below. Keeping this vaguely outlined
 strategy in mind, let us note the identity:

 \begin{equation}
 \label{oneder}
{\nabla}_s( e^{(N-n)\vec{\xi}\cdot\vec{x}})=(N-n){\vec{\xi}}_s (
e^{(N-n)\vec{\xi}\cdot\vec{x}})
 \end{equation}

\par More generally, let us denote by ${\partial}^m_{s_1\dots s_k}$
the coordinate derivative with respect to our coordinate system.
Then, for $k>1$ we have:

\begin{equation}
\label{manyder} {\partial}^k_{s_1\dots s_k}(\vec{\xi}\cdot\vec{x})=0
\end{equation}

for every $x\in B(x_0,\epsilon)$.

 \par Let us consider
 the Christoffel symbols ${\Gamma}^k_{ij}$ with respect to our arbitrary
 coordinate system. Let
  $S{\nabla}^m_{s_1\dots s_m} \vec{\xi}_j$ stand for ${\nabla}^m_{(s_1\dots s_m}
\vec{\xi}_{j)}$ and  $S{\nabla}^p_{r_1\dots r_p}{\Gamma}^k_{ij}$
stand for ${\nabla}^p_{(r_1\dots r_p} {\Gamma}^k_{ij)}$.
\newline

\par Let us now write $I^Z_{g^{N}}({\psi}_1,\dots
 ,{\psi}_Z)$ as a linear combination of complete contractions
 in the following form:

 \begin{equation}
 \label{linisym}
\begin{split}
&contr({\nabla}_{r_1\dots r_{m_1}}^{m_1}R_{i_1j_1k_1l_1}\otimes
\dots \otimes {\nabla}_{v_1\dots
v_{m_s}}^{m_s}R_{i_sj_sk_sl_s}\otimes \\& {\nabla}_{t_1 \dots
t_{p_1}}^{p_1}Ric_{{\alpha}_1 {\beta}_1} \otimes \dots \otimes
{\nabla}_{z_1 \dots z_{p_q}}^{p_q} Ric_{{\alpha}_q{\beta}_q}
\otimes{\nabla}^{{\nu}_1}_ {{\chi}_1\dots
{\chi}_{{\nu}_1}}{\psi}_1 \otimes\dots \otimes
{\nabla}^{{\nu}_Z}_{{\omega}_1\dots {\omega}_{{\nu}_Z}}
 {\psi}_Z)
\end{split}
 \end{equation}

\par Where each of the factors ${\nabla}_{r_1\dots
r_{m_1}}^{m_1}R_{i_1j_1k_1l_1}, \dots ,{\nabla}_{v_1\dots
v_{m_s}}^{m_s}R_{i_sj_sk_sl_s}$ has no two of the indices
 $i,j,k,l$ contracting against each other in (\ref{linisym}).

\par Now, in dimension $N$, we can apply the identities (\ref{curvxi}),
 (\ref{riccixi}),(\ref{scalarxi}) and (\ref{levicivitaxi}) to
find the form of $S^Z_{g^N}({\psi}_1,\dots
 ,{\psi}_Z,\vec{\xi})(x)$. In particular, we write
 it as a linear combination of complete contractions in the
 following two forms:

\begin{equation}
 \label{linisymxi1}
\begin{split}
&contr({\nabla}_{r_1\dots r_{m_1}}^{m_1}R_{i_1j_1k_1l_1}\otimes
\dots \otimes {\nabla}_{v_1\dots
v_{m_s}}^{m_s}R_{i_sj_sk_sl_s}\otimes \\& {\nabla}_{t_1 \dots
t_{p_1}}^{p_1}Ric_{{\alpha}_1 {\beta}_1} \otimes \dots \otimes
{\nabla}_{z_1 \dots z_{p_q}}^{p_q}
Ric_{{\alpha}_q{\beta}_q}\otimes {\nabla}^{{\nu}_1}_
{{\chi}_1\dots {\chi}_{{\nu}_1}}{\psi}_1 \otimes\dots \otimes
{\nabla}^{{\nu}_Z}_{{\omega}_1\dots {\omega}_{{\nu}_Z}} {\psi}_Z
\\& \otimes \vec{\xi}\otimes\dots \otimes \vec{\xi})
\end{split}
 \end{equation}

\begin{equation}
 \label{linisymxi2}
\begin{split}
&contr({\nabla}_{r_1\dots r_{m_1}}^{m_1}R_{i_1j_1k_1l_1}\otimes
\dots \otimes {\nabla}_{v_1\dots
v_{m_s}}^{m_s}R_{i_sj_sk_sl_s}\otimes \\& {\nabla}_{t_1 \dots
t_{p_1}}^{p_1}Ric_{{\alpha}_1 {\beta}_1} \otimes \dots \otimes
{\nabla}_{z_1 \dots z_{p_q}}^{p_q} Ric_{{\alpha}_q{\beta}_q}
\otimes{\nabla}^{{\nu}_1}_ {{\chi}_1\dots
{\chi}_{{\nu}_1}}{\psi}_1 \otimes\dots \otimes
{\nabla}^{{\nu}_Z}_{{\omega}_1\dots {\omega}_{{\nu}_Z}} {\psi}_Z
\\& \otimes \vec{\xi}\otimes\dots \otimes \vec{\xi}\otimes
S[{\nabla}^{w_1}_{u_1\dots u_{w_1}}\vec{\xi}]\otimes \dots
\otimes S[{\nabla}^{w_l}_{q_1\dots q_{w_l}}\vec{\xi}] )
\end{split}
\end{equation}

where each $w_a\ge 1$. We also let $k$ stand for the number
 of factors $\vec{\xi}$ and $l$ for the number of factors
$S{\nabla}^w\vec{\xi}$.

\par We will call complete contractions in the above two forms
 $\vec{\xi}$-contractions.

\par In order to see that we can indeed write
$S^Z_{g^N}({\psi}_1,\dots ,{\psi}_Z)$ as a linear combination of
complete contractions in the above form, we only need the
equation:

\begin{equation}
\label{symunsymgenxi}
\begin{split}
&{\nabla}_aS{\nabla}^m_{r_1\dots r_m}\vec{\xi}_j=
S{\nabla}^m_{ar_1\dots r_m}\vec{\xi}_j +
C_{m-1}\cdot
S^{*}{\nabla}^{m-1}_{r_1\dots r_{m-1}}
R_{aijd}\vec{\xi}^d +
\\&
{\Sigma}_{u\in U^m} a_u
pcontr({\nabla}^{m'} R_{abcd} S{\nabla}^{s_u}\vec{\xi})
\end{split}
\end{equation}
where $S^{*}{\nabla}^{m-1}_{r_1\dots r_{m-1}}
R_{aijd}$ stands for the symmetrization of
${\nabla}^{m-1}_{r_1\dots r_{m-1}}R_{aijd}$ over the indices
 $r_1,\dots ,r_{m-1},i$ and
the symbol $pcontr({\nabla}^{m'} R_{abcd}
S{\nabla}^{s_u}\vec{\xi})$ stands for a partial contraction of at
least one factor ${\nabla}^mR_{ajkl}$ (to one of which the index
$a$ belongs) against a factor $S{\nabla}^{s_u}\vec{\xi}$ with
$s_u\ge 1$.
\newline

\par Our next goal is to answer the following: Given a fixed linear
combination $I^Z_{g^n}({\psi}_1,\dots
 ,{\psi}_Z)(x)$ and its rewriting
$I^Z_{g^N}({\psi}_1,\dots ,{\psi}_Z)(x)$ in any dimension $N\ge n$, how
 does $S^Z_{g^{N}}({\psi}_1,\dots
 ,{\psi}_Z,\vec{\xi})(x)$ depend upon the dimension $N$?
\newline

\par In order
 to answer this question, we will introduce certain definitions. Let us
 for this purpose treat the function $\vec{\xi}\cdot\vec{x}$ as a
 function $\omega (x)$. Hence
$\vec{\xi}_i={\nabla}_i(\vec{\xi}\cdot\vec{x})$ and we can
 speak of the re-writing of a $\vec{\xi}$-contraction in dimension $N$.
 We will consider
 the complete contraction
$C^i_{g^n}({\psi}_1,\dots ,{\psi}_Z, \vec{\xi})$
 together with its rewriting
$C^i_{g^N}({\psi}_1,\dots ,{\psi}_Z,
 \vec{\xi})$ in every dimension $N\ge n$ and call this sequence
  a dimension-independent complete contraction.

 \par On the other hand, we define:

\begin{definition}
\label{xi}
Any factor of the form $\vec{\xi}$ or of the form
$S{\nabla}^m\vec{\xi}$, $m\ge 1$ will be called a
$\vec{\xi}$-factor.
\end{definition}

\begin{definition}
\label{dimdep}
 Consider a sequence $\{ C_{(g,N)}({\psi}_1,\dots
 ,{\psi}_Z,\vec{\xi})\}$ of complete contractions times coefficients in dimensions
$N=n, n+1,\dots$ where the following formula
 holds: There is a fixed complete contraction, say
$C_{g^n}({\psi}_1,\dots ,{\psi}_Z,
 \vec{\xi})$ and a fixed rational function $Q(N)$ so that:

 $$C_{(g,N)}({\psi}_1,\dots
 ,{\psi}_Z,\vec{\xi})= Q(N)\cdot C_{g^N}({\psi}_1,\dots ,{\psi}_Z, \vec{\xi})$$

 where $C_{g^N}({\psi}_1,\dots ,{\psi}_Z,
 \vec{\xi})$ is the rewriting of $C_{g^n}({\psi}_1,\dots ,{\psi}_Z, \vec{\xi})$
  in dimension $N$. In that case, we will say that we
 have a dimension-dependent $\vec{\xi}$-contraction. Furthermore, we
 will say that the three defining numbers of $C_{(g,N)}
({\psi}_1,\dots
 ,{\psi}_Z,\vec{\xi})$ are $(d,k,l)$ where $d$ is the degree
 of the rational function $Q(N)$, $k$ is the number of factors $\vec{\xi}$ and
 $l$ is the number of factors
 $S{\nabla}^m_{i_1\dots i_m}\vec{\xi}_a$, $m\ge 1$.

\par (Given a rational function $Q(N)=\frac{P(N)}{L(N)}$, we
 define the degree of $Q(N)$, $deg[Q(N)]=deg[P(N)]-
deg[L(N)]$. We also define the leading order coefficient of
$Q(N)$ to be $\frac{a_P}{a_L}$, where $a_P$ is the leading
 order coefficient of $P(N)$ and $a_L$ is the leading order
 coefficient of $L(N)$).
\end{definition}

 \par Given a fixed set of numbers $\{a_i\},
 i\in I$, and a set of dimension-dependent
$\vec{\xi}$-contractions
$C^i_{(g,N)}({\psi}_1,\dots ,{\psi}_Z,
 \vec{\xi})$,  we
 can form  in each dimension $N\ge n$ the linear combination:

 $$L_{g^N}({\psi}_1,\dots ,{\psi}_Z,
 \vec{\xi})= {\Sigma}_{i\in I} a_i C^i_{(g,N)}({\psi}_1,\dots ,{\psi}_Z,
 \vec{\xi})$$

 \par Hence we obtain in this way a sequence of
 linear combinations, where the index set for the sequence
 is the set ${\bf N}=\{n,n+1,\dots \}$.

\begin{definition}
\label{lindimdep}
 We will say that a sequence of linear combinations as above
 is dimension-dependent.

  \par We will say that a sequence of linear
 combinations as above is
 {\it suitable} if for each of the $\vec{\xi}$-contractions
 $C^i_{(g,N)}({\psi}_1,\dots ,{\psi}_Z,
 \vec{\xi})$ we have that its three defining numbers satisfy:
 $k+l\ge d$.
\end{definition}

We then have:

 \begin{lemma}
 \label{dimdepindep}
 $S^Z_{g^N}({\psi}_1,\dots
 ,{\psi}_Z,\vec{\xi})(x)$ is a suitable linear combination of
 $\vec{\xi}$-contractions of the form (\ref{linisymxi1}) and
 (\ref{linisymxi2}), with $k+l\ge d\ge 1$.
 \end{lemma}

{\it Proof:} We write

$$I^Z_{g^N}({\psi}_1,\dots ,{\psi}_Z)={\Sigma}_{i\in I} a_i
C^i_{g^N}({\psi}_1,\dots ,{\psi}_Z)$$

where each $C^i_{g^N}({\psi}_1,\dots ,{\psi}_Z)$ is in the
 form (\ref{linicontrpsi}) and has weight $-n$.

\par We introduce some further terminology:
We call the tensors
$({\nabla}^m_{r_1\dots r_m}R_{ijkl})^{g^N}$,
$({\nabla}^p_{r_1\dots r_p} {\psi}_l)^{g^N}$,
$(S{\nabla}^m_{r_1\dots r_m}\vec{\xi}_a)^{g^N}$,
$\vec{\xi}_i$ and $g^N_{ij}$ the {\it free} tensors. We call
partial contractions of those tensors the
 {\it extended} free tensors. (Recall that a {\it partial
 contraction} means a tensor product with some pairs of indices
 contracting against each other.)

\par We see that $e^{-2\vec{\xi}\cdot\vec{x}}
({\nabla}^m_{r_1\dots r_m} R_{ijkl})^{\hat{g}^N}$,
$({\nabla}^p_{r_1\dots r_p} {\psi}_l)^{\hat{g}^N}$ can be written
as linear combinations of extended free
 tensors, after applying the identity (\ref{symunsymgenxi}),
 if necessary.

\par Now, consider any complete contraction
$C^i_{g^N}({\psi}_1,\dots ,{\psi}_Z)$ (in the form
(\ref{linicontrpsi})) in
$I^Z_{g^N}({\psi}_1,\dots ,{\psi}_Z)$ and do the following:
For each of its factors
${\nabla}^m_{r_1\dots r_m}R_{ijkl}$ or
${\nabla}^p_{r_1\dots r_p}{\psi}_l$,
 consider:

$$e^{-2\vec{\xi}\cdot \vec{x}}
({\nabla}^m_{r_1\dots r_m}R_{ijkl})^{\hat{g}^N} = {\Sigma}_{j\in
J'} a_j T^j_{r_1\dots r_mijkl}$$

$$({\nabla}^p_{r_1\dots r_p}{\psi}_l)^{\hat{g}^N}=
{\Sigma}_{j\in J} a_j T^j_{r_1\dots r_p}$$

 where each $T^r_{i_1\dots i_s}$
is an extended free tensor.  We then substitute each factor
${\nabla}^m_{r_1\dots r_m}R_{ijkl}$ by one
$e^{2\vec{\xi}\cdot\vec{x}}a_j T^j_{r_1\dots r_mijkl}$ and each
 factor ${\nabla}^p_{r_1\dots r_p}{\psi}_l$ by one
$a_j T^j_{r_1\dots r_p}$. After this substitution, we perform the
same contractions of indices as in $C^i_{g^N}({\psi}_1,\dots
,{\psi}_Z)$, with respect to the metric $(g^N)$. We do this
according to the following algorithm: Suppose we are contracting
two indices $\alpha,\beta$. If none of them belongs to a tensor
$g^N_{ij}$, we just do that particular contraction. If $\alpha$
but not $\beta$ belongs to a factor
 $g^N_{\alpha\gamma}$, we cross out the index $\beta$ in the other factor
 and substitute it by $\gamma$, and then omit the $g^N_{\alpha\gamma}$.
 Finally, if both the indices $\alpha,\beta$ belong to the same factor
$g^N_{\alpha\beta}$, we cross out that factor and bring out
 a factor of $N$.

\par Adding over all those substitutions, we obtain
$e^{n\vec{\xi}\cdot\vec{x}} C^i_{\hat{g}^N} ({\psi}_1,\dots
,{\psi}_Z)$.

\par We observe that
$e^{n\vec{\xi}\cdot\vec{x}}C^i_{\hat{g}^N}({\psi}_1,\dots
,{\psi}_Z)$ is a dimension-dependent linear
 combination. It follows that $S^Z_{g^N}({\psi}_1,\dots
,{\psi}_Z,\vec{\xi})$ is a dimension-dependent linear
 combination, in the form:

$$S^Z_{g^N}({\psi}_1,\dots ,{\psi}_Z,\vec{\xi})=
{\Sigma}_{l\in L} N^{b_l} C^l_{g^N}({\psi}_1,\dots ,
{\psi}_Z,\vec{\xi})$$
where each complete contraction $C^l_{g^N}({\psi}_1,\dots ,
{\psi}_Z,\vec{\xi})$ is in the form:

\begin{equation}
\label{arxid}
\begin{split}
&contr({\nabla}_{r_1\dots r_{m_1}}^{m_1}
R_{i_1j_1k_1l_1}\otimes
\dots \otimes {\nabla}_{v_1\dots
v_{m_s}}^{m_s}R_{i_sj_sk_sl_s}\otimes
\\&{\nabla}^{{\nu}_1}_ {{\chi}_1\dots
{\chi}_{{\nu}_1}}{\psi}_1 \otimes\dots \otimes
{\nabla}^{{\nu}_Z}_{{\omega}_1\dots {\omega}_{{\nu}_Z}}
{\psi}_Z
 \otimes \vec{\xi}\otimes\dots \otimes \vec{\xi}\otimes
S{\nabla}^{w_1}_{u_1\dots u_{w_1}}\vec{\xi}\otimes \dots
\otimes S{\nabla}^{w_l}_{q_1\dots q_{w_l}}\vec{\xi})
\end{split}
\end{equation}
where $l\ge 0$ and the factors ${\nabla}^mR_{ijkl}$ are allowed
to have internal contractions.

\par In other words,
$C^l_{g^N}({\psi}_1,\dots , {\psi}_Z,\vec{\xi})$
is a linear combination of complete contractions in the form
(\ref{linisymxi1}) or (\ref{linisymxi2}).
\newline

\par Therefore, what remains to be checked is that each
 dimension-dependent $\vec{\xi}$-contraction
$N^{b_i}C^i_{g^N}({\psi}_1,\dots ,{\psi}_Z,\vec{\xi})$ in
$S^Z_{g^N}({\psi}_1,\dots ,{\psi}_Z,\vec{\xi})$ with
$|\vec{\xi}|\ge 1$  satisfies the identity $k+l\ge d$.

\par In order to see this, let us consider any summand
$T^j_{r_1\dots r_mijkl}$ or $T^j_{r_1\dots r_p}$ and denote by
$|g|$ the number of its factors $g^N_{ij}$ and by $|\vec{\xi}|$
the number of its factors $\vec{\xi}_i$ or $S{\nabla}^m_{r_1\dots
r_m}\vec{\xi}_a$. It follows, from
 identities (\ref{curvxi}) and (\ref{levicivitaxi}) that for
 each $T^j_{r_1\dots r_mijkl}$ or $T^j_{r_1\dots r_p}$
 we have that $|\vec{\xi}|\ge |g|$.

\par By virtue of that inequality, the formula (\ref{symunsymgenxi})
(which shows us that if we write a complete contraction in the form
(\ref{arxid}) as a linear combination of complete contractions in the
 forms (\ref{linisymxi1}), (\ref{linisymxi2}) the number of
$\vec{\xi}$-factors remains unaltered) and the algorithm outlined
 above, we observe that for each dimension-dependent
$\vec{\xi}$-contraction
$N^{b_i}C^i_{g^N}({\psi}_1,\dots ,{\psi}_Z,\vec{\xi})$,
 we will have that $b_i$ is less than or equal to the number
 of factors $\vec{\xi}$ or $S{\nabla}^m\vec{\xi}$ in
$C^i_{g^N}({\psi}_1,\dots ,{\psi}_Z,\vec{\xi})$. $\Box$

\begin{definition}
\label{descendants}
Consider any complete contraction
$C^Z_{g^N}({\psi}_1,\dots ,{\psi}_Z)$, in the form
(\ref{linisym}). Consider the quantity:

$$e^{n\vec{\xi}\cdot\vec{x}}C^Z_{e^{2\vec{\xi}\cdot\vec{x}}g^N}
({\psi}_1,\dots ,{\psi}_Z)(x)-
C^Z_{g^N}({\psi}_1,\dots ,{\psi}_Z)(x)$$

The above quantity can be computed by applying the identities
(\ref{curvxi}), (\ref{riccixi}) (\ref{scalarxi}),
(\ref{levicivitaxi}), (\ref{symunsymgenxi}) to each factor
 in $C^Z_{g^N}({\psi}_1,\dots ,{\psi}_Z)$. We write:

$$e^{n\vec{\xi}\cdot\vec{x}}C^Z_{e^{2\vec{\xi}\cdot\vec{x}}g^N}
({\psi}_1,\dots ,{\psi}_Z)(x)-
C^Z_{g^N}({\psi}_1,\dots ,{\psi}_Z)(x)=
{\Sigma}_{t\in T} a_t N^{b_t}
C^t_{g^N}({\psi}_1,\dots ,{\psi}_Z,\vec{\xi})(x)$$

\par Where each dimension-dependent $\vec{\xi}$-contraction
$N^{b_t} C^t_{g^N}({\psi}_1,\dots ,{\psi}_Z,\vec{\xi})(x)$
satisfies $k+l\ge b_t$. Here $C^t_{g^N}({\psi}_1,\dots
,{\psi}_Z,\vec{\xi})(x)$ stands for the rewriting of
\\ $C^t_{g^n}({\psi}_1,\dots ,{\psi}_Z,\vec{\xi})(x)$ in dimension
$N$.

There are many expressions as above for
$e^{n\vec{\xi}\cdot\vec{x}}C^Z_{e^{2\vec{\xi}\cdot\vec{x}}g^N}
({\psi}_1,\dots ,{\psi}_Z)(x)- C^Z_{g^N}({\psi}_1,\dots
,{\psi}_Z)(x)$ which are equal by substitution but not identical.
Once we pick one such expression,
 we will call each dimension-dependent
$\vec{\xi}$-contraction
$N^{b_t} C^t_{g^N}({\psi}_1,\dots ,{\psi}_Z,\vec{\xi})(x)$
 a descendant
 of $C^Z_{g^N}({\psi}_1,\dots ,{\psi}_Z)(x)$.
\end{definition}

 \par We are now near the point where we can integrate by parts
  in the relation (\ref{bignintxi}). At this stage,
 we want to distinguish between descendants of the
 complete contractions
 in $I^Z_{g^N}({\psi}_1,\dots ,{\psi}_Z)$. We define:

 \begin{definition}
 \label{goodbad}
For any complete contraction $C^i_{g^N}({\psi}_1,\dots ,{\psi}_Z)$
in
\\ $I^Z_{g^N}({\psi}_1,\dots ,{\psi}_Z)$, we
will call one of its descendants easy
if $d< l+k$.

\par A descendant in the form (\ref{linisymxi1}) will be called
 good if $d=k>0$ and $l=0$.
A descendant in the form (\ref{linisymxi2}) will be called
undecided if $d=k+l$ and $k,l>0$. (That is, it
contains at least one factor of the form
$S{\nabla}^p\vec{\xi}_i$ with $p\ge 1$
and at least one factor of the form $\vec{\xi}$).

\par Finally, a descendent in the form (\ref{linisymxi2})
with $d=k+l$ will be called hard if $k=0, l>0$. (That is, if
all its $\vec{\xi}$-factors are all of the form
 $S{\nabla}^m\vec{\xi}_j$, with $m\ge 1$.
\end{definition}

\par We now have the equation (\ref{bignintxi}) in any dimension
$N$ and we have $S^Z_{g^N}({\psi}_1,\dots ,{\psi}_Z,
\vec{\xi})$
written out as a linear combination of good, easy,
   undecided and hard complete contractions.

\subsection{The integrations by parts for
$S^Z_{g^N}({\psi}_1,\dots ,{\psi}_Z,\vec{\xi})$.}

\par We want to perform integrations by parts in equation (\ref{bignintxi}).
We will treat the four cases above separately.

\par Let us first treat the easy $\vec{\xi}$-contractions.
Using (\ref{kormaki}), we write out each factor of the form
$S{\nabla}^m\vec{\xi}_s$ as a linear combination of partial
contractions of the Christoffel symbols and their derivatives
(with respect to our arbitrarily chosen coordinate system) and
also of the vector $\vec{\xi}$. We also write out each of the
tensors ${\nabla}^mR_{ijkl}$ as a linear combination of partial
contractions of Christoffel symbols and their derivatives. Hence,
given an {\it easy} $\vec{\xi}$-contraction $P(N)\cdot
C^Z_{g^N}(\psi_1,\dots , \psi_Z,\vec{\xi})$, we express it in our
coordinate system as:

\begin{equation}
\label{chrisana} contr({\partial}^{m_1}{\Gamma}^k_{ij}\otimes\dots
\otimes {\partial}^{m_s}{\Gamma}^k_{ij}\otimes{\nabla}^{p_1}{\psi}_1\otimes\dots \otimes {\nabla}^{p_Z}{\psi}_Z \otimes
\vec{\xi}\otimes\dots
\otimes  \vec{\xi})
\end{equation}

Hence we will have the following identity:

\begin{equation}
\label{easyana} \begin{split} &{\int}_{\mathbb{R}^N}
e^{(N-n)\vec{\xi}\cdot\vec{x}}P(N)\cdot C^Z_{g^N}({\psi}_1,\dots
,{\psi}_Z,\vec{\xi})dV_{g^n}=
\\&{\int}_{\mathbb{R}^N}
e^{(N-n)\vec{\xi}\cdot\vec{x}}P(N)\cdot {\Sigma}_{l\in L}a_l
Contr_l({\partial}^m\Gamma , {\psi}_1,\dots ,{\psi}_Z,
\vec{\xi})dV_{g^n}
\end{split}
\end{equation}

 Where in the above equation the degree of the polynomial
$P(N)$ is
 strictly less than the number of factors $\vec{\xi}$ in the  contraction
 $Contr_l({\partial}^m\Gamma , {\psi}_1,\dots ,{\psi}_Z,
\vec{\xi})$. Now, use the identity (\ref{oneder}) in order to
substitute one factor $\vec{\xi}_i$ in the complete contraction by
the factor $\frac{{\nabla}_ie^{(N-n)\vec{\xi}\cdot\vec{x}}}
{N-n}$. We then integrate by parts with respect to the
 derivative ${\nabla}_i$. Let us note here that this
 integration by parts
is with respect to the Riemannian connection ${\nabla}_i$.

What we will get is the following:

\begin{equation}
\label{easycon}
\begin{split}
&{\int}_{\mathbb{R}^N}
e^{(N-n)\vec{\xi}\cdot\vec{x}}P(N)\cdot C^Z_{g^N}({\psi}_1,\dots
,{\psi}_Z,\vec{\xi})dV_{g^n}=
 \\&-{\int}_{\mathbb{R}^N}
e^{(N-n)\vec{\xi}\cdot\vec{x}} \frac{P(N)}{(N-n)} {\Sigma}_{k\in K}
a_k Contr_k({\partial}^m\Gamma ,
{\psi}_1,\dots,{\psi}_Z,\vec{\xi}) dV_{g^n}
\end{split}
\end{equation}

 Each complete contraction
$Contr_k({\partial}^m\Gamma , {\psi}_1,\dots ,{\psi}_Z\vec{\xi})$
will be in the form (\ref{chrisana}). Also, the number of factors
$\vec{\xi}$ in each contraction
$Contr_k({\partial}^m\Gamma , {\psi}_1,\dots ,{\psi}_Z,\vec{\xi})$
is by one less than the number of such factors in the complete
contraction $C^Z_{g^N}({\psi}_1,\dots ,{\psi}_Z,\vec{\xi})$.
 Hence, inductively repeating the above procedure we will obtain:

\begin{equation}
\label{easycon} \begin{split} & {\int}_{\mathbb{R}^N}
e^{(N-n)\vec{\xi}\cdot\vec{x}}P(N)\cdot C^Z_{g^N}({\psi}_1,\dots
,{\psi}_Z,\vec{\xi})dV_{g^n}=
\\&{\int}_{\mathbb{R}^N}
e^{(N-n)\vec{\xi}\cdot\vec{x}} \frac{P(N)}{(N-n)^w} {\Sigma}_{k\in K}
a_k Contr_k({\partial}^m\Gamma , {\psi}_1,\dots ,{\psi}_Z)
dV_{g^n}
\end{split}
\end{equation}

\par Where we will have $deg[P(N)]=d<w$.
\newline

{\it The good $\vec{\xi}$-contractions.}
\newline

\par Let us now deal with the good complete contractions in
$S^Z_{g^N}({\psi}_1,\dots ,{\psi}_Z,\vec{\xi})$. In this case it
is useful not to write things out in terms of Christoffel symbols
but to work intrinsically on the Riemannian manifold. We have a
good $\vec{\xi}$-contraction $P(N)\cdot C^Z_{g^N}({\psi}_1,\dots
,{\psi}_Z,\vec{\xi})$ in the form (\ref{linisymxi1}) and we want
to perform integrations by
 parts in the integral:

\begin{equation}
\label{goodana}
{\int}_{\mathbb{R}^N}e^{(N-n)\vec{\xi}\cdot\vec{x}}P(N)\cdot C^Z_{g^N}({\psi}_1,\dots
,{\psi}_Z,\vec{\xi})dV_{g^n}
\end{equation}

 We will again use the
 identity (\ref{oneder}). Let us arbitrarily pick out one of
 the $k=d$ factors $\vec{\xi}$ in $C^Z_{g^N}({\psi}_1,\dots
,{\psi}_Z,\vec{\xi})$. Now, use the identity (\ref{oneder}) in
order to substitute the factor $\vec{\xi}_i$ in the complete
contraction by the factor
$\frac{{\nabla}_i[e^{(N-n)\vec{\xi}\cdot\vec{x}}]}{N-n}$. We then
integrate by parts with respect to the derivative ${\nabla}_i$.
Let us again note that this integration by parts is with
respect to the Riemannian connection ${\nabla}_i$.

\par Now, if the $\vec{\xi}$-contraction $C_{g^N}({\psi}_1,\dots
,{\psi}_Z,\vec{\xi})$ in the form (\ref{linisymxi1}) has $L$
factors (including the $k$ factors $\vec{\xi}$), the integration
by parts will produce a sum of $L-1$ complete contractions.
Explicitly, we will have:

\begin{equation}
\label{expostint}
\begin{split}
&{\int}_{\mathbb{R}^N}e^{(N-n)\vec{\xi}\cdot\vec{x}}P(N)\cdot C_{g^N}({\psi}_1,\dots
,{\psi}_Z,\vec{\xi})dV_{g^n}=
\\& -{\int}_{\mathbb{R}^N}e^{(N-n)\vec{\xi}\cdot\vec{x}}\frac{P(N)}{N-n}\cdot
{\Sigma}_{\alpha =1}^{L-1} C^{\alpha}_{g^N} ({\psi}_1,\dots
,{\psi}_Z,\vec{\xi})dV_{g^n}
\end{split}
\end{equation}

 Let us separate these $\vec{\xi}$-contractions
$C^{\alpha}_{g^N} ({\psi}_1,\dots
,{\psi}_Z,\vec{\xi})$ into two categories: A
$\vec{\xi}$-contraction belongs to the first category if the
derivative ${\nabla}_i$ has hit one of the factors
${\nabla}^mR_{ijkl}$, ${\nabla}^pRic$ or
${\nabla}^p{\psi}_k$. Hence,
 we see that  $\frac{P(N)}{N-n}\cdot
C^{\alpha}_{g^N} ({\psi}_1,\dots ,{\psi}_Z,\vec{\xi})$
is a linear combination of
$\vec{\xi}$-contractions in the form (\ref{linisymxi1}) with $k-1$ factors $\vec{\xi}$. If $k=1$, each will be in the form
 (\ref{linisym}). Otherwise, each of them will be a good
$\vec{\xi}$-contraction.

\par On the other hand, a $\vec{\xi}$-contraction
$\frac{P(N)}{N-n}\cdot
C^{\alpha}_{g^N} ({\psi}_1,\dots ,{\psi}_Z,\vec{\xi})$
belongs to the
second category if the derivative ${\nabla}_i$ hit one of the
$k-1$ factors $\vec{\xi}$. In that case, we get a
$\vec{\xi}$-contraction in the form (\ref{linisymxi2}) with
 $k-2$ factors
$\vec{\xi}$ and one factor ${\nabla}_i\vec{\xi}$. It will be an undecided or a hard $\vec{\xi}$-contraction.

\par Now, we can repeat the above intrinsic integration by parts
for each of the good $\vec{\xi}$-contractions
$\frac{P(N)}{N-n}\cdot
C^{\alpha}_{g^N} ({\psi}_1,\dots ,{\psi}_Z,\vec{\xi})$. They are each of the form (\ref{linisymxi1}) with
$k-1$ factors $\vec{\xi}$. Each of these integrations by parts
will give a sum of $\vec{\xi}$-contractions, $L-k+1$ of which
are in
 the form (\ref{linisymxi1}) with $k-2$ factors
$\vec{\xi}$ and $k-2$ of them will be of the form
(\ref{linisymxi2}) (either undecided or hard).
 Hence, we can form a procedure of $k$ steps,
starting from $C^Z_{g^N}({\psi}_1,\dots ,{\psi}_Z,\vec{\xi})$ and
integrating by parts one factor $\vec{\xi}$ at a time. At each
stage we get a sum of good and of undecided or hard
$\vec{\xi}$-contractions out of this integration by parts.
 We then focus on the good $\vec{\xi}$-contractions that we have
 obtained and we repeat the integration by parts. Thus, after this
sequence of integrations by parts we will have:

\begin{equation}
\label{goodana}
\begin{split}
&{\int}_{\mathbb{R}^N}e^{(N-n)\vec{\xi}\cdot\vec{x}}P(N)\cdot C^Z_{g^N}({\psi}_1,\dots
,{\psi}_Z,\vec{\xi})dV_{g^n}=
\\&{\int}_{\mathbb{R}^N}e^{(N-n)\vec{\xi}\cdot\vec{x}}[
\frac{P(N)}{(N-n)^k}
{\Sigma}_{j\in J}a_j C^j_{g^N}({\psi}_1,\dots ,{\psi}_Z)+
\\&{\Sigma}_{h\in H}\frac{P_h(N)}{(N-n)^{s_h}}C^h_{g^N}({\psi}_1,\dots ,{\psi}_Z,\vec{\xi})]dV_{g^n}
\end{split}
\end{equation}

where the complete contractions $C^j_{g^N}({\psi}_1,\dots
,{\psi}_Z)$ are in the form (\ref{linisym}) (they are
independent of the variable $\vec{\xi}$)
and the $\vec{\xi}$-contractions
$C^h_{g^N}({\psi}_1,\dots ,{\psi}_Z,\vec{\xi})$ are
in the form (\ref{linisymxi2}) and are undecided or hard.
 Each of the
undecided $\vec{\xi}$-contractions will have at most $k-1$
$\vec{\xi}$-factors. For each complete contraction
$C^j_{g^N}({\psi}_1,\dots ,{\psi}_Z)$ we have that
$deg[P(N)]=k$.
For each complete contraction $C^h_{g^N}({\psi}_1,\dots
,{\psi}_Z,\vec{\xi})$, with $l_h$ factors ${\nabla}\vec{\xi}$ and $k_h$ factors $\vec{\xi}$,
 we have $k_h+ l_h+s_h=deg[P_h(N)]$.
\newline

{\it The undecided $\vec{\xi}$-contractions.}
\newline

\par We now proceed to integrate by parts the undecided
$\vec{\xi}$-contractions. Let
\\ $C^Z_{g^N}({\psi}_1,\dots
,{\psi}_Z,\vec{\xi})$ be  an undecided $\vec{\xi}$-contraction in
the form (\ref{linisymxi2}).
 We want to perform integrations by parts in the integral:

$${\int}_{R^N}e^{(N-n)\vec{\xi}\cdot\vec{x}}
\frac{P_h(N)}{(N-n)^{m_h}}C^Z_{g^N}({\psi}_1,\dots ,{\psi}_Z,
\vec{\xi})(x)dV_{g^n}$$

\par Let us suppose that the length of the
$\vec{\xi}$-contraction
(including the $k$ factors $\vec{\xi}$ and the $l$ factors
$S{\nabla}^m\vec{\xi}$) is $L$. We will first integrate by parts the
factors $\vec{\xi}$. Let us pick out one of them at random and
integrate by parts as before, using the familiar formula
(\ref{oneder}). We will then get a sum of
$\vec{\xi}$-contractions as follows:

\begin{equation}
\label{expostintunde}
\begin{split}
&{\int}_{R^N}e^{(N-n)\vec{\xi}\cdot\vec{x}}
\frac{P_h(N)}{(N-n)^{m_h}}\cdot C^Z_{g^N}({\psi}_1,\dots
,{\psi}_Z,\vec{\xi})dV_{g^n}=
\\& -{\int}_{R^N}e^{(N-n)\vec{\xi}\cdot\vec{x}}
\frac{P_h(N)}{(N-n)^{m_h+1}}\cdot {\Sigma}_{\alpha =1}^{L-1}
C^{h,\alpha}_{g^N} ({\psi}_1,\dots ,{\psi}_Z,\vec{\xi})dV_{g^n}
\end{split}
\end{equation}

\par We will sort out the complete contractions according to
 what sort of factor was hit by the derivative ${\nabla}_i$.

\par If $C^{h,\alpha}_{g^N} ({\psi}_1,\dots
,{\psi}_Z,\vec{\xi})$ arises when
the covariant derivative ${\nabla}_i$ hits a factor of the form
${\nabla}^mR_{ijkl}$ or ${\nabla}^pRic$ or
 ${\nabla}^p{\psi}_l$, we get a
 $\vec{\xi}$-contraction with $k-1$ factors $\vec{\xi}$
and $l$ factors $S{\nabla}^m\vec{\xi}$, $m\ge 1$.

\par If $C^{h,\alpha}_{g^N} ({\psi}_1,\dots
,{\psi}_Z,\vec{\xi})$ arises when ${\nabla}_i$
 hits a factor $\vec{\xi}$, we get another
 $\vec{\xi}$-contraction with $k-2$ factors $\vec{\xi}$
and $l+1$ factors $S{\nabla}^m\vec{\xi}$ where $m\ge 1$.

\par Finally, if $C^{h,\alpha}_{g^N} ({\psi}_1,\dots
,{\psi}_Z,\vec{\xi})$ arises when ${\nabla}_i$ hits
a factor $S{\nabla}^m\vec{\xi}$, we get a factor
${\nabla}_iS{\nabla}^m\vec{\xi}$.
We then decompose that factor
 according to equation (\ref{symunsymgenxi}).
In either case, we have reduced by 1 the number of
$\vec{\xi}$-factors.

\par  The good $\vec{\xi}$-contractions
 we have already seen how to
 treat. Finally, if we get an undecided
$\vec{\xi}$-contraction, we have reduced the number of
$\vec{\xi}$-factors.
\newline

{\it The hard $\vec{\xi}$-contractions.}
\newline

\par Let us suppose that $\frac{P(N)}{(N-n)^m}
C^j_{g^n}({\psi}_1,\dots ,
{\psi}_Z,\vec{\xi})$ is a hard $\vec{\xi}$-contraction in
the form (\ref{linisymxi2}) with $k=0$.

\par We pick out one of the $l$ factors
$S{\nabla}^m_{r_1\dots r_m}\vec{\xi}_j$ and write it as
$S{\nabla}^{m-1}_{(r_1\dots r_{m-1}}{\Gamma}^k_{r_mj)}
\vec{\xi}_k$. We then integrate by parts the factor
$\vec{\xi}_k$ and obtain a formula:

\begin{equation}
\label{intbphard}
\begin{split}
&{\int}_{R^N}[e^{(N-n)\vec{\xi}\cdot\vec{x}}
\frac{P(N)}{(N-n)^m}C^j_{g^n} ({\psi}_1,\dots
{\psi}_Z,\vec{\xi})]dV_{g^n}=
\\&
-{\int}_{R^N}[e^{(N-n)\vec{\xi}\cdot\vec{x}}{\Sigma}_{h\in H^j}
\frac{P(N)}{(N-n)^{m+1}}C^h_{g^n} ({\psi}_1,\dots
{\psi}_Z,\vec{\xi})]dV_{g^n}
\end{split}
\end{equation}

where each complete contraction $C^h_{g^n}
({\psi}_1,\dots {\psi}_Z,\vec{\xi})$ is
either in the form (\ref{linisymxi1}) or in the form
(\ref{linisymxi2}) or in the form:

\begin{equation}
 \label{symunintrcontr1}
\begin{split}
&contr({\nabla}_{r_1\dots r_{m_1}}^{m_1}R_{i_1j_1k_1l_1}\otimes
\dots \otimes {\nabla}_{v_1\dots
v_{m_s}}^{m_s}R_{i_sj_sk_sl_s}\otimes \\& {\nabla}_{t_1 \dots
t_{p_1}}^{p_1}Ric_{{\alpha}_1 {\beta}_1} \otimes \dots \otimes
{\nabla}_{z_1 \dots z_{p_q}}^{p_q} Ric_{{\alpha}_q{\beta}_q}
\otimes
{\nabla}^{{\nu}_1}_ {{\chi}_1\dots
{\chi}_{{\nu}_1}}{\psi}_1 \otimes\dots \otimes
{\nabla}^{{\nu}_Z}_{{\omega}_1\dots {\omega}_{{\nu}_Z}}
{\psi}_Z
\\& \otimes
S{\nabla}^{z_1}{\Gamma}^k_{ij}\otimes
S{\nabla}^{w_1}\vec{\xi}\otimes \dots\otimes
S{\nabla}^{w_a}\vec{\xi}\otimes (\vec{\xi}))
\end{split}
\end{equation}

where the symbol $(\vec{\xi})$
means that there may or may not be a factor $\vec{\xi}$.

\par We see that each $C^h_{g^n}
({\psi}_1,\dots {\psi}_Z,\vec{\xi})$ can be taken to be in
 the form (\ref{symunintrcontr1}), by the
 following reasoning:

 If the covariant derivative ${\nabla}_k$ hits a
 factor ${\nabla}^mR_{ijkl}$ or ${\nabla}^pRic$
or ${\nabla}^{\nu}{\psi}_l$,
then  we
 will get a  $\vec{\xi}$-contraction
 in the form (\ref{symunintrcontr1}).
 If it hits a factor $S{\nabla}^m
\vec{\xi}_j$, we apply the formula (\ref{symunsymgenxi})
and we get a linear combination of $\vec{\xi}$-contractions
 in the form  (\ref{symunintrcontr1}). Finally, if
 it hits
 the factor $S{\nabla}^m{\Gamma}^k_{ij}$,
 we will get a complete contraction
 as in (\ref{symunintrcontr1}) with $l-1$ factors
$S{\nabla}^m\vec{\xi}$,
   and with an extra factor ${\nabla}_kS{\nabla}^{m-1}
{\Gamma}^k_{r_mj}$. We then apply the formula:

\begin{equation}
\label{symunsymgengamma}
\begin{split}
&{\nabla}_a S{\nabla}^m_{r_1\dots r_m}{\Gamma}^k_{ij}=
S{\nabla}^{m+1}_{ar_1\dots r_m}{\Gamma}^k_{ij}+
C_m\cdot S^{*}{\nabla}^m_{r_1\dots r_m} {R_{aij}}^k+
\\& {\Sigma}_{u\in U^m} a_u pcontr({\nabla}^{m'}R_{fghj},
S{\nabla}^{x_u} {\Gamma}^k_{bc})
\end{split}
\end{equation}

where the symbol $pcontr({\nabla}^{m'}R_{fghj}, S{\nabla}^{x_u}
{\Gamma}^k_{bc})$ (we call that sublinear combination the {\it
correction terms}) stands for a partial
 contraction of at least one factor ${\nabla}^{m'}R_{fghj}$
against a factor $S{\nabla}^{x_u}{\Gamma}^k_{bc}$ {\it or} a
 partial
contraction of $a\ge 2$ factors $\nabla^{m''}R_{f'g'h'j'}$. We
recall that $S^{*}\nabla^m_{r_1\dots r_m}R_{aijk}$ stands for
 the symmetrization of
 the tensor $\nabla^m_{r_1\dots r_m}R_{aijk}$ over the indices $r_1,
\dots ,r_m,i,j$.

 Furthermore, we have that in each such partial contraction,
 the index $a$ appears in a factor ${\nabla}^{m'}R_{fghj}$.

\par In order to check that in each correction term there can be
at most one factor $S\nabla^p\Gamma$, we only have to observe that
in order to symmetrize a tensor ${\nabla}_kS{\nabla}^{m-1}
{\Gamma}^k_{r_mj}$, we only introduce correction
terms by virtue of the formula $[\nabla_a\nabla_b-\nabla_b\nabla_a
]X_c=R_{abcd}X^d$, and the formula
$\nabla_a\Gamma^k_{bc}-\nabla_b\Gamma^k_{ac}={R_{abc}}^k$. Hence,
for each application of the above formulas, we may decrease the
number of factors $\nabla^p\Gamma$, but we can not increase it.

\par Thus we see that our $\vec{\xi}$-contraction will be
a linear combination of $\vec{\xi}$-contractions in the form
(\ref{symunintrcontr1}) or (\ref{linisymxi2}).

\par So, in general, we must integrate by
 parts expressions of the following form:

$${\int}_{R^N}[e^{(N-n)\vec{\xi}\cdot\vec{x}}
\frac{P(N)}{(N-n)^{m}}C^j_{g^n} ({\psi}_1,\dots
{\psi}_Z,\vec{\xi})]dV_{g^n}$$

where the complete contraction $C^j_{g^n}
({\psi}_1,\dots {\psi}_Z,\vec{\xi})$ is in the form:

\begin{equation}
 \label{symunintrcontr}
\begin{split}
&contr({\nabla}_{r_1\dots r_{m_1}}^{m_1}R_{i_1j_1k_1l_1}\otimes
\dots \otimes {\nabla}_{v_1\dots
v_{m_s}}^{m_s}R_{i_sj_sk_sl_s}\otimes \\& {\nabla}_{t_1 \dots
t_{p_1}}^{p_1}Ric_{{\alpha}_1 {\beta}_1} \otimes \dots \otimes
{\nabla}_{z_1 \dots z_{p_q}}^{p_q} Ric_{{\alpha}_q{\beta}_q}
\otimes
{\nabla}^{{\nu}_1}_ {{\chi}_1\dots
{\chi}_{{\nu}_1}}{\psi}_1 \otimes\dots \otimes
{\nabla}^{{\nu}_Z}_{{\omega}_1\dots {\omega}_{{\nu}_Z}}
{\psi}_Z
\\& \otimes
S{\nabla}^{z_1}{\Gamma}^k_{ij}\otimes\dots\otimes
S{\nabla}^{z_v}{\Gamma}^k_{ij}\otimes S{\nabla}^{u_1}\vec{\xi}\otimes \dots
 \otimes S{\nabla}^{u_d}\vec{\xi}\otimes \vec{\xi}
\otimes\dots \otimes \vec{\xi})
\end{split}
\end{equation}

\par  The integration by parts of such complete contractions
 can be done as before: If there is a factor $\vec{\xi}$ then we
integrate by parts using it, and symmetrize and anti-symmetrize as
will be explained below. If there is no factor $\vec{\xi}$, we
pick out one
 factor $S{\nabla}^m_{r_1\dots r_m}\vec{\xi}_j$, we write it as
$S{\nabla}^{m-1}_{r_1\dots r_{m-1}}{\Gamma}^k_{r_m j}
\vec{\xi}_k$. We then integrate by parts with respect to the factor
$\vec{\xi}_k$, using the formula (\ref{oneder}). If the derivative
 ${\nabla}_k$ hits a factor
${\nabla}^mR_{ijkl}$, or ${\nabla}^p{\psi}_l$, or
${\nabla}^pRic$, we leave them as they are. If it hits a
 factor $S{\nabla}^x{\Gamma}_{ij}^k$ or a factor
$S{\nabla}^m\vec{\xi}$, we apply the formulas,
(\ref{symunsymgenxi}), (\ref{symunsymgengamma})
respectively.

\par In the end, we will have the following formula for the integration
by parts of a hard $\vec{\xi}$-contraction
$C_{g^N}({\psi}_1,\dots ,{\psi}_Z,\vec{\xi})$:

\begin{equation}
\label{exanalund}\begin{split}
&{\int}_{R^N}e^{(N-n)\vec{\xi}\cdot\vec{x}}
 \frac{P(N)}{(N-n)^{m}}
C_{g^N}({\psi}_1,\dots ,{\psi}_Z,\vec{\xi})dV_{g^n}=
\\&{\int}_{R^N}e^{(N-n)\vec{\xi}\cdot\vec{x}}\frac{P(N)}{(N-n)^{m'}}
{\Sigma}_{s\in S} a_s C^s_{g^N}({\psi}_1,\dots ,{\psi}_Z)dV_{g^n}
\end{split}
\end{equation}

where the degree of the rational function
$\frac{P(N)}{(N-n)^{m'}}$ will be zero and the complete
contractions $C^s_{g^N}({\psi}_1,\dots ,{\psi}_Z)$ will be in the
following general form:

\begin{equation}
 \label{vgenlinisymxi}
\begin{split}
&contr({\nabla}_{r_1\dots r_{m_1}}^{m_1}R_{i_1j_1k_1l_1}\otimes
\dots \otimes {\nabla}_{v_1\dots
v_{m_s}}^{m_s}R_{i_sj_sk_sl_s}\otimes \\& {\nabla}_{t_1 \dots
t_{p_1}}^{p_1}Ric_{{\alpha}_1 {\beta}_1} \otimes \dots \otimes
{\nabla}_{z_1 \dots z_{p_q}}^{p_q} Ric_{{\alpha}_q{\beta}_q}
\otimes
\\&{\nabla}^{{\nu}_1}_ {{\chi}_1\dots
{\chi}_{{\nu}_1}}{\psi}_1 \otimes\dots \otimes
{\nabla}^{{\nu}_Z}_{{\omega}_1\dots {\omega}_{{\nu}_Z}}
 {\psi}_Z \otimes S{\nabla}^{x_1}{\Gamma}^{k_1}_{ij}
\otimes\dots\otimes S{\nabla}^{x_u}{\Gamma}^{k_u}_{ij})
\end{split}
\end{equation}

where $u\ge 0$. Therefore, by virtue of (\ref{symunsymgengamma}),
we see that if $C_{g^n}({\psi}_1,\dots ,{\psi}_Z)$ is hard, then
the integrand on the right hand side of (\ref{exanalund}) may
 apriori contain complete contractions in the form
(\ref{linicontrpsi}). We accept this for the time being,
 although we will later show, in Lemma \ref{finecanc} that,
 in fact, there will be cancellation among such complete
 contractions.

\subsection{The simple divergence formula}

\par Therefore, after a series of integrations by parts, the relation
(\ref{bignintxi}) can be brought into the form:

\begin{equation}
\label{laplzero}\begin{split} &{\int}_{\mathbb{R}^N}
e^{(N-n)\vec{\xi}\cdot\vec{x}}[I^Z_{g^N}({\psi}_1,\dots ,{\psi}_Z)
-{\Sigma}_{a\in A} {\alpha}_a \frac{P_a(N)}{(N-n)^{r_a}}
C^a_{g^N}({\psi}_1,\dots ,{\psi}_Z)
\\& - {\Sigma}_{b\in B}
{\beta}_b \frac{P_b(N)}{(N-n)^{r_b}} C^b_{g^N}({\psi}_1,\dots
,{\psi}_Z)]dV_{g^n}=0
\end{split}
\end{equation}

Where $deg[P_a(N)]=r_a$ and $deg[P_b(N)]<r_b$. The complete
contractions

\noindent $\frac{P_a(N)}{(N-n)^{r_a}} C^{Z,a}_{g^N}({\psi}_1,\dots
,{\psi}_Z)(x)$ have arisen from iterated integrations by parts of
the good, the hard and the  undecided complete contractions. They
are in the form (\ref{linisym}) or (\ref{vgenlinisymxi}).
 We may assume with no loss of generality
that the leading order coefficient of each of the polynomials
$P_a(N)$ is 1, incorporating it in ${\alpha}_a$.

\par The complete contractions $C^b_{g^N}({\psi}_1,\dots
,{\psi}_Z)(x)$ have arisen from the easy complete
contractions. All of the complete contractions in the formula
(\ref{laplzero}) have arisen according to the procedure we
outlined in the previous subsection.

\par Now, relation (\ref{laplzero}) shows us that the quantity
between brackets is zero for every $x\in
B(\tilde{x}_0, \epsilon)$. In
particular,

\begin{equation}
\label{prsimp}\begin{split} &I^Z_{g^N}({\psi}_1,\dots
,{\psi}_Z)(\tilde{x}_0) -{\Sigma}_{a\in A} {\alpha}_a
\frac{P_a(N)}{(N-n)^{r_a}} C^a_{g^N}({\psi}_1,\dots
,{\psi}_Z)(\tilde{x}_0)-
\\&{\Sigma}_{b\in B} {\beta}_b
\frac{P_b(N)}{(N-n)^{r_b}} C^b_{g^N}({\psi}_1,\dots
,{\psi}_Z)(\tilde{x}_0)=0
\end{split}
\end{equation}

for every Riemannian manifold  $(M^N,g^N)$, any functions
${\psi}_1,\dots ,{\psi}_Z$ around $\tilde{x}_0 \in M^N$ and
 any coordinate
system around $\tilde{x}_0\in M^N$. Now pick any $(M^n,g^n)$,
 any $x_0\in M^n$ and any coordinate system around $x_0$. We
define $M^N=M^n\times S^1\times \dots \times S^1$ ($S^1$ has the standard flat metric and $g^N$ is the product metric).
 We pick $\tilde{x}_0=(x_0,0,\dots ,0)$ and consider the induced coordinate system around $\tilde{x}_0$. It follows
 that:

\begin{equation}
\label{presimp}
\begin{split}
& I^Z_{g^n}({\psi}_1,\dots ,{\psi}_Z)(x_0)
-{\Sigma}_{a\in A} {\alpha}_a \frac{P_a(N)}{(N-n)^{r_a}}
C^a_{g^n}({\psi}_1,\dots ,{\psi}_Z)(x_0)-
\\&{\Sigma}_{b\in B}
{\beta}_b \frac{P_b(N)}{(N-n)^{r_b}} C^b_{g^n}({\psi}_1,\dots
,{\psi}_Z)(x_0)=0
\end{split}
\end{equation}

for every Riemannian manifold  $(M^n,g^n)$, any functions
${\psi}_1,\dots ,{\psi}_Z$ around $x_0 \in M^n$ and any coordinate
system around $x_0\in M^n$.

\par In equation (\ref{presimp}), $N$ is just a free variable.
Hence, we can take the limit as $N\longrightarrow \infty$ in
(\ref{presimp}) and obtain the {\it simple divergence formula}:

\begin{equation}
\label{simp}
 I^Z_{g^n}({\psi}_1,\dots ,{\psi}_Z)(x_0)
-{\Sigma}_{a\in A} {\alpha}_a C^a_{g^n}({\psi}_1,\dots
,{\psi}_Z)(x_0) =0
\end{equation}

\par So we have disposed of the integrations by parts
of the easy complete contractions.

\section{The three refinements of the simple divergence
 formula.}
\label{threerefs}
\subsection{The first refinement: Separating intrinsic from
un-intrinsic complete contractions.}

\par  We recall from the previous section that some of the
 complete contractions in
(\ref{simp}) will be in the form (\ref{linisym}). On the
other hand, we have also found that there will be complete
contractions in the general form (\ref{vgenlinisymxi}),
with $u\ge 1$.  Accordingly, we
introduce the following dichotomy:

\begin{definition}
\label{intrinsic} Complete contractions in the form
(\ref{linicontrpsi}) or (\ref{linisym}) will be called
{\it intrinsic}. Complete
contractions in the general form (\ref{vgenlinisymxi}) with
$u>0$ will be called un-intrinsic.
\end{definition}

 Let us consider, in (\ref{simp}), the two sub-linear
combinations of the intrinsic  and of the
un-intrinsic complete contractions. Written that way,
(\ref{simp}) will be:

\begin{equation}
\label{tofilter}
I^Z_{g^n}({\psi}_1,\dots ,{\psi}_Z)
-{\Sigma}_{l\in L} {\alpha}_l C^l_{g^n}({\psi}_1,\dots
,{\psi}_Z) -{\Sigma}_{r\in R} {\alpha}_r
C^r_{g^n}({\psi}_1,\dots ,{\psi}_Z) =0
\end{equation}

where the complete contractions $C^l_{g^n}({\psi}_1,\dots
,{\psi}_Z)(x)$ are the intrinsic ones and the complete
contractions $C^r_{g^n}({\psi}_1,\dots ,{\psi}_Z)(x)$ are the
un-intrinsic ones. We have, of course, that $L\bigcup R=A$.

\par Our next goal will be to show that:

\begin{equation}
\label{intrinsic}
I^Z_{g^n}({\psi}_1,\dots ,{\psi}_Z)(x)
-{\Sigma}_{l\in L} {\alpha}_l C^l_{g^n}({\psi}_1,\dots
,{\psi}_Z)(x)=0
\end{equation}

which is equivalent to proving:

\begin{equation}
\label{unintrinsic}
{\Sigma}_{r\in R} {\alpha}_r
C^r_{g^n}({\psi}_1,\dots ,{\psi}_Z)(x)=0
\end{equation}

\par So let us
focus on showing (\ref{intrinsic}). Let us treat the
value of the left hand side of (\ref{simp}) as a function of the
coordinate system. We want to show, roughly speaking, that the
tensors $S{\nabla}^m{\Gamma}^k_{ij}(x_0)$ are not independent
 of the coordinate
system in which they are expressed. In other words, they are not
intrinsic tensors of the Riemannian manifold $(M^n,g^n)$.

\begin{lemma}
\label{killunintrinsic}
 We have that (\ref{intrinsic}) holds.
\end{lemma}

{\it Proof:} We consider the tensors $S{\nabla}^m_{s_1\dots
s_m}{\Gamma}^k_{ij}(x_0)$, ${\Gamma}^k_{ij}(x_0)$, written out in
any coordinate system.  We want to see what their values can be,
given our metric $g^n$ around $x_0$.

\par We need to recall the following fact from \cite{e:rg}: Consider a coordinate
transformation around the point $x_0\in M^n$. Let us say we
had coordinates $\{x^1,\dots ,x^n\}$ and now we have coordinates
$\{y^1,\dots ,y^n\}$. Then the Christoffel symbols ${\Gamma}^k_{ij}$
will transform as follows:

\begin{equation}
\label{Chris}
\tilde{\Gamma}^{\lambda}_{\mu\nu}\frac{\partial x^l}{\partial y^{\lambda}}(x_0)={\Gamma}^l_{ij}(x_0)
\frac{\partial x^i}{\partial y^{\mu}}
\frac{\partial x^j}{\partial y^{\nu}}+
\frac{{\partial}^2x^l}{\partial y^{\mu}\partial y^{\nu}}(x_0)
\end{equation}
($\tilde{\Gamma}^{\lambda}_{\mu\nu}(x_0)$ stands for the Christoffel symbols in the new coordinate system).

\par Now, it is an elementary fact that the tensors ${\nabla}^mR_{ijkl}$
 are intrinsic tensors of the Riemannian manifold.
That means that they satisfy the intrinsic transformation law
under coordinate changes, as in \cite{e:rg}.

\par We will need the following Lemma:

\begin{lemma}
\label{prescr} Consider a point $x_0\in M^n$ and a coordinate
system $\{x^1,\dots ,x^n\}$ around $x_0$ for which
$g^n_{ij}(x_0)={\delta}_{ij}$. Then, given any list of special
tensors $T^k_{r_1\dots r_{p+2}}$, which are symmetric in the
indices
 $r_1,\dots r_{p+2}$, there
 is a coordinate system $\{y^1,\dots ,y^n\}$ around
$x_0\in M^n$ so
 that the tensors $S{\nabla}^p_{r_1\dots r_p}
{\Gamma}^k_{r_{p+1}r_{p+2}}$ have the values of the
 arbitrarily chosen tensors $T^k_{r_1\dots r_{p+2}}$ at
$x_0$ and furthermore we have that
$[\frac{\partial y}{\partial x}](x_0)=Id^{n\times n}$ and $g^n_{ij}={\delta}_{ij}$
 (this is with respect to the new coordinate system).
\end{lemma}

{\it Proof:} Let us observe that by \cite{e:rg} when we change the coordinate system $\{ x^1,\dots x^n\}$ into
$\{y^1,\dots ,y^n\}$, the tensors
${\nabla}^p_{r_1\dots r_p}{\Gamma}^l_{r_{p+1}r_{p+2}}$ will
transform as follows:

\begin{equation}
\label{trgamma}
\begin{split}
&{\nabla}^p_{r'_1\dots r'_p}\tilde{\Gamma}^{\lambda}_{r'_{p+1}
r'_{p+2}}\frac{\partial x^l}{\partial y^{\lambda}}(x_0)=
{\nabla}^p_{r_1\dots r_p}{\Gamma}^l_{r_{p+1}r_{p+2}}
\frac{\partial x^{r_1}}{\partial y^{r'_1}}\dots
\frac{\partial x^{r_{p+2}}}{\partial y^{r'_{p+2}}}(x_0)
\\&+\frac{{\partial}^{p+2} x^l}{\partial y^{r'_1}\dots
\partial  y^{r'_{p+2}}}(x_0)+ {\Sigma}({\partial}^f \Gamma ,
\frac{{\partial}^h x}{{\partial}^h y})(x_0)
\end{split}
\end{equation}

where $\tilde{\Gamma}$ stands for the Christoffel symbols in the new coordinate system and ${\Sigma}({\partial}^f \Gamma ,
\frac{{\partial}^h x}{{\partial}^h y})$ stands for a linear
 combination of partial contractions of factors against
 factors $\frac{{\partial}^h x}{{\partial}^h y}$
with $h<p+2$.

\par Now, we can prescribe $\frac{{\partial}^{p+2}x^l}{\partial y^{r'_1}\dots \partial y^{r'_{p+2}}}(x_0)$ to have
 any symmetric value in the indices $r_1,\dots ,r_{p+2}$.
 Thus, if we write out the transformation law for
$S{\nabla}^p_{r_1\dots r_p}{\Gamma}^k_{r_{p+1}r_{p+2}}(x_0)$
 under coordinate changes, then the linearized part
 of its transformation law will be precisely
$\frac{{\partial}^{p+2}y^l}{\partial x^{r_1}\dots \partial
x^{r_{p+2}}}(x_0)$.
Hence, by induction on $p$, we have our Lemma.
 $\Box$
\newline

\par Let us call these arbitrary tensors $T^k_{r_1\dots r_{p+2}}$ the
{\it un-intrinsic free variables}. By construction, they satisfy
$T^k_{r_1\dots r_{p+2}}(t^2g^n)= T^k_{r_1\dots r_{p+2}}(g^n)$.
Thus, they are special tensors.

\par But then it is straightforward to check Lemma
\ref{killunintrinsic}.
We can break equation (\ref{tofilter}) into two summands: The left hand side of (\ref{intrinsic})
 plus the left hand side of (\ref{unintrinsic}).
We may then pick any $\lambda\in \mathbb{R}$ and pick a new
coordinate system so that $(S\nabla^p\Gamma^k_{ij})'(x_0)=
\lambda\cdot (S\nabla^p\Gamma^k_{ij})(x_0)$. (Here
$(S\nabla^p\Gamma^k_{ij})'(x_0)$ stands for the value of
$S\nabla^p\Gamma^k_{ij}(x_0)$ with respect to the new coordinate system). We can then see the left hand side of
(\ref{tofilter}) as a
 polynomial in $\lambda$, $\Pi(\lambda)$. We have that the constant term of $\Pi(\lambda)$ must be zero. Also, the constant term of
$\Pi(\lambda)$ is precisely the left hand side of
(\ref{intrinsic}). We have shown our Lemma. $\Box$

\subsection{The second refinement: An intrinsic divergence formula.}

\par We begin this subsection with one more convention.
When we have an equation of the form:

\begin{equation}
\label{christodoulou}
\Sigma_{l\in L} a_l C^l_{g^n}(\psi_1,\dots ,\psi_Z,\vec{\xi})=0
\end{equation}
where each $C^l_{g^n}(\psi_1,\dots ,\psi_Z,\vec{\xi})$ is a
complete contraction in the form (\ref{linisymxi2}), we will be
thinking of the factors $S\nabla^m_{r_1\dots r_m}\vec{\xi}_j$
($m\ge 0$) as symmetric $(m+1)$-tensors in the indices $r_1,\dots
,r_m,j$ so that $S\nabla^m_{r_1\dots r_m}\vec{\xi}_j(t^2g^n)=
S\nabla^m_{r_1\dots r_m}\vec{\xi}_j(g^n)$. This condition
trivially holds since $\nabla^m_{r_1\dots r_m}\vec{\xi}_j=\nabla^{m+1}_{r_1\dots r_mj}(\vec{x}\cdot\vec{\xi})$.
Moreover, we will be implying that the above equation
holds for every $x_0\in (M^n,g^n)$ ($g^n$ can be any Riemannian
metric), any functions $\psi_1,\dots ,\psi_Z$ defined around
$x_0$, any vector $\vec{\xi}\in \mathbb{R}^n$ and any coordinate
system defined around
 $x_0$.

\par Now, we define $C^l_{g^n}(\psi_1,\dots ,\psi_Z,\Omega)$ to stand
 for complete contraction that arises from
$C^l_{g^n}(\psi_1,\dots ,\psi_Z,\vec{\xi})$ by replacing each factor
$S\nabla^m_{r_1\dots r_m}\vec{\xi}_j$ by an auxiliary symmetric tensor
$\Omega_{r_1\dots r_m j}$.
We claim:

\begin{lemma}
\label{cristodoulou2}
Assuming (\ref{christodoulou}) (as explained above), we have that:

\begin{equation}
\label{christodoulou3}
\Sigma_{l\in L} a_l C^l_{g^n}(\psi_1,\dots ,\psi_Z,\Omega)=0
\end{equation}
will hold for every
$x_0\in (M^n,g^n)$ ($g^n$ can be any Riemannian metric), any functions
$\psi_1,\dots ,\psi_Z$ defined around $x_0$ and {\it any} symmetric tensors $\Omega_{i_1\dots i_s}$.
\end{lemma}

{\it Proof:} Firstly, we observe that for every sequence $\Omega_{i_1},
\dots ,\Omega_{i_1\dots i_s},\dots $ of symmetric tensors for which
$\Omega_{i_1}\ne 0$, we have that there is vector $\vec{\xi}\in
\mathbb{R}^n$ and also a coordinate system around $x_0\in M^n$ so that:
$\vec{\xi}_{i_1}=\Omega_{i_1},\dots ,S\nabla^m_{r_1\dots r_m}
\vec{\xi}_j=\Omega_{r_1\dots r_mj},\dots$. This is clear by virtue of the formula $S\nabla^m_{r_1\dots r_m}\vec{\xi}_j=
S\nabla^{m-1}_{r_1\dots r_{m-1}}\Gamma^k_{r_mj}\vec{\xi}_k$ and by Lemma \ref{prescr}.

\par Now, for any sequence
$\Omega_{i_1},\dots ,\Omega_{i_1\dots i_s},\dots $ where
$\Omega_{i_1}= 0$, we only have to consider any vector
$\vec{\epsilon}_i$ where $|\vec{\epsilon}_i|$ is small. We then
have
 that there is a coordinate system so that
$\vec{\xi}_i=\vec{\epsilon}_i$ and $S\nabla^m_{r_1\dots r_m}\vec{\xi}_j=
\Omega_{r_1\dots r_m j}$, for every $m\ge 1$.
 Letting $\vec{\epsilon}_i\longrightarrow 0$, we obtain
 our Lemma. $\Box$
\newline

\par Now, the aim of this subsection is to further refine Lemma
\ref{killunintrinsic}. We will need certain preliminary observations.

\par Notice the following: Let us pick out one
$\vec{\xi}$-contraction $Q(N)\cdot C^l_{g^N}({\psi}_1,\dots ,
{\psi}_Z,\vec{\xi})$ of the form (\ref{linisymxi2}) with $k+l\le
|\vec{\xi}|$.   We have then treated the integrals

$${\int}_{R^N}e^{(N-n)\vec{\xi}\cdot\vec{x}}Q(N)\cdot C^l_{g^N}({\psi}_1,\dots
,{\psi}_Z,\vec{\xi})dV_{g^n}$$

and performed integrations by parts, obtaining a relation

\begin{equation}
\label{intbp}
\begin{split}
&{\int}_{R^N}e^{(N-n)\vec{\xi}\cdot\vec{x}}Q(N)\cdot C^l_{g^N}({\psi}_1,\dots
,{\psi}_Z,\vec{\xi})dV_{g^n}=
\\&{\int}_{R^N}e^{(N-n)\vec{\xi}\cdot\vec{x}}{\Sigma}_{s\in S^l}
 Q_s(N) C^{l,s}_{g^N}({\psi}_1,\dots ,{\psi}_Z)dV_{g^n}
\end{split}
\end{equation}

where the degree of the rational function $Q_s(N)$ is zero. Adding
up all the integrations by parts, writing things in dimension $n$
and taking the limit $N\longrightarrow \infty$ gives us the
formula (\ref{simp}). Let us call this procedure by which we
integrate by parts one $\vec{\xi}$-factor at a time the {\it
iterative procedure} of integrating by parts.

\par After we do all the integrations by parts for a $\vec{\xi}$-contraction
  as in (\ref{intbp}), we will call the
 quantity:

$$lim_{N\longrightarrow\infty}{\Sigma}_{s\in S} Q_s(N)
C^{l,s}_{g^n}({\psi}_1,\dots ,{\psi}_Z)$$

the {\it final outcome} of the iterative integration by parts. We
will denote the final outcome by $F[Q(N)C^l_{g^n}({\psi}_1,
\dots ,{\psi}_Z,\vec{\xi})]$.

\par Recall that we are assuming the leading order coefficient of
$Q_s(N)$ to be 1. We make a further notational convention: When we
write out the good or undecided or hard $\vec{\xi}$-contractions
and also when we integrate by parts, we will be omitting the
 dimensional rational function $Q_i(N)$.
 This is justified by the fact that we eventually take a
 limit $N\longrightarrow \infty$.
So all the formulas that will appear in this section
 will be true after we take the limit  $N\longrightarrow \infty$. We will refer to this
notational convention as the $N$-cancelled notation.

\par As an example of this notational convention, we will
be saying that we apply the third summand on the right hand side
of the formula (\ref{levicivitaxi}) to the pair $({}^m,{}_m)$ in
$\nabla^mR_{mjkl}$ and bring out $\vec{\xi}^mR_{mjkl}$ instead of
saying that we bring out $N\vec{\xi}^mR_{mjkl}$. Also, we will say
that we substitute a factor $Ric_{ij}$ by $-{\nabla}_i\vec{\xi}_j$
or a factor $R$ by $-|\vec{\xi}|^2$ (instead of
$-N{\nabla}_i\vec{\xi}_j$ or
 $-N^2|\vec{\xi}|^2$ respectively). What is meant
 will be clear.

\begin{observation}
\label{toremark}
We note that the formal expression
for $F[C^l_{g^n}({\psi}_1, \dots ,{\psi}_Z,\vec{\xi})]$
 depends on the order in which we perform the integrations by
 parts. In general, whenever we make reference to the integrations by parts,
 we will be assuming that we arbitrarily
 pick an order in which we perform integrations by parts
-subject to the restrictions imposed in the corresponding
 section or any extra restrictions we wish to impose.
\end{observation}

\par Now, some conventions in order to state and prove our
 Lemma for this subsection:

\begin{definition}
\label{outgrowth}
In $N$-cancelled notation:
Consider any good or undecided $\vec{\xi}$-contraction
$C_{g^n}({\psi}_1,\dots ,{\psi}_Z,\vec{\xi})$,
 in the form (\ref{linisymxi1}) or (\ref{linisymxi2}).
The number of $\vec{\xi}$-factors in
$C_{g^n}({\psi}_1,\dots ,{\psi}_Z,\vec{\xi})$ is $k+l$.
We perform the iterative integrations by parts, subject to the following restriction:
 In each step of the iterative integration by parts, suppose
 we start off with $X$ $\vec{\xi}$-factors. We integrate by parts with respect to a factor $\vec{\xi}_i$ and we obtain a linear combination of $\vec{\xi}$-contractions (each in the form (\ref{linisymxi1}) or (\ref{linisymxi2})),
each with $X-1$ $\vec{\xi}$-factors. In that linear
 combination we cross
out the hard $\vec{\xi}$-contractions.
 We then pick out one of the $\vec{\xi}$-contractions we are left with
 (it will either be good or undecided) and again integrate by parts with
 respect to a factor $\vec{\xi}$. After $k+l$ steps, this procedure will
 terminate and we will be left with an expression:

$${\int}_{\mathbb{R}^N}e^{(N-n)\vec{\xi}\cdot\vec{x}}
{\Sigma}_{h\in H} a_h Q_h(N)C^h_{g^N}({\psi}_1,\dots ,
{\psi}_Z)dV_{g^N}$$

Each complete contraction $C^h_{g^N}({\psi}_1,\dots ,
{\psi}_Z)$ is in the form (\ref{linisym}) and the rational
 function $Q_h(N)$ has degree $0$ and leading order
 coefficient $1$.

We then define ${\Sigma}_{h\in H} a_h C^h_{g^n}({\psi}_1,
\dots ,{\psi}_Z)$ to be the
 {\it outgrowth} of $C_{g^n}({\psi}_1,\dots ,{\psi}_Z,\vec{\xi})$.
 We will denote it by $O[C_{g^n}({\psi}_1,\dots ,{\psi}_Z,\vec{\xi})]$.
\end{definition}

\par We want to prove the following:

\begin{proposition}
\label{2ndref}
Consider the sublinear combination
$\tilde{S}^Z_{g^n}({\psi}_1,\dots ,{\psi}_Z,\vec{\xi})$ of
\\ $S^Z_{g^n}({\psi}_1,\dots ,{\psi}_Z,\vec{\xi})$
which consists of the good and the undecided
$\vec{\xi}$-contractions.
\par We then claim that if in our $N$-cancelled notation we
 have

$$\tilde{S}^Z_{g^n}({\psi}_1,\dots ,{\psi}_Z,\vec{\xi})=
{\Sigma}_{l\in L}a_l C^l_{g^n}({\psi}_1,\dots ,{\psi}_Z,\vec{\xi})$$

then we have the following:

\begin{equation}
\label{Ointrinsic}
I^Z_{g^n}({\psi}_1,\dots ,{\psi}_Z)+ {\Sigma}_{l\in L}a_lO[C^l_{g^n}({\psi}_1,\dots ,{\psi}_Z,\vec{\xi})]=0
\end{equation}
\end{proposition}

{\it Proof of Proposition \ref{2ndref}:}
\newline

\par Let us suppose that the linear combination of the
 hard $\vec{\xi}$-contractions that we encounter, along
  the iterative integration by parts of $C^l_{g^n}({\psi}_1,
\dots ,{\psi}_Z,\vec{\xi})$ is
${\Sigma}_{b\in B^l} a_b C^b_{g^n}({\psi}_1,\dots ,
{\psi}_Z,\vec{\xi})$. We can then write:

\begin{equation}
\label{breakf}
F[C^l_{g^n}({\psi}_1,\dots ,{\psi}_Z,\vec{\xi})]=
O[C^l_{g^n}({\psi}_1,\dots ,{\psi}_Z,\vec{\xi})]+
{\Sigma}_{b\in B^l} a_b
F[C^b_{g^n}({\psi}_1,\dots ,{\psi}_Z,\vec{\xi})]
\end{equation}

\par Now, let us note the following fact for the final
 outcome of the iterative integration by parts of a
 hard $\vec{\xi}$-contraction
$C^u_{g^n}(\phi,\vec{\xi})$, in the form (\ref{linisymxi2})
 with $k=0$, $l>0$.

\begin{lemma}
\label{fardtod}
Let us suppose that:

$$F[C^u_{g^n}(\phi,\vec{\xi})]={\Sigma}_{y\in Y^u} a_y
C^y_{g^n}(\phi)$$

Then, there will be one complete contraction
$a_y C^y_{g^n}(\phi)$
(along with its coefficient) which is obtained from
$C^u_{g^n}(\phi,\vec{\xi})$
by substituting each of the $l$ factors
$S{\nabla}^m_{r_1\dots r_m}
\vec{\xi}_j$ by $-S{\nabla}^m_{dr_1\dots r_{m-1}}
{\Gamma}^d_{r_m j}$. That complete contraction arises when each derivative ${\nabla}^d$ which arises from a factor
$S{\nabla}^m_{r_1\dots r_m}
\vec{\xi}_j= S{\nabla}^{m-1}_{r_1\dots r_{m-1}}
{\Gamma}^d_{r_mj}\vec{\xi}_d$
 hits the factor $S{\nabla}^{m-1}_{r_1\dots r_{m-1}}
{\Gamma}^d_{r_mj}$, and then we symmetrize using
(\ref{symunsymgengamma}). We will denote this complete
 contraction by
$DF[C_{g^n}^l({\psi}_1,\dots ,{\psi}_Z,\vec{\xi})]$.

\par Furthermore, we claim
that each other complete contraction $a_y C^y_{g^n}(\phi)$ in
$F[C^u_{g^n}(\phi,\vec{\xi})]$ will have strictly less than $l$
un-intrinsic free
 variables $S{\nabla}^m_{r_1\dots r_m}{\Gamma}^d_{ij}$ for which $d$ contracts against one of the indices
$r_1,\dots ,r_m,i,j$ in $C^y_{g^n}(\phi)$.
\end{lemma}

{\it Proof:} This follows from the procedure by which we integrate by parts and also from
 the formula (\ref{symunsymgengamma}). $\Box$
\newline

\par Now, let us consider the sublinear combination of
 good, hard and undecided $\vec{\xi}$-contractions in
$S^Z_{g^n}({\psi}_1,\dots ,{\psi}_Z,\vec{\xi})$:

$${\Sigma}_{l\in L} a_l C^l_{g^n}({\psi}_1,\dots ,
{\psi}_Z,\vec{\xi})$$

\par Now, we break up the index set $L$ as follows:
$l\in L^1$ if and only if $C^l_{g^n}({\psi}_1,\dots ,{\psi}_Z,
\vec{\xi})$ is good or undecided. $l\in L^2$ if and only if
$C^l_{g^n}({\psi}_1,\dots ,{\psi}_Z,\vec{\xi})$ is hard.

\par For any hard $\vec{\xi}$-contraction
$C_{g^n}(\phi,\vec{\xi})$, we break up the linear combination
$F[C_{g^n}(\phi,\vec{\xi})]$ into the
 sublinear combination $F^{Intr}[C_{g^n}(\phi,\vec{\xi})]$
of intrinsic complete contractions and
 the sublinear combination
$F^{UnIntr}[C_{g^n}(\phi,\vec{\xi})]$ of unintrinsic
 complete contractions.
\newline

\par We can then rewrite (\ref{intrinsic}) as:

\begin{equation}
\label{reformintrinsic}
\begin{split}
&I^Z_{g^n}({\psi}_1,\dots ,{\psi}_Z) +{\Sigma}_{l\in L^1}a_l
O[C^l_{g^n}({\psi}_1,\dots ,{\psi}_Z,\vec{\xi})]+
\\& {\Sigma}_{l\in L^1}a_l
{\Sigma}_{b\in B^l} a_b F^{Intr}
[C^b_{g^n}({\psi}_1,\dots ,{\psi}_Z,\vec{\xi})] +
{\Sigma}_{l\in L^2} a_l F^{Intr}[C^l_{g^n}({\psi}_1,\dots
,{\psi}_Z,\vec{\xi})]=0
\end{split}
\end{equation}

and also (\ref{unintrinsic}) as:

\begin{equation}
\label{reformunintrinsic}
{\Sigma}_{l\in L^1}a_l
{\Sigma}_{b\in B^l} a_b F^{UnIntr}
[C^b_{g^n}({\psi}_1,\dots ,{\psi}_Z,\vec{\xi})]+
{\Sigma}_{l\in L^2} a_l F^{UnIntr}[C^l_{g^n}({\psi}_1,\dots
,{\psi}_Z,\vec{\xi})]=0
\end{equation}

\par We then claim:

\begin{lemma}
\label{allo2ndref}
We have that:

$${\Sigma}_{l\in L^1}a_l
{\Sigma}_{b\in B^l} a_b F^{Intr}
[C^b_{g^n}({\psi}_1,\dots ,{\psi}_Z,\vec{\xi})] +
 {\Sigma}_{l\in L^2} a_l F^{Intr}[C^l_{g^n}({\psi}_1,\dots
,{\psi}_Z,\vec{\xi})]=0$$
for every $(M^n,g^n)$, every ${\psi}_1,\dots ,{\psi}_Z$.
\end{lemma}

Notice that if we prove the above we will have our
 Proposition.

{\it Proof:} We will in fact prove a stronger statement.
We claim:

\begin{lemma}
\label{finecanc} We will have that:

\begin{equation}
\label{todo}
{\Sigma}_{l\in L^1}a_l
{\Sigma}_{b\in B^l} a_b
C^b_{g^n}({\psi}_1,\dots ,{\psi}_Z,\vec{\xi}) +
 {\Sigma}_{l\in L^2} a_l C^l_{g^n}({\psi}_1,\dots
,{\psi}_Z,\vec{\xi})=0
\end{equation}

for any point $x_0$, any metric $g^n$ around $x_0$, any functions ${\psi}_1,\dots ,{\psi}_Z$ and any coordinate
 system.
\end{lemma}

{\it Proof that Lemma \ref{allo2ndref} follows from Lemma
\ref{finecanc}:}
\newline

\par We consider the linear combination
$${\Sigma}_{l\in L^1}a_l
{\Sigma}_{b\in B^l} a_b N^{p_b}Q_b(N) C^b_{g^N}({\psi}_1,\dots
,{\psi}_Z,\vec{\xi}) +
 {\Sigma}_{l\in L^2} a_l N^{p_l}Q_l(N) C^l_{g^N}({\psi}_1,\dots
,{\psi}_Z,\vec{\xi})$$ of the hard $\vec{\xi}$-contractions that
we have put aside, {\it without} the $N$-cancelled notation.
 Here, the rational functions $Q_b(N), Q_l(N)$ have degree zero
 and leading order coefficient 1. Moreover, if we denote by
 $|\vec{\xi}|_b, |\vec{\xi}|_l$ the number of $\vec{\xi}$-factors
 in $C^b_{g^n}({\psi}_1,\dots
,{\psi}_Z,\vec{\xi})$, $C^l_{g^n}({\psi}_1,\dots
,{\psi}_Z,\vec{\xi})$ respectively, we will have that $p_b
=|\vec{\xi}|_b$ and $p_l=|\vec{\xi}|_l$.

\par For the purposes of this proof, we will consider any hard or easy
$\vec{\xi}$-contraction $N^a Q_a(N) C_{g^N}({\psi}_1,\dots
,{\psi}_Z,\vec{\xi})$, where $Q_a(N)$ has degree zero and leading
order coefficient one. We perform the iterative integrations by
parts, as explained in the previous subsection, and obtain a
relation:

\begin{equation}
\label{kudos} \begin{split} & \int_{\mathbb{R}}
e^{(N-n)\vec{\xi}\cdot\vec{x}} N^a Q_a(N) C_{g^N}({\psi}_1,\dots
,{\psi}_Z,\vec{\xi})dV_{g^n}
\\& =\int_{\mathbb{R}}
e^{(N-n)\vec{\xi}\cdot\vec{x}} [\Sigma_{u\in U} a_u
\frac{N^a}{(N-n)^{p_u}} Q_a(N)  C^u_{g^N}({\psi}_1,\dots
,{\psi}_Z)]dV_{g^n}
\end{split}
\end{equation}
 where either $a=p_u$ for every $u\in U$  or $a<p_u$ for every $u\in U$,
  depending whether $C_{g^N}({\psi}_1,\dots
,{\psi}_Z,\vec{\xi})$ is hard or easy, respectively. We then
denote the expression between brackets
 by $E[N^a Q_a(N)C_{g^N}({\psi}_1,\dots ,{\psi}_Z,\vec{\xi})]$.

\par As before, we break up
$E[N^a Q_a(N) C_{g^N}({\psi}_1,\dots ,{\psi}_Z,\vec{\xi})]$ into
two  sublinear combinations $E^{Intr}[N^a Q_a(N)
C_{g^N}({\psi}_1,\dots ,{\psi}_Z,\vec{\xi})]$, $E^{Unintr}[N^a
Q_a(N) C_{g^N}({\psi}_1,\dots ,{\psi}_Z,\vec{\xi})]$, that consist
of the intrinsic and un-intrinsic complete contractions,
respectively.

\par Now, in view of Lemma \ref{finecanc}, it follows that:

\begin{equation}
\label{giax}
\begin{split}
&{\Sigma}_{l\in L^1}a_l {\Sigma}_{b\in B^l} a_b N^{p_b}Q_b(N)
C^b_{g^N}({\psi}_1,\dots ,{\psi}_Z,\vec{\xi}) +
 {\Sigma}_{l\in L^2} a_l N^{p_l}Q_l(N) C^l_{g^N}({\psi}_1,\dots
,{\psi}_Z,\vec{\xi})=
\\& \Sigma_{w\in W} a_w N^{p_w}Q_w(N)C^w_{g^N}({\psi}_1,\dots
,{\psi}_Z,\vec{\xi})
\end{split}
\end{equation}
where each $\vec{\xi}$-contraction
$N^{p_w}Q_w(N)C^w_{g^N}({\psi}_1,\dots ,{\psi}_Z,\vec{\xi})$ is
easy, and moreover the rational function $Q_w(N)$ has degree zero
and leading order coefficient one. We deduce that:

\begin{equation}
\label{giax2}
\begin{split}
&{\Sigma}_{l\in L^1}a_l {\Sigma}_{b\in B^l} a_b E[N^{p_b}Q_b(N)
 C^b_{g^N}({\psi}_1,\dots ,{\psi}_Z,\vec{\xi})] +
 {\Sigma}_{l\in L^2} a_l E[N^{p_l}Q_l(N) \cdot
\\& C^l_{g^N}({\psi}_1,\dots
,{\psi}_Z,\vec{\xi})]=
 \Sigma_{w\in W} a_w E[N^{p_w}Q_w(N)C^w_{g^N}({\psi}_1,\dots
,{\psi}_Z,\vec{\xi})]
\end{split}
\end{equation}
and therefore:

\begin{equation}
\label{giax3}
\begin{split}
&{\Sigma}_{l\in L^1}a_l {\Sigma}_{b\in B^l} a_b
E^{Intr}[N^{p_b}Q_b(N) C^b_{g^N}({\psi}_1,\dots
,{\psi}_Z,\vec{\xi})] +
 {\Sigma}_{l\in L^2} a_l E^{Intr}[N^{p_l}Q_l(N)\cdot
\\& C^l_{g^N}({\psi}_1,\dots
,{\psi}_Z,\vec{\xi})]=
 \Sigma_{w\in W} a_w E^{Intr}[N^{p_w}Q_w(N)C^w_{g^N}({\psi}_1,\dots
,{\psi}_Z,\vec{\xi})]
\end{split}
\end{equation}

\par We then define a new operation $Oplim$ that acts on linear
combinations

\noindent $\Sigma_{h\in H} a_h
E[N^{p_l}Q_l(N)C^h_{g^N}({\psi}_1,\dots ,{\psi}_Z,\vec{\xi})]$,
(where $N^{p_l}Q_l(N)C^h_{g^N}({\psi}_1,\dots
,{\psi}_Z,\vec{\xi})$ may be either hard or undecided), by
rewriting them in dimension $n$ (thus the coefficients $N$ are now
independent of the dimension $n$) and letting $N\longrightarrow
\infty$. We act on the linear combinations on the left and right
hand sides of the above by the operation $Oplim$ and deduce:

\begin{equation}
\label{giax4}
\begin{split}
&{\Sigma}_{l\in L^1}a_l {\Sigma}_{b\in B^l} a_b
F^{Intr}[C^b_{g^N}({\psi}_1,\dots ,{\psi}_Z,\vec{\xi})] +
 {\Sigma}_{l\in L^2} a_l F^{Intr}[ C^l_{g^N}({\psi}_1,\dots
,{\psi}_Z,\vec{\xi})]=
\\& \Sigma_{w\in W} a_w Oplim\{ E^{Intr}[N^{p_w}Q_w(N)C^w_{g^N}({\psi}_1,\dots
,{\psi}_Z,\vec{\xi})]\}=0
\end{split}
\end{equation}

\par Thus, we indeed have that Lemma \ref{allo2ndref} follows
from Lemma \ref{finecanc}. $\Box$
\newline

{\it Proof of Lemma \ref{finecanc}:}
 Let rewrite (\ref{todo}) in the form:

$${\Sigma}_{l\in L} a_l (C^l)_{g^n}({\psi}_1,\dots ,
{\psi}_Z,\vec{\xi})=0$$

We will say $l\in L_\mu$ if and only if $C^l_{g^n}({\psi}_1,\dots
, {\psi}_Z,\vec{\xi})$ has $\mu>0$ factors $S\nabla^a\vec{\xi}$.
We prove the following: Let us suppose that for some $M>0$, we
have that for every $\mu> M$:

$${\Sigma}_{l\in L_\mu} a_l (C^l)_{g^n}({\psi}_1,\dots ,
{\psi}_Z,\vec{\xi})=0$$ We will then show that:

\begin{equation}
\label{tomeg}
{\Sigma}_{l\in L_M} a_l (C^l)_{g^n}({\psi}_1,\dots ,
{\psi}_Z,\vec{\xi})=0
\end{equation}

If we can show the above claim, our Lemma will follow by
induction. Now, recall that if for some linear combination of hard
$\vec{\xi}$-contractions we have (in $N$-cancelled notation) that
$$\Sigma_{r\in R} a_r C^r_{g^n}({\psi}_1,\dots ,
{\psi}_Z,\vec{\xi})=0$$

Then, by the argument above it follows that:

$$\Sigma_{r\in R} a_r F^{UnIntr}[C^r_{g^n}({\psi}_1,\dots ,
{\psi}_Z,\vec{\xi})]=0$$

Therefore, in view of our induction hypothesis, we will have that:

\begin{equation}
\label{gitchin} \Sigma_{\mu >M}\Sigma_{l\in L_\mu} a_l
F^{UnIntr}[(C^l)_{g^n}({\psi}_1,\dots , {\psi}_Z,\vec{\xi})]=0
\end{equation}
Therefore, for the proof of our inductive statement we may assume
that we have crossed out the above sublinear combination from
(\ref{reformunintrinsic}). Whenever we refer to
(\ref{reformunintrinsic}), we will be making that assumption.

\par In order to show (\ref{tomeg}), we will initially show:

\begin{equation}
\label{tomeg2} {\Sigma}_{l\in L_M} a_l
DF[(C^l)_{g^n}({\psi}_1,\dots , {\psi}_Z,\vec{\xi})]=0
\end{equation}

\par Firstly, for any complete contraction in the form (79), let us call
a factor $S\nabla^\nu_{r_1\dots r_\nu}\Gamma^k_{ij}$ where the
index $k$ contracts against one of the indices $r_1,\dots ,j$ a
{\it useful} factor.

\par Now, we notice that any complete contraction in
 (\ref{reformunintrinsic}) which does not belong to the sublinear
 combination (\ref{tomeg}) will have strictly less
 than $M$ useful factors.
This follows from our definition of the index set $L_M$ and
 Lemma \ref{fardtod}.

\par Now, let us denote by $Special(\Sigma_{l\in L_M})$ the
 sublinear combination in (\ref{tomeg2}) that
 consists of complete contractions all of whose
 factors in the form
$S\nabla^p\Gamma^k_{ij}$ satisfy $p\ge 1$. It follows that:

\begin{equation}
\label{avi}
Special(\Sigma_{l\in L_M})=0
\end{equation}
by substitution. But then, in view of Lemma \ref{corUMxi} we have that
(\ref{avi}) holds formally.
Then notice that under all the permutation identities in
 definition \ref{realformxi}, the number of factors
$S{\nabla}^p_{f_1\dots f_p}{\Gamma}^k_{ij}$ where the index
$k$ contracts against one of the indices $f_1,\dots ,j$
 remains invariant. Hence, since the left hand
 side of (\ref{tomeg2}) is the sublinear combination in
(\ref{reformunintrinsic}) with the maximum number of useful
factors, (\ref{tomeg2}) follows.

\par But then (\ref{tomeg2}) holds formally (again by Lemma
 \ref{corUMxi}). Hence,
let us imitate the permutations of factors that we do in
(\ref{tomeg2}) to make it formally zero for the
$\vec{\xi}$-contractions
 in (\ref{tomeg}). We only have to observe that if we can
 permute the indices of two tensors
$S{\nabla}^{m+1}_{d r_1\dots r_m}{\Gamma}^d_{ij}(x_0)$,
$S{\nabla}^{m+1}_{dr'_1\dots r'_m}{\Gamma}^d_{i'j'}(x_0)$
 to make them formally identical, we can then also permute
 the indices of the tensors $S{\nabla}^m_{r_1\dots r_m}
{\Gamma}^k_{ij}(x_0)\vec{\xi}_k$, $S{\nabla}^m_{r'_1\dots r'_m}
{\Gamma}^k_{i'j'}(x_0)\vec{\xi}_k$
to make them formally identical.
 We have shown our Lemma. $\Box$
\newline

We have proven our Proposition \ref{2ndref}. $\Box$

\subsection{The third refinement: The super divergence
 formula.}
\label{superdiv}

\par We begin this subsection with a few definitions.

\begin{definition}
\label{defstigma}
A $\vec{\xi}$-contraction $C_{g^n}({\psi}_1,\dots ,
{\psi}_Z,\vec{\xi})$ will be called stigmatized if the following holds: $C_{g^n}({\psi}_1,\dots ,
{\psi}_Z,\vec{\xi})$ is in the form (\ref{linisymxi2})
and each of its factors $\vec{\xi}$ contracts against another
 factor $\vec{\xi}$. We note that $C_{g^n}({\psi}_1,\dots ,
{\psi}_Z,\vec{\xi})$ is allowed to contain factors
$S{\nabla}^m\vec{\xi}$, $m\ge 1$.
\end{definition}

\par Now, consider any good or undecided complete contraction
$C_{g^n}({\psi}_1,\dots ,{\psi}_Z,\vec{\xi})$ and consider
 its iterative integration by parts.

\begin{definition}
\label{po}
We will define the pure outgrowth of
$C_{g^n}({\psi}_1,\dots ,{\psi}_Z)$.
We impose additional restrictions on the integration by parts:
Firstly, whenever we encounter a hard
$\vec{\xi}$-contraction we discard it. Secondly, whenever we
 encounter a stigmatized $\vec{\xi}$-contraction we also
 discard it.  Lastly, if we have a $\vec{\xi}$-contraction which
 is neither hard nor stigmatized, we will choose to integrate
  by parts with respect to a factor $\vec{\xi}$ that
 does not contract against another factor $\vec{\xi}$.

\par In the end, we will be left with a linear combination:

$${\int}_{\mathbb{R}^N}e^{(N-n)\vec{\xi}\cdot\vec{x}}
{\Sigma}_{h\in H} Q_h(N)\cdot a_h
C^h_{g^N}({\psi}_1,\dots ,{\psi}_Z)dV_{g^N}$$

Each complete contraction $C^h_{g^N}({\psi}_1,\dots ,
{\psi}_Z)$ is in the form (\ref{linisym}) and the rational
 function $Q_h(N)$ has degree $0$ and leading order
 coefficient $1$.

\par We define: $PO[C_{g^n}({\psi}_1,\dots ,{\psi}_Z)]=
 {\Sigma}_{h\in H} a_h
C^h_{g^n}({\psi}_1,\dots ,{\psi}_Z)$.
\end{definition}

\par Our goal for this subsection will be to show that:

\begin{proposition}
\label{bigcanc} If the sublinear combination of good and undecided
$\vec{\xi}$-contractions in $S^Z_{g^n}({\psi}_1,\dots ,
{\psi}_Z,\vec{\xi})$ is ${\Sigma}_{l\in L}a_l
C^l_{g^n}({\psi}_1,\dots ,{\psi}_Z,\vec{\xi})$, then:

\begin{equation}
\label{bigeqn}
I^Z_{g^n}({\psi}_1,\dots ,{\psi}_Z)+ {\Sigma}_{l\in L}
a_l PO[C^l_{g^n}({\psi}_1,\dots ,{\psi}_Z,\vec{\xi})]=0
\end{equation}
\end{proposition}

\par Before proving this proposition, we will need some preliminary Lemmas.

\begin{lemma}
\label{xiandr} Let us consider a good or undecided
$\vec{\xi}$-contraction $C_{g^n}({\psi}_1,\dots
,{\psi}_Z,\vec{\xi})$, in the form (\ref{linisymxi1}) or
(\ref{linisymxi2}). Let us suppose that $C_{g^n}({\psi}_1,\dots
,{\psi}_Z,\vec{\xi})$ has $\alpha$ factors $|\vec{\xi}|^2$ and
$\beta$ factors $R$ (scalar curvature). Let us consider the
iterative integration by parts (as in the previous subsection) of
$C_{g^n}({\psi}_1,\dots ,{\psi}_Z,\vec{\xi})$. Then, at each step
along the iterative integration by parts of
$C_{g^n}({\psi}_1,\dots ,{\psi}_Z,\vec{\xi})$,
 the number of factors $|\vec{\xi}|^2$ and the number of
 factors $R$ does not increase.
\end{lemma}

{\it Proof:} The proof is by induction, following the iterative integration by parts. $\Box$
\newline

\par We also define:

\begin{definition}
\label{xilength}
 Given any $\vec{\xi}$-contraction
$C_{g^n}({\psi}_1,\dots ,{\psi}_Z,\vec{\xi})$
in the general form (\ref{linisymxi2}), let $A$ be the
number of its factors ${\nabla}^mR_{ijkl}$, ${\nabla}^pRic$,
 $Z$ be the number of factors
${\nabla}^p{\psi}_l$, $C$ the number of its factors
 $S{\nabla}^m\vec{\xi}$ (with $m\ge 1$),
 $D$ the number of its factors
$|\vec{\xi}|^2$ and $E$ the number of its factors $\vec{\xi}$ that
do not contract against another factor $\vec{\xi}$. We then define
the $\vec{\xi}$-length of $C_{g^n}({\psi}_1,\dots ,{\psi}_Z)$ to
be $A+Z+C+D$. For any partial contraction in the form
(\ref{linisymxi1}) or (\ref{linisymxi2}), or any
$\vec{\xi}$-contraction with factors $\nabla^u\vec{\xi}$
(non-symmetrized), we define
its $\vec{\xi}$-length in the same way.
\end{definition}

\par We now want to see how any given complete contraction
$C_{g^n}({\psi}_1,\dots ,{\psi}_Z)$ in the form (\ref{linisym})
can give rise to good, undecided or hard $\vec{\xi}$-contractions
under the re-scaling $\hat{g}^N=e^{2\vec{\xi}\cdot\vec{x}}
g^N$.

\begin{definition}
\label{acceunacce}
 Let us consider any dimension-dependent
complete contraction $Q(N)\cdot C_{g^N}({\psi}_1,\dots
,{\psi}_Z,\vec{\xi})$, where $C_{g^N}({\psi}_1,\dots
,{\psi}_Z,\vec{\xi})$ is in the form (\ref{linisymxi2}), with
factors ${\nabla}^m\vec{\xi}$ instead of
 $S{\nabla}^m\vec{\xi}$. We will call such a
$\vec{\xi}$-contraction de-symmetrized. Recall that
$|\vec{\xi}|$ stands for the number of $\vec{\xi}$-factors. We
will call such a dimension-dependent $\vec{\xi}$-contraction {\it
acceptable}
 if $deg[Q(N)]=|\vec{\xi}|$ and {\it unacceptable} if
\\ $deg[Q(N)]<|\vec{\xi}|$.
\end{definition}

\par Now, consider any $C_{g^n}({\psi}_1,\dots ,{\psi}_Z)$, which is in the
 form (\ref{linisym}). We want to understand how the sublinear
 combination of acceptable complete contractions arises in
$e^{n\vec{\xi}\cdot\vec{x}}
C_{e^{2\vec{\xi}\cdot\vec{x}g^N}}({\psi}_1,\dots ,{\psi}_Z)$.
 We need one small convention before making our definition:
Whenever we have a factor ${\nabla}^m_{r_1\dots r_m}Ric_{ij}$ with
$m\ge 1$, we
 will assume that $i,j$ are not contracting between
 themselves. This can be done with no loss of generality by
 virtue of the formula ${\nabla}_aR=2{\nabla}^bRic_{ab}$.
Thus, we think of our complete contraction as being in the
 form:

\begin{equation}
\label{desym}
\begin{split}
&contr({\nabla}^{m_1}_{r_1\dots r_{m_1}}R_{ijkl}\otimes\dots
\otimes {\nabla}^{m_s}_{t_1\dots t_{m_s}}R_{ijkl}\otimes
\nabla^{p_1}_{r_1\dots r_{p_1}}Ric_{ij}\otimes\dots \otimes
\nabla^{p_q}_{t_1\dots t_{p_q}}Ric_{ij}\otimes \\& R^\alpha\otimes
 {\nabla}^{p_1}_{a_1\dots a_{p_1}}{\psi}_1\otimes \dots \otimes
{\nabla}^{p_Z}_{b_1\dots b_{p_Z}}{\psi}_Z)
\end{split}
\end{equation}

where the factors $\nabla^mR_{ijkl}$ do not have internal
contractions between the indices $i,j,k,l$, the factors
$\nabla^pRic_{ij}$ do not have internal contractions between the
indices $i,j$.

\par We are now ready for our definition.

\begin{definition}
\label{goodsub}
We consider internally contracted tensors in one of the
 following forms: ${\nabla}^p_{r_1\dots r_p}{\psi}_l$,
${\nabla}^p_{r_1\dots r_p}Ric_{ij}$, ${\nabla}^p_{r_1\dots
r_p}\vec{\xi}_j$ or ${\nabla}^p_{r_1\dots r_p}R_{ijkl}$. The
indices that are not internally contracted are considered
 to be free.

 We will call a pair of internally contracting indices, at
 least one of which is a derivative index,  an internal
 derivative contraction.
 We now want to define the {\it good substitutions} of each
tensor above.

\par For the tensor ${\nabla}^p_{r_1\dots r_p}
{\psi}_l$, we denote the pairs of internal contractions by
 $(r_{a_1},r_{b_1}),\dots (r_{a_l},r_{b_l})$.
The ordering of the indices $r_a, r_b$ in $(r_a, r_b)$ is
arbitrarily chosen. We define the set of {\it good
 substitutions} of the tensor ${\nabla}^p_{r_1\dots
r_p}{\psi}_l$ as follows: For any subset $\{ w_1,\dots ,
w_j\}\subset \{1,\dots ,l\}$ (including the empty set) the tensor
$\vec{\xi}^{r_{b_{w_1}}}\dots \vec{\xi}^{r_{b_{w_j}}}
{\nabla}^{p-j}_{r_1\dots \hat{r}_{a_{w_1}}\dots r_p}{\psi}_l$
 is a good substitution of ${\nabla}^p_{r_1\dots r_p}
{\psi}_l$.

\par We similarly define the set of good substitutions of any
tensor ${\nabla}^p_{r_1\dots r_p}Ric_{r_{p+1}r_{p+2}}$,
${\nabla}^p_{r_1\dots r_p}\vec{\xi}_{r_{p+1}}$ or
${\nabla}^p_{r_1\dots r_p}R_{r_{p+1}r_{p+2}r_{p+3}r_{p+4}}$
(this last is allowed to have internal contractions, but not among
 the set $r_{p+1},\dots ,r_{p+4}$):
 For any tensor above, let the
set of pairs of internal derivative contractions be
$(r_{a_1},r_{b_1}),\dots (r_{a_l},r_{b_l})$. The order of $r_a,
r_b$ in $(r_a,r_b)$ is arbitrarily chosen, but $r_a$ must be a
derivative index. Also, for the factor ${\nabla}^p_{r_1\dots r_p}
Ric_{r_{p+1}r_{p+2}}$, if $p\ge 1$, we assume that the indices
$r_{p+1},r_{p+2}$ do not contract against  each other.

\par Then, we  define the set of {\it good
 substitutions} of any tensor as above as follows:
For any subset $\{ w_1,\dots , w_j\}\subset \{1,\dots ,l\}$
(including the empty set), the tensor

\noindent $\vec{\xi}^{r_{b_{w_1}}}\dots \vec{\xi}^{r_{b_{w_j}}}
{\nabla}^{p-j}_{r_1\dots \hat{r}_{a_{w_1}}\dots r_p}
Ric_{r_{p+1}r_{p+2}}$ or
 $\vec{\xi}^{r_{b_{w_1}}}\dots \vec{\xi}^{r_{b_{w_j}}}
{\nabla}^{p-j}_{r_1\dots \hat{r}_{a_{w_1}}\dots r_p}
\vec{\xi}_{r_{p+1}}$ or

\noindent $e^{-2\vec{\xi}\cdot\vec{x}}\vec{\xi}^{b_{w_1}}\dots
\vec{\xi}^{b_{w_j}} {\nabla}^{p-j}_{r_1\dots
\hat{r}_{a_{w_1}}\dots r_p} R_{r_{p+1}r_{p+2}r_{p+3}r_{p+4}}$,
respectively, is a good substitution.
\end{definition}

\par We define any partial contraction
$C^{i_1\dots i_s}_{g^n}({\psi}_1,\dots ,{\psi}_Z,\vec{\xi})$ in
the form (\ref{linisymxi1}) or (\ref{linisymxi2}) to be
{\it nice} if in no
 factor $\vec{\xi}_i$ is the index $i$ free and no factor
$\vec{\xi}$ contracts against another factor $\vec{\xi}$
in $C^{i_1\dots i_s}_{g^n}({\psi}_1,\dots ,{\psi}_Z,
\vec{\xi})$.

\par We are now ready for our Lemma. We want to study the transformation law of any
$C_{g^n}({\psi}_1,\dots ,{\psi}_Z)$
 in the form (\ref{linisym}) under the re-scaling
$g^N\longrightarrow \hat{g}^N$. We do this in steps: Pick out any
factor $T^s_{a_1\dots a_j}$ in $C_{g^n}({\psi}_1,\dots ,{\psi}_Z)$
and make the indices $a_i$ that contract against any other factor
in $C_{g^n}({\psi}_1,\dots , {\psi}_Z)$ into free indices. So we
have a factor $(T^s_{a_1\dots a_j})_{a_{h_1}\dots a_{h_l}}$, which
we will call the {\it liberated form} of the factor $T^s_{a_1\dots
a_j}$. We view  $C_{g^n}({\psi}_1,\dots ,{\psi}_Z)$ as a complete
contraction among those tensors
 $T^s_{a_{h_1}\dots a_{h_l}}$. We
then consider each tensor $(T^s_{a_{h_1}\dots
a_{h_l}})^{(g^N)^{\vec{\xi}}}$. It will be a tensor of rank $l$.
 It follows that if we substitute each
$(T^s_{a_{h_1}\dots a_{h_l}})^{g^N}$ by $(T^s_{a_{h_1}\dots
a_{h_l}})^{\hat{g}^N}$ and then take the same complete contraction
of factors as for $C_{g^n}({\psi}_1,\dots ,{\psi}_Z)$, with
respect to the
 metric $(g^N)$, we will obtain
$e^{n\vec{\xi}\cdot\vec{x}}C_{e^{2\vec{\xi}\cdot\vec{x}}g^N}
({\psi}_1,\dots ,{\psi}_Z)$.

\begin{lemma}[The acceptable descendants]
\label{Itosupdiv}

\par Given a complete contraction
\\$C_{g^n}({\psi}_1,\dots ,{\psi}_Z)$ in the form
(\ref{desym}), the sublinear combination of
the acceptable $\vec{\xi}$-contractions in
$e^{n\vec{\xi}\cdot\vec{x}}C_{e^{2\vec{\xi}\cdot\vec{x}}g^N}
({\psi}_1,\dots ,{\psi}_Z)$
 arises as follows, in
$N$-cancelled notation:

\par Each of its liberated factors
$(T^s_{a_1\dots a_j})_{a_{h_1}\dots a_{h_l}}^{g^N}$ can be
substituted according to the pattern:

\begin{enumerate}
\item{
Any factor ${\nabla}^m
R_{ijkl}$ with no internal contractions must be left
 unaltered. }

\item{ Any factor of the form ${\nabla}^m_{r_1\dots r_m}
R_{ijkl}$
(where the indices $i,j,k,l$ do not contract between
 themselves) can be substituted by a good substitution of
${\nabla}^m_{r_1\dots r_m}R_{ijkl}$ or by a nice partial
contraction of $\vec{\xi}$-length $\ge 2$. }

\item{ Any factor ${\nabla}^p{\psi}_l$ can be substituted by a
good substitution
 of ${\nabla}^p{\psi}_l$ or by a nice
partial contraction of
$\vec{\xi}$-length $\ge 2$. }

\item{
Any factor ${\nabla}^p_{r_1\dots r_p}Ric_{ij}\ne R$ can
  be substituted either by a good substitution of
${\nabla}^p_{r_1\dots r_p}Ric_{ij}$
or a good substitution of $-{\nabla}^{p+1}_{r_1\dots r_pi}
\vec{\xi}_j$
 or by an nice
partial contraction of $\vec{\xi}$-length $\ge 2$.}

\item{
Any factor $R$ can be left unaltered or be
substituted by
$-2{\nabla}^i\vec{\xi}_i$ or by $-|\vec{\xi}|^2$.}
\end{enumerate}

We then claim that the sublinear combination of acceptable
$\vec{\xi}$-contractions in   $e^{n\vec{\xi}\cdot\vec{x}}
C_{e^{2\vec{\xi}\cdot\vec{x}}g^N}({\psi}_1,\dots ,{\psi}_Z)$
arises by substituting each liberated factor
$(T_{a_1\dots a_j})^{g^N}_{a_{h_1}\dots a_{h_l}}$
 in $C_{g^n}({\psi}_1,\dots ,{\psi}_Z)$ as explained above
 and then
performing the same particular contractions among the liberated
factors as in $C_{g^n}({\psi}_1,\dots , {\psi}_Z)$, with respect
to the metric $g^N$.
\end{lemma}

{\it Proof:} The proof of this Lemma is a matter of
applying formulas (\ref{levicivitaxi}), (\ref{scalarxi}),
(\ref{riccixi}) and (\ref{curvxi}) as well as
(\ref{symunsymgenxi}).
\newline

\par Consider any sequence of tensors times coefficients:
$a(N)\cdot(T_{i_1\dots i_j})^{g^N}$, where
$N=n,n+1,\dots$ and the tensors $(T_{i_1\dots i_j})^{g^n}$ are
partial contractions of the form:

\begin{equation}
\label{partcontra}
\begin{split}
&contr({\nabla}_{r_1\dots r_{m_1}}^{m_1}R_{i_1j_1k_1l_1}\otimes
\dots \otimes {\nabla}_{v_1\dots
v_{m_s}}^{m_s}R_{i_sj_sk_sl_s}\otimes
\\&{\nabla}^{{\nu}_1}_ {{\chi}_1\dots
{\chi}_{{\nu}_1}}{\psi}_l \otimes {\nabla}^m_{u_1\dots u_m}\vec{\xi}_z\otimes \dots \otimes
{\nabla}^m_{u_1\dots u_m}\vec{\xi}_z\otimes g^N_{ij}\otimes
\dots \otimes g^N_{ij})
\end{split}
\end{equation}
where there is at least one factor ${\nabla}^{\nu}{\psi}_l$ or
${\nabla}^mR_{ijkl}$ or ${\nabla}^m\vec{\xi}$, but not necessarily
one of each kind.  $a(N)$ is a rational function in $N$ and
$(T_{i_1\dots i_j})^{g^N}$ is the rewriting of $(T_{i_1\dots
i_j})^{g^n}$ in dimension $N$.

\par For any such partial contraction let $|g|$ stand for the number
of factors $g^N_{ij}$, $|\vec{\xi}|$ stand for the number of
factors ${\nabla}^m\vec{\xi}$ and $deg[a(N)]$ stand for the degree
of the rational function $a(N)$.

We also
consider linear combinations:

\begin{equation}
\label{usenouse}
{\Sigma}_{t\in T} a_t(N) (T^t_{i_1\dots i_s})^{g^N}
\end{equation}

where each sequence $a_t(N) (T^t_{i_1\dots i_s})^{g^N}$ is as
above. From now on we will just speak of the partial contraction
$a_t(N) (T^t_{i_1\dots i_s})^{g^N}$, rather than the sequence of
partial contractions times coefficients.

\par We say that such a partial contraction is {\it useful} if $|g|=0$,
$|\vec{\xi}|=deg[a_t(N)]$ and the index $k$ in each factor
$\vec{\xi}_k$ is not free and there are no factors
$|\vec{\xi}|^2$. We will call a partial contraction
 {\it useless} if $deg[a_t(N)]+|g|< |\vec{\xi}|$
or if $deg[a_t(N)]+|g|= |\vec{\xi}|$ and $|g|>0$. Note that
``useless'' is not the negation of ``useful''.
\newline

\par Consider any tensor
$({\nabla}^m_{r_1\dots r_m} R_{ijkl})^{g^N}$ or
${({\nabla}^p{\psi}_l)}^{g^N}$ or $(\nabla^p_{t_1\dots
t_p}Ric_{ij})^{g^N}$ with internal contractions. Let us suppose
that the free indices are $i_1,\dots i_s$. We will write those
tensors out as $({\nabla}^m_{r_1\dots r_m}
R_{ijkl})^{g^N}_{i_1,\dots i_s}$, ${({\nabla}^p
{\psi}_l)}^{g^N}_{i_1,\dots i_s}$, $(\nabla^p_{t_1\dots
t_p}Ric_{ij})^{g^N}_{i_1\dots i_s}$.
\newline

\noindent We claim that any tensor $e^{-2\vec{\xi}\cdot\vec{x}}
({\nabla}^m_{r_1\dots r_m}R_{ijkl})^{\hat{g}^N}_{i_1 \dots i_s}$
or ${({\nabla}^p{\psi}_l)_{i_1\dots i_s}}^{\hat{g}^N}$ or
\\ $(\nabla^p_{t_1\dots t_p}Ric_{ij}){i_1\dots
i_s}^{\hat{g}^N}$ (where $i,j,k,l$ in the first case and $i,j$ in
the second do not contract against each other) is a linear
combination of useful
 and useless tensors, as in (\ref{usenouse}).
Furthermore, we claim that each useful partial contraction of
$\vec{\xi}$-length 1 (in the expression for $({\nabla}^m_{r_1\dots
r_m}R_{ijkl})^{\hat{g}^N}_{i_1 \dots i_s}$ or
$({\nabla}^p{\psi}_l)_{i_1\dots i_s}^{\hat{g}^N}$ or
$(\nabla^p_{t_1\dots t_p}Ric_{ij})_{i_1\dots i_s}^{\hat{g}^N}$)
will be one of the good substitutions described in Definition
\ref{goodsub}. We refer to this claim as claim A.
\newline

We will check the above by an induction on $m$ or $p$,
 respectively. For $m=0$ or $p=1$, the fact is
 straightforward from (\ref{riccixi}) and (\ref{scalarxi}).
So, let us assume that we know that fact for $p=K$ or $m=K$
and show it for $p=K+1$ or $m=K+1$. Let us first consider the
 case of a tensor $({\nabla}^{K+1}_{r_1\dots r_{K+1}}
{\psi}_l)_{i_1\dots i_s}^{\hat{g}^N}$. We inquire whether the
index $r_1$ is free. If so, we then use our inductive hypothesis
for $p=K$: We know that the tensor $({\nabla}^K_{r_2\dots r_{K+1}}
{\psi}_l)_{i_2\dots i_s}^{\hat{g}^N}$ satisfies the induction
hypothesis. We want to use this fact to find
$({\nabla}^{K+1}_{r_1\dots r_{K+1}} {\psi}_l)_{i_1\dots
i_s}^{\hat{g}^N}$.

\par We write out:

\begin{equation}
\label{rita} ({\nabla}^K_{r_2\dots r_{K+1}} {\psi}_l)_{i_2\dots
i_s}^{\hat{g}^N}=\Sigma_{t\in T_1} a_t(N)
T^t_{g^N}(\psi_l,\vec{\xi})_{i_2\dots i_s} + \Sigma_{t\in T_2}
a_t(N) T^t_{g^N}(\psi_l,\vec{\xi})_{i_2\dots i_s}
\end{equation}
where the first sublinear combination stands for the useful
tensors and the second stands for the  useless tensors.

We only have to apply the transformation law (\ref{levicivitaxi})
to each pair $(r_1,i_2),\dots ,(r_1,i_s)$. We easily observe
 that if any summand in the expression of
$({\nabla}^K_{r_2\dots r_{K+1}} {\psi}_l)_{i_2\dots
i_s}^{\hat{g}^N}$ is useless, then any application of the identity
(\ref{levicivitaxi}) to any pair of indices $(r_1,i_2),\dots ,
(r_1,i_s)$ will give rise to a useless partial contraction.

\par On the other hand, consider any factor
$T^t_{g^N}(\psi_l,\vec{\xi})_{i_2\dots i_s}$ in
$({\nabla}^K_{r_2\dots r_{K+1}} {\psi}_l)_{i_2\dots
i_s}^{\hat{g}^N}$ which is
 useful. We then observe that when we apply any of the last three
 summands in
 (\ref{levicivitaxi}) to any pair of indices $(r_1,i_2),
\dots ,(r_1,i_s)$ and bring out a factor $\vec{\xi}$, we
 obtain a useless partial contraction. Finally, if we substitute
$({\nabla}_{r_1}X_{i_l})^{\hat{g}^N}$ by
$({\nabla}_{r_1}X_{i_l})^{(g^N)}$ (the first summand on the right
hand side of (\ref{levicivita}),
 we will get a linear combination of useful
$\vec{\xi}$-contractions, by applying the rule

$${\nabla}_i[A_{k_1\dots k_s}\otimes B_{u_1\dots u_h}]=
{\nabla}_iA_{k_1\dots k_s}\otimes B_{u_1\dots u_h}+
A_{k_1\dots k_s}\otimes {\nabla}_i B_{u_1\dots u_h}$$

 Furthermore,
we observe that if a partial contraction $a_t(N)T^t_{g^N}
(\psi_l,\vec{\xi})$ in $({\nabla}^K_{r_2\dots r_{K+1}}
{\psi}_l)_{i_2\dots i_s}^{\hat{g}^N}$ does not contain factors
$\vec{\xi}_k$ where the index $k$ is free, nor factors
$|\vec{\xi}|^2$, then we will have no such factors in
${\nabla}_{r_1}a_t(N)T^t_{g^N}(\psi_l,\vec{\xi})$ either.

\par Finally, any partial contraction in
${\nabla}_{r_1}a_t(N)T^t_{g^N}(\psi_l,\vec{\xi})$ of
$\vec{\xi}$-length 1 will arise if $T^t$ has $\vec{\xi}$-length 1
and provided the derivative does not hit any factor $\vec{\xi}_k$
. So, by our inductive hypothesis, we observe
 that any useful partial contraction in
${({\nabla}^{K+1}_{r_1\dots r_{K+1}}
{\psi}_l)_{i_1\dots i_s}}^{(g^N)^{\vec{\xi}}}$
of $\vec{\xi}$-length $1$ is a good substitution.
\newline

\par Next, let us consider the case where the index $r_1$ in
$({\nabla}^{K+1}_{r_1\dots r_{K+1}} {\psi}_l)_{i_1\dots
i_s}^{\hat{g}^N}$ is not a free index. Let us suppose that $r_1$
contracts against $r_j$. We consider the tensor
$({\nabla}^K_{r_2\dots r_{K+1}} {\psi}_l)_{i_2\dots r_j\dots
i_s}^{\hat{g}^N}$ which is obtained from
${({\nabla}^{K+1}_{r_1\dots r_{K+1}} {\psi}_l)_{i_1\dots
i_s}}^{(g^N)^{\vec{\xi}}}$ by erasing the derivative
${\nabla}_{r_1}$
 and making the index $r_j$ into a free index. We consider
 the transformation law for $({\nabla}^K_{r_2\dots r_{K+1}}
{\psi}_l)_{i_1\dots r_j\dots i_s}^{\hat{g}^N}$. Our inductive
hypothesis applies. So, in order to determine
$({\nabla}^{K+1}_{r_1\dots r_{K+1}} {\psi}_l)_{i_1\dots
i_s}^{\hat{g}^N}$, we have to apply (\ref{levicivitaxi}) to each
pair $(r_1, i_2),\dots , (r_1,i_s), (r_1,r_j)$ and then contract
$r_1$ and $r_j$.
 We observe that if we consider any partial contraction in
$({\nabla}^K_{r_2\dots r_{K+1}} {\psi}_l)_{i_1\dots r_j\dots
i_s}^{\hat{g}^N}$ which is useless, then any application of the
law (\ref{levicivitaxi})
 to any above pair will give us a useless partial contraction.

\par Now, let us consider any useful partial contraction
$a_t(N)\cdot T^t_{i_1\dots r_j\dots i_s}$ in
 $({\nabla}^K_{r_2\dots r_{K+1}}
{\psi}_l)_{i_1\dots r_j\dots i_s}^{\hat{g}^N}$. We notice that if
we apply the identity (\ref{levicivitaxi})
 to any pair of indices $(r_1, i_2),\dots ,(r_1,i_s),
(r_1,r_j)$ without bringing out a factor $\vec{\xi}$ (meaning that
we apply the first summand on the right hand side of
(\ref{levicivitaxi})), then by
 the same reasoning as before we have our claim. On
 the other hand, if we apply the identity (\ref{levicivitaxi})
 to any pair of indices $(r_1, i_1),\dots ,(r_1,i_s)$ and we
bring out a factor $\vec{\xi}$, then after we contract
$r_1,r_j$ we will obtain a useless partial contraction.
Also, if we apply the transformation law (\ref{levicivitaxi})
 to the pair $(r_1,r_j)$ and
 bring out a factor $\vec{\xi}$ but not a factor $g_{ij}$,
 then after we contract $r_1, r_j$ we will again obtain a
 useless partial contraction. Finally, if we apply the
 transformation law (\ref{levicivitaxi}) to
$({\nabla}_{r_1}X_{r_j})^{\hat{g}^N}$ and bring out
 $g^N_{r_1 r_j}\vec{\xi}^sX_s$, then after we contract
$r_1, r_j$ we will bring out a factor $N$. We observe that we thus
obtain another useful $\vec{\xi}$-contraction.

\par Finally, notice that if $a_t(N)\cdot T^t_{i_1\dots r_j\dots i_s}$ had
$\vec{\xi}$-length $1$, then by our inductive hypothesis it
 was a good substitution of
$({\nabla}^p{\psi}_l)_{i_1\dots r_j\dots i_s}$. Hence, for each
such good substitution, we now have the option of either
substituting ${\nabla}^{r_1}X_{r_j}$ by $N\vec{\xi}^sX_s$ or
leaving it unaltered. Therefore, we see that the set of useful
$\vec{\xi}$-contractions of $\vec{\xi}$-length 1 in
$({\nabla}^{K+1}_{r_1\dots r_{K+1}} {\psi}_l)_{i_1\dots
i_s}^{\hat{g}^N}$ is
 indeed contained in the set of good substitutions of
$({\nabla}^{K+1}_{r_1\dots r_{K+1}}
{\psi}_l)_{i_1\dots i_s}$. Moreover, we again observe that since any
 useful tensor $a_t(N) T^t_{i_1\dots r_j\dots i_s}(\psi_l,\vec{\xi})$
in $({\nabla}^K_{r_2\dots r_{K+1}} {\psi}_l)_{i_1\dots r_j\dots
i_s}^{\hat{g}^N}$ does not have factors $\vec{\xi}$ or
$|\vec{\xi}|^2$, then there will
 be no such factors in either
$a_t(N) \nabla^{r_j}T^t_{i_1\dots r_j\dots i_s}(\psi_l,\vec{\xi})$
or $a_t(N) \vec{\xi}^{r_j}T^t_{i_1\dots r_j\dots i_s}
(\psi_l,\vec{\xi})$. Hence, we have completely shown our
 inductive step.

\par The same proof applies to show our claim for the tensors
$({\nabla}^m_{r_1\dots r_m} R_{ijkl})^{\hat{g}^N}_{i_1\dots i_s}$
and $({\nabla}^p_{r_1\dots r_p} R_{ij})^{\hat{g}^N}_{i_1\dots
i_s}$. The cases $m=0, p=0$ follow by equations (\ref{curvxi}),
(\ref{riccixi}). We then see that the same proof as above
 still
 applies, since it is only an iterative application of the formula
(\ref{levicivitaxi}). We have proven claim A.
\newline

\par Now, in order to complete the proof of our Lemma, we only have to
observe that if we substitute any liberated factor $T\ne R$ from
$C_{g^n}({\psi}_1,\dots ,{\psi}_Z)$
 by a useless partial contraction, and then proceed to replace the
other factors by either useful or useless partial contractions and
then perform the same contractions for those replacements as for
$C_{g^N}({\psi}_1,\dots , {\psi}_Z)$, we will obtain an
unacceptable complete contraction in
$e^{n\vec{\xi}\cdot\vec{x}}C_{e^{2\vec{\xi}\cdot\vec{x}}g^N}
({\psi}_1,\dots , {\psi}_Z)$. This follows by the same reasoning
as for Lemma \ref{dimdepindep}. Regarding the substitutions of
scalar curvature we note that we can replace it by either a factor
$-(2-N)\nabla^a\vec{\xi}_a$ (in which case $deg[-(2-N)]=1,
|\vec{\xi}|=1$) or by $-(N-1)(N-2)|\vec{\xi}|^2$ (in which case
$deg[-(N-1)(N-2)]=2, |\vec{\xi}|=2$). $\Box$
\newline

\par Let us now state a corollary of the above Lemma regarding
the linear combination of good, hard and undecided
descendants of a complete contraction
$C_{g^n}({\psi}_1,\dots ,{\psi}_Z)$, in the form
(\ref{linisym}).

\par We consider any complete contraction
$C_{g^n}({\psi}_1,\dots ,{\psi}_Z)$, in the form (\ref{linisym}), and write it in the form (\ref{desym}).
We then consider the sublinear combination of its
 acceptable descendants, in the form:

\begin{equation}
\label{desymxi}
\begin{split}
&contr({\nabla}^{m_1}_{r_1\dots r_{m_1}}R_{ijkl}\otimes\dots
\otimes {\nabla}^{m_s}_{t_1\dots t_{m_s}}R_{ijkl}\otimes
\\& \nabla^{p_1}_{r_1\dots r_{p_1}}Ric_{ij}\otimes\dots \otimes
\nabla^{p_q}_{t_1\dots t_{p_q}}Ric_{ij}\otimes R^\alpha\otimes
\\& {\nabla}^{p_1}_{a_1\dots a_{p_1}}{\psi}_1\otimes
\dots \otimes {\nabla}^{p_Z}_{b_1\dots b_{p_Z}}{\psi}_Z
\otimes\nabla^{b_1}\vec{\xi}\otimes \dots\otimes
\nabla^{b_v}\vec{\xi})
\end{split}
\end{equation}

 Then, by repeated application of the formula
(\ref{symunsymgenxi}), we write each such de-symmetrized
 descendant as a linear combination of good, hard and undecided
 $\vec{\xi}$-contractions of the form
(\ref{linisymxi1}) or (\ref{linisymxi2}).

\par We then claim the following:

\begin{lemma}
\label{apogonoi}
\par Given any complete contraction
$C_{g^n}({\psi}_1,\dots ,{\psi}_Z)$ in the form
(\ref{linisym}), of length $L$,
then each of its good or hard or undecided
descendants constructed above
will have $\vec{\xi}$-length $\ge L$.

\par Furthermore, if the complete contraction
$C_{g^n}({\psi}_1,\dots ,{\psi}_Z)$
has no factors $R$, then none of its descendants will
contain a factor $-|\vec{\xi}|^2$. On the other hand, if
$C_{g^n}({\psi}_1,\dots ,{\psi}_Z)$
contains $A>0$ factors $R$, then we can write the
sublinear combination of its good, undecided and
hard descendants as follows:

\begin{equation}
\label{sxesh}
{\Sigma}_{l\in L} a_l[C^l_{g^n}({\psi}_1,\dots ,{\psi}_Z,
\vec{\xi})+ {\Sigma}_{r\in R^l} C^r_{g^n}({\psi}_1,\dots ,{\psi}_Z,\vec{\xi})]
\end{equation}

where each $\vec{\xi}$-contraction
$C^l_{g^n}({\psi}_1,\dots ,{\psi}_Z,\vec{\xi})$
arises from $C_{g^n}({\psi}_1,\dots ,{\psi}_Z)$
by doing all the substitutions explained in Lemma
 \ref{Itosupdiv} but leaving all the factors $R$ unaltered.
The linear combination ${\Sigma}_{r\in R^l} C^r_{g^n}({\psi}_1,
\dots ,{\psi}_Z,\vec{\xi})$
arises from
\\ $C^l_{g^n}({\psi}_1,\dots ,{\psi}_Z,\vec{\xi})$ by
 substituting a nonzero number of factors $R$
 by either $-2{\nabla}^i\vec{\xi}_i$ or
$-|\vec{\xi}|^2$ and then summing over all those different
 substitutions.
\end{lemma}

{\it Proof:} This Lemma follows straightforwardly from Lemma
\ref{Itosupdiv}: We only have to make note that $\vec{\xi}$-length
is additive and that the correction terms that we introduce in the
symmetrization of factors $\nabla^p\vec{\xi}$
(using (\ref{symunsymgenxi})) may increase
the $\vec{\xi}$-length but not decrease it. So, since we
 are substituting each factor in
$C_{g^n}({\psi}_1,\dots ,{\psi}_Z)$ by a tensor of
$\vec{\xi}$-length $\ge 1$, the first claim of our Lemma will
follow.
\newline

\par Our second claim will follow from the transformation law
(\ref{scalarxi}), provided we can show that no factors
$|\vec{\xi}|^2$ arise when we symmetrize and anti-symmetrize the
factors $\nabla^p\vec{\xi}$ and then repeat the same particular
contractions as for $C_{g^n}(\psi_1,\dots ,\psi_s)$. In order to
see this, we only have to observe that for each factor of the form
$\nabla^p\vec{\xi}_j$, $p\ge 1$, none of the correction terms
involve a factor $\vec{\xi}_a$ with the index $a$ being free.
 We can see this if we can show that in each correction term
 from the
 symmetrization of $\nabla^p\vec{\xi}$, there is no factor
$\vec{\xi}_a$ where the index $a$ is free.

 This follows
because in order to symmetrize the factor $\nabla^p\vec{\xi}_j$ we only
use the identities $[\nabla_a\nabla_b-\nabla_b\nabla_a]\vec{\xi}_j=
R_{abjd}\vec{\xi}^d$ and, if $k\ge 1$:

$$\nabla^u\{ [\nabla_a\nabla_b-\nabla_b\nabla_a]\nabla^k\vec{\xi}\}=
\Sigma(\nabla^tR\nabla^y\vec{\xi})$$ where
$\Sigma(\nabla^tR\nabla^y\vec{\xi})$ stands for a linear
combination of partial contractions of the form $\nabla^\alpha
R_{ijkl}\nabla^y\vec{\xi}$, where $1\le y<k+u+2$. $\Box$
\newline

{\it The proof of Proposition \ref{bigcanc}.}
\newline

 Recall that we have defined a stigmatized
$\vec{\xi}$-contraction to be in the form:

\begin{equation}
 \label{stigma}
\begin{split}
&contr({\nabla}_{r_1\dots r_{m_1}}^{m_1}R_{i_1j_1k_1l_1}\otimes
\dots \otimes {\nabla}_{v_1\dots
v_{m_s}}^{m_s}R_{i_sj_sk_sl_s}\otimes \\& {\nabla}_{t_1 \dots
t_{p_1}}^{p_1}Ric_{{\alpha}_1 {\beta}_1} \otimes \dots \otimes
{\nabla}_{z_1 \dots z_{p_q}}^{p_q} Ric_{{\alpha}_q{\beta}_q}
\otimes{\nabla}^{{\nu}_1}_ {{\chi}_1\dots
{\chi}_{{\nu}_1}}{\psi}_1 \otimes\dots \otimes
{\nabla}^{{\nu}_Z}_{{\omega}_1\dots {\omega}_{{\nu}_Z}}
{\psi}_Z
\\& \otimes
S{\nabla}^{{\mu}_1}\vec{\xi}_{j_1}\otimes \dots\otimes
S{\nabla}^{{\mu}_r}\vec{\xi}_{j_s}\otimes |\vec{\xi}|^2
\otimes \dots \otimes |\vec{\xi}|^2)
\end{split}
\end{equation}

  where each ${\mu}_i\ge 1$ and
there are $r$ factors $S{\nabla}^{\nu}\vec{\xi}$ and
 $s>0$ factors $|\vec{\xi}|^2$. If $r=0$, we will call the
the above $\vec{\xi}$-contraction stigmatized of type 1 and
 if $r>0$, we will call it stigmatized of type 2.
\newline

\par Let us, for each good or undecided
$\vec{\xi}$-contraction $C^l_{g^n}(\phi,\vec{\xi})$ (with $X$
$\vec{\xi}$-factors) break up
 its outgrowth $O[C^l_{g^n}(\phi,\vec{\xi})]$ as follows:
After each integration by parts of a factor $\vec{\xi}$, we
 discard any hard $\vec{\xi}$-contractions that arise, but
moreover, when we encounter any stigmatized complete contractions of type 1 or
 type 2 we put them aside.
We denote by ${\Sigma}_{k\in K^l_1} a_k C^k_{g^n}({\psi}_1, \dots
,{\psi}_Z,\vec{\xi})$ the sublinear combination of
$\vec{\xi}$-contractions that we are left with after $X-1$
integrations by parts, after we have discarded all the hard
$\vec{\xi}$-contractions we encounter and after we have put aside
all the stigmatized $\vec{\xi}$-contractions we encounter. We also
denote by ${\Sigma}_{k\in K^l_2} a_k C^k_{g^n}({\psi}_1, \dots
,{\psi}_Z,\vec{\xi})$, ${\Sigma}_{k\in K^l_3} a_k
C^k_{g^n}({\psi}_1, \dots ,{\psi}_Z,\vec{\xi})$ the sublinear
combinations of stigmatized $\vec{\xi}$-contractions of types 1
and 2, respectively, that we have put aside along our iterative
integrations by parts.

  We will then have:

\begin{equation}
\label{outbreak}
\begin{split}
&O[C^l_{g^n}(\phi,\vec{\xi})]= {\Sigma}_{k\in K^l_1} a_k
O[C^k_{g^n}({\psi}_1, \dots ,{\psi}_Z,\vec{\xi})] + {\Sigma}_{k\in
K_2^l} a_k O[C^k_{g^n}({\psi}_1, \dots ,{\psi}_Z,\vec{\xi})]+
\\& {\Sigma}_{k\in K_3^l} a_k O[C^k_{g^n}({\psi}_1,
\dots ,{\psi}_Z,\vec{\xi})]
\end{split}
\end{equation}

We observe that the $\vec{\xi}$-contractions $C^k_{g^n}({\psi}_1,
\dots ,{\psi}_Z)$, $k\in K_1^l$ are good, in the from
(\ref{linisymxi1}) with one factor $\vec{\xi}$.

\par Hence, we can rewrite (\ref{Ointrinsic}) as follows:

\begin{equation}
\label{reintr}
\begin{split}
& I^Z_{g^n}(\phi)+ {\Sigma}_{l\in L} a_l \{
{\Sigma}_{k\in K_l^1} a_k O[C^k(\phi,\vec{\xi})]+
\\& {\Sigma}_{k\in K_l^2} a_k O[C^k(\phi,\vec{\xi})]
+ {\Sigma}_{k\in K_l^3} a_k O[C^k(\phi,\vec{\xi})] \}=0
\end{split}
\end{equation}

\par Our Proposition will follow from the following equation:

\begin{equation}
\label{Rcanc}
{\Sigma}_{l\in L} a_l \{ {\Sigma}_{k\in K^l_2} a_k O[C^k_{g^n}({\psi}_1,
\dots ,{\psi}_Z,\vec{\xi})]+
{\Sigma}_{k\in K^l_3} a_k O[C^k_{g^n}({\psi}_1,
\dots ,{\psi}_Z,\vec{\xi})]\} =0
\end{equation}

\par In fact, we will show that:

\begin{equation}
\label{Laplbreak1}
{\Sigma}_{l\in L} a_l \{{\Sigma}_{k\in K^l_2} a_k C^k_{g^n}({\psi}_1,
\dots ,{\psi}_Z,\vec{\xi})\}=0
\end{equation}

and

\begin{equation}
\label{Laplbreak2}
{\Sigma}_{l\in L} a_l\{ {\Sigma}_{k\in K^l_3} a_k C^k_{g^n}({\psi}_1,\dots ,{\psi}_Z,\vec{\xi})\}=0
\end{equation}

\par We see that (\ref{Rcanc}) follows from the above two
equations by the same reasoning by which Lemma \ref{allo2ndref}
follows from Lemma \ref{finecanc}.

\par We first show (\ref{Laplbreak2}).
Let us define a procedure which we will call the {\it sieving}
 integration by parts.

\begin{definition}
\label{sieve}
Consider any good or undecided $\vec{\xi}$-contraction
$C_{g^n}({\psi}_1,\dots ,{\psi}_Z,\vec{\xi})$. We consider its
iterative integrations by parts and we impose the
 following rules:  Whenever along the iterative integration by
 parts we encounter a hard $\vec{\xi}$-contraction, we
erase it and we put it in the linear combination
$H[C_{g^n}({\psi}_1,\dots ,{\psi}_Z,\vec{\xi})]$. Whenever we
 encounter a $\vec{\xi}$-contraction which is stigmatized of type 2,
 we erase it and put it in the linear combination
 $Stig^2[C_{g^n}({\psi}_1,\dots ,{\psi}_Z,\vec{\xi})]$.
 Also, whenever we encounter a stigmatized $\vec{\xi}$-contraction
 of type 1, we erase it and put it in the linear combination
 $Stig^1[C_{g^n}({\psi}_1,\dots ,{\psi}_Z,\vec{\xi})]$.

\par Furthermore, consider any complete contraction
$C_{g^n}({\psi}_1,\dots ,{\psi}_Z)$ which is in the form
(\ref{linisym}). We consider the linear combination of its good
or undecided or hard descendants, say
${\Sigma}_{d\in D} a_d C^d_{g^n}({\psi}_1,\dots ,{\psi}_Z,
\vec{\xi})$. We define
$$H[C_{g^n}({\psi}_1,\dots ,{\psi}_Z)]={\Sigma}_{d\in D} a_d
H[C^d_{g^n}({\psi}_1,\dots ,{\psi}_Z,\vec{\xi})]$$
We also define:

$$PO[C_{g^n}({\psi}_1,\dots ,{\psi}_Z)]=
{\Sigma}_{d\in D} a_d PO[C_{g^n}^d({\psi}_1,\dots ,{\psi}_Z,
\vec{\xi})]$$

$$Stig^2[C_{g^n}({\psi}_1,\dots ,{\psi}_Z)]=
{\Sigma}_{d\in D} a_d Stig^2[C_{g^n}^d({\psi}_1,\dots ,{\psi}_Z,
\vec{\xi})]$$

$$Stig^1[C_{g^n}({\psi}_1,\dots ,{\psi}_Z)]=
{\Sigma}_{d\in D} a_d Stig^1[C^d_{g^n}({\psi}_1,\dots ,{\psi}_Z,\vec{\xi})]$$
\end{definition}

\par We claim the following:

\begin{lemma}
\label{htostig}
Let us consider any complete contraction
$C_{g^n}({\psi}_1,\dots ,{\psi}_Z)$ in the form
(\ref{linisym}) of weight $-n$. We claim that there is a way to perform our sieving integration by parts, so that we can
express $PO[C_{g^n}({\psi}_1,\dots ,{\psi}_Z)]$,
$H[C_{g^n}({\psi}_1,\dots ,{\psi}_Z)]$,
$Stig^1[C_{g^n}({\psi}_1,\dots ,{\psi}_Z)]$ and
$Stig^2[C_{g^n}({\psi}_1,\dots ,{\psi}_Z)]$ as follows:

\begin{equation}
\label{toagno}
C_{g^n}({\psi}_1,\dots ,{\psi}_Z)+ PO[C_{g^n}({\psi}_1,\dots ,{\psi}_Z)]={\Sigma}_{v\in V}
a_v C^v_{g^n}({\psi}_1,\dots ,{\psi}_Z)R^{{\alpha}_v}
\end{equation}

where each $C^v_{g^n}({\psi}_1,\dots ,{\psi}_Z)$
is of weight $-n+2{\alpha}_v$, in the
 form (\ref{linisym}), with no factors $R$
(they are pulled out on the right). We then claim that
$H[C_{g^n}({\psi}_1,\dots ,{\psi}_Z)]$ can be expressed as follows:

\begin{equation}
\label{h}
\begin{split}
&H[C_{g^n}({\psi}_1,\dots ,{\psi}_Z)]=
{\Sigma}_{v\in V}
a_v C^v_{g^n}({\psi}_1,\dots ,{\psi}_Z)\cdot G(R,{\alpha}_v,
-2{\nabla}^i
\vec{\xi}_i)+
\\& {\Sigma}_{f\in F} a_f C^f_{g^n}({\psi}_1,\dots ,{\psi}_Z,
\vec{\xi})\cdot R^{{\alpha}_f}+
{\Sigma}_{f\in F} a_f C^f_{g^n}({\psi}_1,\dots ,{\psi}_Z,
\vec{\xi})\cdot G(R, {\alpha}_f, -2{\nabla}^i\vec{\xi}_i)
\end{split}
\end{equation}

where each $C^f_{g^n}({\psi}_1,\dots ,{\psi}_Z,\vec{\xi})$
is of weight $-n+2{\alpha}_f$, in the
 form (\ref{linisymxi2}) with $k=0$ and with no factors $R$
(they are pulled out on the right). The symbol
$G(R,\lambda, B)$ stands for the sum over
all the possible substitutions of $\lambda$ factors $R$ by a factor
$B$, so that we make at least one such substitution.

\par Finally, we claim that
$Stig^1[C_{g^n}({\psi}_1,\dots ,{\psi}_Z)]$ and
$Stig^2[C_{g^n}({\psi}_1,\dots ,{\psi}_Z)]$ can be
 expressed as follows:

\begin{equation}
\label{stig1}
Stig^1[C_{g^n}({\psi}_1,\dots ,{\psi}_Z)]=
{\Sigma}_{v\in V}
a_v C^v_{g^n}({\psi}_1,\dots ,{\psi}_Z)\cdot
G(R,{\alpha}_v, -|\vec{\xi}|^2)
\end{equation}

\begin{equation}
\label{stig2}
\begin{split}
&Stig^2[C_{g^n}({\psi}_1,\dots ,{\psi}_Z)]=
{\Sigma}_{v\in V}
a_v C^v_{g^n}({\psi}_1,\dots ,{\psi}_Z)\cdot T^{*}({\alpha}_v,R
,-2{\nabla}^i\vec{\xi}_i, -|\vec{\xi}|^2)+
\\&{\Sigma}_{f\in F} a_f C^f_{g^n}({\psi}_1,\dots ,{\psi}_Z,
\vec{\xi})\cdot T({\alpha}_f, R,-2{\nabla}^i\vec{\xi}_i,
-|\vec{\xi}|^2)
\end{split}
\end{equation}

where $T(j,R, -2{\nabla}^i\vec{\xi}_i,
-|\vec{\xi}|^2)$
stands for the sum over all the possible ways to substitute a nonzero
number of factors in $R^j$ by either
$-2{\nabla}^i\vec{\xi}_i$ or
$-|\vec{\xi}|^2$, so that at least one factor is substituted
 by $-|\vec{\xi}|^2$. $T^{*}(j,R, -2{\nabla}^i\vec{\xi}_i,
-|\vec{\xi}|^2)$ stands for the same thing, with the additional
 restriction that at least one factor $R$ must be substituted by $-2{\nabla}^i\vec{\xi}_i$.
\end{lemma}

{\it Proof:} In order to prove the above, we will consider
 the linear combination of
$C_{g^n}({\psi}_1,\dots ,{\psi}_Z)$, together with its
good, undecided and hard
 descendants in

\noindent $e^{N\vec{\xi}\cdot\vec{x}}C_{e^{2\vec{\xi}\cdot
\vec{x}}g^n}({\psi}_1,\dots ,{\psi}_Z)$,
 grouped up as in (\ref{sxesh}).
Given any $l\in L$, we pick any
$C^r_{g^n}({\psi}_1,\dots ,{\psi}_Z,\vec{\xi})$, $r\in R^l$
 and we identify any factor $T\ne -2{\nabla}^i\vec{\xi}_i,
-|\vec{\xi}|^2$ in
$C^r_{g^n}({\psi}_1,\dots ,{\psi}_Z,\vec{\xi})$
with a factor in
$C^l_{g^n}({\psi}_1,\dots ,{\psi}_Z,\vec{\xi})$. We then say
that such a factor in $C^r_{g^n}({\psi}_1,\dots ,{\psi}_Z,
\vec{\xi})$, $r\in R^l$
 {\it corresponds} to a factor in
$C^l_{g^n}({\psi}_1,\dots ,{\psi}_Z,\vec{\xi})$.

\par We will now perform integrations by parts among the
sublinear combination of good, hard and undecided descendants
in $\int_{\mathbb{R}^N}e^{(N-n)\vec{\xi}\cdot\vec{x}}
C^a_{e^{2\vec{\xi}\cdot\vec{x}}g^N}
({\psi}_1,\dots ,{\psi}_Z)dV_{g^N}$, so that after any number
 of integrations by parts we will be left with an integrand
 of $\vec{\xi}$-contractions as in (\ref{sxesh}).

\par Now, for any
$C^l_{g^n}({\psi}_1,\dots ,{\psi}_Z,\vec{\xi})$, we pick out
 a factor $\vec{\xi}_i$ (which does not contract against
 another $\vec{\xi}$) and perform an integration by parts.
We will obtain a formula:

\begin{equation}
\label{intbp}
\begin{split}
&{\int}_{\mathbb{R}^N}e^{(N-n)\vec{\xi}\cdot\vec{x}}Q(N)C^l_{g^N}
({\psi}_1,\dots ,{\psi}_Z,\vec{\xi})dV_{g^N}=
\\&{\int}_{\mathbb{R}^N}e^{(N-n)\vec{\xi}\cdot\vec{x}}
\frac{Q(N)}{N-n}[{\Sigma}_{\alpha =1}^L C^{l,\alpha}_{g^N}
({\psi}_1,\dots ,{\psi}_Z,\vec{\xi})]dV_{g^N}
\end{split}
\end{equation}

 Consider any $\vec{\xi}$-contraction $C^{l,\alpha}_{g^N}
({\psi}_1,\dots ,{\psi}_Z,\vec{\xi})$
which arises when ${\nabla}_i$ hits a
 factor $T$ in $C^l_{g^n}({\psi}_1,\dots ,
{\psi}_Z,\vec{\xi})$ with $T\ne R$. Then consider  any
$\vec{\xi}$-contraction $C^r_{g^n}({\psi}_1,\dots ,
{\psi}_Z,\vec{\xi})$, $r\in R^l$ and integrate by parts the
corresponding factor $\vec{\xi}_i$. Consider the
$\vec{\xi}$-contraction $C^{r,\alpha}_{g^n}({\psi}_1,\dots ,
{\psi}_Z,\vec{\xi})$ which arises when ${\nabla}_i$ hits the
corresponding factor $T$ as before. It is then clear that each
linear combination $C^{l,\alpha}_{g^n}({\psi}_1,\dots ,
{\psi}_Z,\vec{\xi})+ {\Sigma}_{r\in
R^l}C^{r,\alpha}_{g^n}({\psi}_1,\dots ,{\psi}_Z,\vec{\xi})$ is of
the form of equation (\ref{sxesh}). Notice that by Observation
\ref{toremark} we are free to impose this restriction on the order
of integrations by parts of the factors $\vec{\xi}$ in
$\Sigma_{r\in R^l}C^{r,\alpha}_{g^n}({\psi}_1,\dots
,{\psi}_Z,\vec{\xi})$. We note that the order in which we
integrate by parts is consistent with our rules on dropping
$\vec{\xi}$-contractions into the sublinear combinations $PO[\dots
], H[\dots ], Stig^1[\dots ]$, $Stig^2[\dots ]$. This will follow
from the arguments below.
\newline

\par Now let us consider any $\vec{\xi}$-contraction that
arises in the integration by parts of
$C^l_{g^n}({\psi}_1,\dots ,{\psi}_Z,\vec{\xi})$ when ${\nabla}_i$
hits a factor $T=R$. Let us restrict our attention to the
$\vec{\xi}$-contractions $C^r_{g^n}({\psi}_1,\dots ,
{\psi}_Z,\vec{\xi})$, $r\in R^l$, which arise from
$C^l_{g^n}({\psi}_1,\dots ,{\psi}_Z,\vec{\xi})$
by leaving the factor $T(=R$) unaltered. Suppose their
 index set is $R^l_{\alpha,+}$. We then observe that
 the linear combination  $C^{l,\alpha}_{g^n}({\psi}_1,\dots ,
{\psi}_Z,\vec{\xi})+ {\Sigma}_{r\in R^l_{\alpha,+}}
C^{r,\alpha}_{g^n}({\psi}_1,\dots ,{\psi}_Z,\vec{\xi})$
is of the form (\ref{sxesh}).

\par Finally, consider the $\vec{\xi}$-contractions
$C^{r_1}_{g^n}({\psi}_1,\dots ,{\psi}_Z,\vec{\xi})$,
$C^{r_2}_{g^n}({\psi}_1,\dots ,{\psi}_Z,\vec{\xi})$
 which arise from
$C^l_{g^n}({\psi}_1,\dots ,{\psi}_Z,\vec{\xi})$
by substituting the factor $T(=R$) by
$-2{\nabla}^i\vec{\xi}_i$  and $-|\vec{\xi}|^2$ respectively.
Also, define $R^l_1,R^l_2\subset R^l$ to be the index sets of
all the $\vec{\xi}$-contractions
$C^r_{g^n}({\psi}_1,\dots ,{\psi}_Z,\vec{\xi})$ which arise from
$C^l_{g^n}({\psi}_1,\dots ,{\psi}_Z,\vec{\xi})$
by substituting the factor $T(=R$) by
$-2{\nabla}^i\vec{\xi}_i$  and $-|\vec{\xi}|^2$, respectively,
and by substituting at least one more factor $R$. We then
consider the $\vec{\xi}$-contractions
$C^{r_1,\alpha}_{g^n}({\psi}_1,\dots ,{\psi}_Z,\vec{\xi})$ and
 $C^{r_2,\alpha}_{g^n}({\psi}_1,\dots ,{\psi}_Z,\vec{\xi})$
which arise
from the integration by parts of
$C^{r_1}_{g^n}({\psi}_1,\dots ,{\psi}_Z,\vec{\xi})$ and
 $C^{r_2}_{g^n}({\psi}_1,\dots ,{\psi}_Z,\vec{\xi})$,
 respectively, when ${\nabla}_i$ hits the factors
$-2{\nabla}^i\vec{\xi}_i$, $-|\vec{\xi}|^2$, respectively.
We also consider the $\vec{\xi}$-contractions
$C^{r,\alpha}_{g^n}({\psi}_1,\dots ,{\psi}_Z,\vec{\xi})$ and
 $C^{r,\alpha}_{g^n}({\psi}_1,\dots ,{\psi}_Z,\vec{\xi})$
in the integration by parts of each
$C^r_{g^n}({\psi}_1,\dots ,{\psi}_Z,\vec{\xi})$, $r\in R^l_1$
 or $r\in R^l_2$ when ${\nabla}_i$ hits
 the  factors $-2{\nabla}^i\vec{\xi}_i$
$-|\vec{\xi}|^2$, respectively, which correspond to the factors $-2{\nabla}^i\vec{\xi}_i$ or
$-|\vec{\xi}|^2$ in
$C^{r_1}_{g^n}({\psi}_1,\dots ,{\psi}_Z,\vec{\xi})$ and
 $C^{r_2}_{g^n}({\psi}_1,\dots ,{\psi}_Z,\vec{\xi})$.
It follows by construction that the sublinear combinations
$C^{r_1,\alpha}_{g^n}({\psi}_1,\dots ,{\psi}_Z,\vec{\xi})+
{\Sigma}_{r\in R^l_1}C^{r,\alpha}_{g^n}({\psi}_1,\dots ,{\psi}_Z,\vec{\xi})$
and $C^{r_2,\alpha}_{g^n}({\psi}_1,\dots ,{\psi}_Z,\vec{\xi})+
{\Sigma}_{r\in R^l_2}C^{r,\alpha}_{g^n}({\psi}_1,\dots ,
{\psi}_Z,\vec{\xi})$ are in the form of equation (\ref{sxesh}).

\par Hence, we have shown that if we start with a
 linear combination of $\vec{\xi}$-contractions in the form
(\ref{sxesh}), then for each integration by parts in any
$C^l_{g^n}({\psi}_1,\dots ,{\psi}_Z,\vec{\xi})$, we can
 consider the corresponding integrations by parts of each
$C^r_{g^n}({\psi}_1,\dots ,{\psi}_Z,\vec{\xi})$, $r\in R^l$
and we have that at the next step we will be left with a
 linear combination of $\vec{\xi}$-contractions in the
 form (\ref{sxesh}).

\par If at any stage $C^{l,\alpha}_{g^n}({\psi}_1,\dots ,{\psi}_Z,\vec{\xi})$
 is a complete contraction in the form (\ref{linisym}),
 we put it into $PO[C_{g^n}({\psi}_1,\dots ,{\psi}_Z)]$.
It then also follows that the $\vec{\xi}$-contraction in
${\Sigma}_{r\in R^l} C^{r,\alpha}_{g^n}({\psi}_1,\dots ,
{\psi}_Z,\vec{\xi})$
 which arises from $C^{l,\alpha}_{g^n}({\psi}_1,\dots ,
 {\psi}_Z,\vec{\xi})$ by substituting factors $R$ only by
$-|\vec{\xi}|^2$ is stigmatized of type 1, and it is put
into $Stig^1[C_{g^n}({\psi}_1,\dots ,{\psi}_Z,\vec{\xi})]$.
The $\vec{\xi}$-contraction in
${\Sigma}_{r\in R^l} C^{r,\alpha}_{g^n}({\psi}_1,\dots ,
{\psi}_Z,\vec{\xi})$
which arises from $C^{l,\alpha}_{g^n}({\psi}_1,\dots ,{\psi}_Z,\vec{\xi})$
by substituting factors $R$ only by $-2{\nabla}^i\vec{\xi}_i$ is a hard $\vec{\xi}$-contraction and we put it into
$H[C_{g^n}({\psi}_1,\dots ,{\psi}_Z,\vec{\xi})]$.
 Finally, any $\vec{\xi}$-contraction in
${\Sigma}_{r\in R^l} C^{r,\alpha}_{g^n}({\psi}_1,\dots ,
{\psi}_Z,\vec{\xi})$ which arises from
$C^{l,\alpha}_{g^n}({\psi}_1,\dots ,{\psi}_Z,\vec{\xi})$ by
substituting at least one factor $R$ by $-|\vec{\xi}|^2$ and at
least another factor $R$ by $-2{\nabla}^i\vec{\xi}_i$ is
stigmatized of type 2 and we put it into
$Stig^2[C_{g^n}({\psi}_1,\dots ,{\psi}_Z)]$.

\par Let us also notice that the $\vec{\xi}$-contraction
$C^{r_1,\alpha}_{g^n}({\psi}_1,\dots ,
{\psi}_Z,\vec{\xi})$
will always be undecided (it contains a factor
${\nabla}_i\vec{\xi}_k\vec{\xi}^k$). For the
$\vec{\xi}$-contraction $C^{r_2,\alpha}_{g^n}({\psi}_1,\dots
,{\psi}_Z,\vec{\xi})$
(and also for its followers),
we decompose the factor ${\nabla}_{ik}\vec{\xi}^k$
into $S{\nabla}_{ik}\vec{\xi}^k$ and $Ric_{ik}\vec{\xi}^k$.
 We notice that substituting the factor
${\nabla}_{ik}\vec{\xi}^k$ by $Ric_{ik}\vec{\xi}^k$ will give us
either a good or an undecided $\vec{\xi}$-contraction.

\par Now, let us suppose that the $\vec{\xi}$-contraction
 $C^{r_1,\alpha}_{g^n}({\psi}_1,\dots ,{\psi}_Z,\vec{\xi})$
(after the symmetrization $-2({\nabla}_i{\nabla}^k
\vec{\xi}_k)\longrightarrow -2(S{\nabla}_{ik}\vec{\xi}^k)$)
 is hard. We then observe that the $\vec{\xi}$-contraction
in ${\Sigma}_{r\in R^l_1} C^{r,\alpha}_{g^n}({\psi}_1,\dots ,
{\psi}_Z,\vec{\xi})$ which arises from
$C^{r_1,\alpha}_{g^n}({\psi}_1,\dots ,{\psi}_Z,\vec{\xi})$ by
performing an integration by parts of $\vec{\xi}_i$ and hitting
$-2\nabla^k\vec{\xi}_k$ and symmetrizing by
$-2({\nabla}_i{\nabla}^k \vec{\xi}_k)\longrightarrow
-2(S{\nabla}_{ik}\vec{\xi}^k)$ and then by
 substituting factors $R$
only by $-2{\nabla}^i\vec{\xi}_i$ is also hard. Furthermore, any
$\vec{\xi}$-contraction which
 arises from $C^{l,\alpha}_{g^n}({\psi}_1,\dots ,
{\psi}_Z,\vec{\xi})$ or from
$C^{r_1,\alpha}_{g^n}({\psi}_1,\dots ,{\psi}_Z,\vec{\xi})$
by substituting factors $R$ only by
 $-|\vec{\xi}|^2$ is stigmatized of type 2.

\par So we notice that for each $\vec{\xi}$-contraction
that we put into $PO[C_{g^n}({\psi}_1,\dots ,{\psi}_Z)]$ or
$H[C_{g^n}({\psi}_1,\dots ,{\psi}_Z)]$, the
$\vec{\xi}$-contractions that we will put into
$Stig^1[C_{g^n}({\psi}_1,\dots ,{\psi}_Z)]$ or
$Stig^2[C_{g^n}({\psi}_1,\dots ,{\psi}_Z)]$
will be of the form described in (\ref{stig1}) and
(\ref{stig2}).

\par We have shown our Lemma. $\Box$
\newline

\par We now want to apply the above Lemma in order to prove
 equations (\ref{Laplbreak1}) and (\ref{Laplbreak2}).

\par We make a notational convention: Given any
contraction $C^z_{g^n}({\psi}_1,\dots ,{\psi}_Z)$ in the form
(\ref{linisym}), let us write it as
${C^z}'_{g^n}({\psi}_1,\dots ,{\psi}_Z)\cdot R^{\alpha}$, where
${C^z}'_{g^n}({\psi}_1,\dots ,{\psi}_Z)$ does not contain factors $R$.

We then define:

\begin{equation}
\label{wadham1} \begin{split} &{\Sigma}_{r\in
R^z}C^r_{g^n}({\psi}_1,\dots ,{\psi}_Z, \vec{\xi}) =
{C^z}'_{g^n}({\psi}_1,\dots ,{\psi}_Z)\cdot[G(R,\alpha
,-2{\nabla}^i\vec{\xi}_i)+ \\&T(R,\alpha
,-2{\nabla}^i\vec{\xi}_i,-|\vec{\xi}|^2)]
\end{split}
\end{equation}

In the above, each summand on the right hand side arises from one
of
 the substitutions described in the definitions of
$G(R,\alpha, -2{\nabla}^i\vec{\xi}_i)$ and
$T(R,\alpha ,-2{\nabla}^i\vec{\xi}_i,-|\vec{\xi}|^2)$.

Also, given any hard $\vec{\xi}$-contraction
$C^h_{g^n}({\psi}_1,\dots ,{\psi}_Z,\vec{\xi})$ in the form
(\ref{linisymxi2}), we write it as
${C^h}'_{g^n}({\psi}_1,\dots ,{\psi}_Z,\vec{\xi})\cdot R^{\alpha}$, where
  ${C^h}'_{g^n}({\psi}_1,\dots ,{\psi}_Z,\vec{\xi})$ is in
 the form (\ref{linisymxi2}) and does not contain factors
$R$. We then define:

\begin{equation}
\label{wadham2}
\begin{split}
&{\Sigma}_{w\in W^h} C^w_{g^n}({\psi}_1,\dots ,{\psi}_Z,
\vec{\xi})= {C^h}'_{g^n}({\psi}_1,\dots ,{\psi}_Z)\cdot
[G(R,\alpha, -2{\nabla}^i\vec{\xi}_i)+
\\& T(R,\alpha
,-2{\nabla}^i\vec{\xi}_i,-|\vec{\xi}|^2)]
\end{split}
\end{equation}
In the above, each summand on the right hand side arises from one
of
 the substitutions described in the definition of
$G(R,\alpha, -2{\nabla}^i\vec{\xi}_i)$ and
$T(R,\alpha ,-2{\nabla}^i\vec{\xi}_i,-|\vec{\xi}|^2)$.

\par We now prove equations
(\ref{Laplbreak1}) and (\ref{Laplbreak2})
 through an inductive argument. We first recall the terminology and notation used
 in the previous Lemma. Consider
$$I^Z_{g^n}({\psi}_1,\dots ,{\psi}_Z)={\Sigma}_{a\in A} b_a
 C^a_{g^n}({\psi}_1,\dots ,{\psi}_Z)$$
For any complete contraction
$C^a_{g^n}({\psi}_1,\dots ,{\psi}_Z)$ consider the
 sublinear combination of its
 good, hard or undecided descendants, say
${\Sigma}_{x\in X^a} c_x C^x_{g^n}({\psi}_1,\dots
,{\psi}_Z,\vec{\xi})$. We perform integrations by parts in the
expression

$${\int}_{\mathbb{R}^N} e^{(N-n)\vec{\xi}\cdot\vec{x}}\{
{\Sigma}_{a\in A} b_a [C^a_{g^n}({\psi}_1,\dots ,{\psi}_Z)+
{\Sigma}_{x\in X^a} c_x C^x_{g^n}({\psi}_1,\dots ,{\psi}_Z,\vec{\xi})]\}dV_{g^N}=0$$

as explained in the previous Lemma. Whenever we encounter hard or stigmatized $\vec{\xi}$-contractions, we stop (and do
 not discard). In the end, we are
left with a linear combination of sums of complete
 contractions:

\begin{equation}
\label{writeout}
\begin{split}
&{\int}_{\mathbb{R}^N}e^{(N-n)\vec{\xi}\cdot\vec{x}}\{
{\Sigma}_{z\in Z} Q^z(N) a_z C^z_{g^N}
({\psi}_1,\dots ,{\psi}_Z) + {\Sigma}_{h\in H} Q^h(N) a_h
C^h_{g^N}({\psi}_1,\dots ,{\psi}_Z,\vec{\xi})
\\& +{\Sigma}_{z\in Z} Q^z(N) a_z [ {\Sigma}_{r\in R^Z}
C^r_{g^N}({\psi}_1,\dots ,{\psi}_Z,\vec{\xi})] +
\\&{\Sigma}_{h\in H} Q^h(N) a_h [{\Sigma}_{w\in W^h} C^w_{g^N}
({\psi}_1,\dots ,{\psi}_Z,\vec{\xi})]\} dV_{g^N}
=0
\end{split}
\end{equation}

In the above equation each rational function has degree zero
 and leading order coefficient $1$. Moreover, ${\Sigma}_{z\in Z} Q^z(N) a_z C^z_{g^N}
({\psi}_1,\dots ,{\psi}_Z)$ is the sublinear combination that is
 dropped into $PO[\dots ]$, while
${\Sigma}_{h\in H} Q^h(N) a_h C^h_{g^N}({\psi}_1,\dots
,{\psi}_Z,\vec{\xi})$ is the
 sublinear combination that arises by summing over all the sublinear combinations
 of hard $\vec{\xi}$-contractions in the form
 $\Sigma_{f\in F} a_f C^f_{g^n}(\psi_1,\dots ,\psi_Z,\vec{\xi})$
 on the right hand side of (\ref{h}). Then
${\Sigma}_{r\in R^Z} Q^z(N) C^r_{g^N}({\psi}_1,\dots
,{\psi}_Z,\vec{\xi})$ and
\\${\Sigma}_{w\in W^h} Q^h(N) C^w_{g^N}({\psi}_1,\dots ,{\psi}_Z,
\vec{\xi})$ are the sublinear combinations of hard and stigmatized (of
 both types) $\vec{\xi}$-contractions that arise from
${\Sigma}_{z\in Z} a_z C^z_{g^N}({\psi}_1,\dots ,{\psi}_Z)$ and
${\Sigma}_{h\in H} a_h
C^h_{g^N}({\psi}_1,\dots ,{\psi}_Z,\vec{\xi})$ respectively, by
 performing the
 substitutions for the factors $R$ that are explained in
(\ref{wadham1}), (\ref{wadham2}).

\par Our inductive assumption is the following:
For any $T$, We define $Z^T\subset Z$ to be the index set of complete
 contractions $C^z_{g^N} ({\psi}_1,\dots ,{\psi}_Z)$ with $T$
 factors $R$. Furthermore,
we define $Z^{|T}$ to be the index set of complete contractions
$C^z_{g^N} ({\psi}_1,\dots ,{\psi}_Z)$ with more than $T$ factors $R$. We also define $H^T\subset Z$ to be
be the index set of complete contractions $C^z_{g^N}
({\psi}_1,\dots ,{\psi}_Z)$ with $T$ factors $R$.
Also, we define
$H^{|T}\subset H$ to be the index set of
$\vec{\xi}$-contractions $C^h_{g^N}({\psi}_1,\dots ,{\psi}_Z
,\vec{\xi})$ with more than $T$ factors $R$.
We now inductively assume that for some $T$:

\begin{equation}
\label{indT1}
{\Sigma}_{z\in Z^{|T}} a_z C^z_{g^N}
({\psi}_1,\dots ,{\psi}_Z)=0
\end{equation}

and

\begin{equation}
\label{indT2}
{\Sigma}_{h\in H^{|T}} a_h
C^h_{g^N}({\psi}_1,\dots ,{\psi}_Z,\vec{\xi})=0
\end{equation}

We furthermore assume that:

\begin{equation}
\label{indT3}{\Sigma}_{z\in Z^{|T}} a_z[{\Sigma}_{r\in R^z}
C^r_{g^N} ({\psi}_1,\dots ,{\psi}_Z,\vec{\xi})]=0
\end{equation}

and also that:

\begin{equation}
\label{indT4}
{\Sigma}_{h\in H^{|T}} a_h
[{\Sigma}_{w\in W^z}
C^w_{g^N} ({\psi}_1,\dots ,{\psi}_Z,\vec{\xi})]=0
\end{equation}

\par Our goal will be to prove:

\begin{equation}
\label{corindT1}
{\Sigma}_{z\in Z^T} a_z C^z_{g^N}
({\psi}_1,\dots ,{\psi}_Z)=0
\end{equation}

\begin{equation}
\label{corindT2}{\Sigma}_{h\in H^T} a_h
C^h_{g^N}({\psi}_1,\dots ,{\psi}_Z,\vec{\xi})=0
\end{equation}

and furthermore:

\begin{equation}
\label{corindT3}
{\Sigma}_{z\in Z^T} a_z [
{\Sigma}_{r\in R^z}
C^r_{g^N}({\psi}_1,\dots ,{\psi}_Z,\vec{\xi})]=0
\end{equation}

\begin{equation}
\label{corindT4}
{\Sigma}_{h\in H^T} a_z[{\Sigma}_{r\in R^h}
C^r_{g^N} ({\psi}_1,\dots ,{\psi}_Z,\vec{\xi})]=0
\end{equation}

\par We first state and prove a Lemma that will be useful for this purpose:

\begin{lemma}
\label{killRs} Suppose we are given a set of hard
$\vec{\xi}$-contractions $\{C^l_{g^n}({\psi}_1,\dots
,{\psi})_Z,\vec{\xi})\}_{l\in L}$, each in the form
(\ref{linisymxi2}) with $k=0$ (meaning no
 factors $\vec{\xi}$) and of weight $-n$. We suppose that:

\begin{equation}
\label{toana} {\Sigma}_{l\in L} a_l C^l_{g^n}({\psi}_1,\dots
,{\psi}_Z, \vec{\xi})=0
\end{equation}

for every $(M^n,g^n)$, for every ${\psi}_1,\dots ,{\psi}_Z \in
C^{\infty}(M^n)$ and every coordinate system.

\par We define the subsets $L^m\subset L$ as follows:
$l\in L^m$ if and only if $C^l_{g^n}({\psi}_1,\dots , {\psi}_Z,
\vec{\xi})$ has $m$ factors $R$. We then claim that for each $L^m$
for which $L^m\ne \emptyset$:

\begin{equation}
\label{gettoana} {\Sigma}_{l\in L^m} a_l C^l_{g^n}({\psi}_1,\dots
,{\psi}_Z, \vec{\xi})=0
\end{equation}

\par The same result is true if we have complete contractions
$C^l_{g^n}({\psi}_1,\dots ,{\psi}_Z)$ instead of hard
$\vec{\xi}$-contractions $C^l_{g^n}({\psi}_1,\dots
,{\psi}_Z,\vec{\xi})$.
\end{lemma}

{\it Proof:} We will think of the $\vec{\xi}$-contractions
$C^l_{g^n}({\psi}_1,\dots ,{\psi}_Z,\vec{\xi})$ as being in
 the form (\ref{partlinisymxi}). It is straightforward to
 notice that any $\vec{\xi}$-contraction in the form
(\ref{linisymxi2}) with $m$ factors $R$ will give rise to
 $\vec{\xi}$-contractions in the form (\ref{partlinisymxi})
 with $m$ factors $R$.

\par Let us suppose that for some $M>0$ we have that for each $\mu>M$:

$$\Sigma_{l\in L^\mu} a_lC^l_{g^n}(\psi_1,\dots ,\psi_s,\vec{\xi})=0$$

Notice that if we can prove that:

\begin{equation}
\label{toanaM} {\Sigma}_{l\in L^M} a_l C^l_{g^n}({\psi}_1,\dots
,{\psi}_Z, \vec{\xi})=0
\end{equation}

 then the
whole Lemma will follow by induction. In view of our induction
 hypothesis, we erase the sublinear combination
$\Sigma_{\mu >M} \Sigma_{l\in L^\mu} a_lC^l_{g^n}(\psi_1,\dots ,
\psi_s,\vec{\xi})$ from (\ref{toana}).

\par Recall that equation (\ref{toana}) holds for any
Riemannian metric, any functions ${\psi}_1,\dots {\psi}_Z$,
 any coordinate system and any $\vec{\xi}$.
 Hence, equation (\ref{toana})
 must hold {\it formally}.

\par If we can prove that the number of factors $R$ in a complete contraction
 of the form (\ref{partlinisymxi}) remains invariant under the permutations of Definition
\ref{realform}, we will have our Lemma.

\par For any complete contraction
$C^l_{g^n}({\psi}_1,\dots ,{\psi}_Z,\vec{\xi})$ of the form
(\ref{linisymxi2}), we will call one of its
 factors
${\nabla}^m_{r_1\dots r_m}R_{ijkl}$ {\it connected} if one of the
indices $r_1,\dots ,l$ contracts against
 another factor in
$C^l_{g^n}({\psi}_1,\dots ,{\psi}_Z,\vec{\xi})$. By virtue of the
identities in Definition \ref{realform}, we observe
 that any permutation of indices in any connected factor
${\nabla}^m_{r_1\dots r_m}R_{ijkl}$ in $C^l_{g^n}({\psi}_1,\dots
,{\psi}_Z,\vec{\xi})$ will give rise to a complete contraction
$C^{l'}_{g^n}({\psi}_1,\dots ,{\psi}_Z,\vec{\xi})$, which is
 obtained from $C^l_{g^n}({\psi}_1,\dots ,{\psi}_Z,\vec{\xi})$
by substituting its factor ${\nabla}^m_{r_1\dots r_m}R_{ijkl}$ by
a number of
 factors ${\nabla}^pR_{ijkl}$, each of which is connected in
$C^{l'}_{g^n}({\psi}_1,\dots ,{\psi}_Z,\vec{\xi})$.

\par For any complete contraction of the form
$C^l_{g^n}({\psi}_1,\dots ,{\psi}_Z,\vec{\xi})$, we will call one
of its factors ${\nabla}^m_{r_1\dots r_m}R_{ijkl}$
$m$-self-contained if all the indices $r_1,\dots ,l$ contract
against another index in ${\nabla}^m_{r_1\dots r_m}R_{ijkl}$. We
now observe that any application of the identities of Definition
\ref{realform} to a factor ${\nabla}^m_{r_1\dots r_m}R_{ijkl}$
will give rise to a complete contraction
$C^{l'}_{g^n}({\psi}_1,\dots ,{\psi}_Z,\vec{\xi})$, which is
 obtained from $C^l_{g^n}({\psi}_1,\dots ,{\psi}_Z,\vec{\xi})$
by substituting its factor ${\nabla}^m_{r_1\dots r_m}R_{ijkl}$ by
a number of
 factors ${\nabla}^pR_{ijkl}$, each of which is either
$m$-self-contained or connected in
$C^{l'}_{g^n}({\psi}_1,\dots
,{\psi}_Z,\vec{\xi})$.

\par Hence we have shown our Lemma. $\Box$
\newline

\par We now show (\ref{corindT1}).We observe that if
 a complete contraction $C^z_{g^N} ({\psi}_1,\dots ,{\psi}_Z)$
has $\gamma$ factors $R$, then each
$\vec{\xi}$-contraction
$C^r_{g^N} ({\psi}_1,\dots ,{\psi}_Z,\vec{\xi})$ with
$r\in R^z$ has strictly less than $\gamma$
 factors $R$. Furthermore, if
$C^h_{g^N} ({\psi}_1,\dots ,{\psi}_Z,\vec{\xi})$
has $\epsilon$ factors $R$ then each
$C^w_{g^N} ({\psi}_1,\dots ,{\psi}_Z), w\in W^h$
 has strictly less than $\epsilon$ factors $R$.
 Finally, we notice that along the iterative integrations by
 parts the number of factors $R$ either decreases or remains the
 same; it cannot increase.
Now, we want to apply Lemma \ref{finecanc} and (\ref{todo}) to the
case at hand. For any $\vec{\xi}$-contraction
$C_{g^n}({\psi}_1,\dots ,{\psi}_Z,\vec{\xi})$, we have defined
$O[C_{g^n}({\psi}_1,\dots ,{\psi}_Z,\vec{\xi}]$ to stand for its
outgrowth. We also define $H[C_{g^n}({\psi}_1,\dots
,{\psi}_Z,\vec{\xi})]$ to stand for the sublinear combination of
the hard $\vec{\xi}$-contractions that arise along its iterative
 integration by parts. We then re-express the equation in
 Proposition \ref{2ndref} as follows:

\begin{equation}
\label{out3ana}
\begin{split}
&{\Sigma}_{m=0}^T \{ {\Sigma}_{z\in Z^m} a_z C^z_{g^n}
({\psi}_1,\dots ,{\psi}_Z,\vec{\xi}) + {\Sigma}_{r\in R^z} O[
C^r_{g^n}({\psi}_1,\dots ,{\psi}_Z,\vec{\xi})]\} +
\\& {\Sigma}_{m=0}^T \{ {\Sigma}_{h\in H^m} a_h
{\Sigma}_{w\in W^h} O[C^w_{g^n}
({\psi}_1,\dots ,{\psi}_Z)]\}=0
\end{split}
\end{equation}

\par  Let us consider the sublinear combination of complete
contractions in (\ref{out3ana}) with $T$ factors
$R$. It follows from our reasoning above and from Lemma
\ref{xiandr} that it is
 precisely the left hand side of (\ref{corindT1}).
Hence, invoking Lemma \ref{killRs}, we will have
(\ref{corindT1}). Therefore, by the construction of
${\Sigma}_{r\in R^z}
C^r_{g^n}({\psi}_1,\dots ,{\psi}_Z,\vec{\xi})$, we obtain
(\ref{corindT3}).

\par Furthermore, we re-express (\ref{todo}) as follows:

\begin{equation}
\label{hard3ana}
\begin{split}
&{\Sigma}_{m=0}^T \{ {\Sigma}_{h\in H^m} a_h C^h_{g^n}
({\psi}_1,\dots ,{\psi}_Z) + {\Sigma}_{w\in W^h} H[
C^w_{g^n}({\psi}_1,\dots, {\psi}_Z,\vec{\xi})] +
\\& {\Sigma}_{m=0}^T \{ {\Sigma}_{z\in Z^m} a_z [
{\Sigma}_{r\in R^z} H[C^r_{g^n}
({\psi}_1,\dots ,{\psi}_Z,\vec{\xi})]\}=0
\end{split}
\end{equation}

\par Now, we consider the sublinear combination of
$\vec{\xi}$-contractions in the above equation with $T$
 factors $R$. From our reasoning
above, from Lemma \ref{xiandr} and
 also from equation (\ref{corindT1}), we have that
that sublinear combination is precisely the left hand side
 (\ref{corindT2}). Hence, invoking Lemma \ref{killRs},
we have (\ref{corindT2}).
Finally, (\ref{corindT4}) follows from
(\ref{corindT2}) and from its definition.

\par Hence, in view of (\ref{corindT1}) and Lemma \ref{htostig}, we obtain
(\ref{Laplbreak1}), (\ref{Laplbreak2}).

This completes the proof of our Proposition
 \ref{bigcanc}. $\Box$
\newline

We now state a fact that illustrates its usefulness.

\begin{lemma}
\label{lxitol}
Consider  a good or undecided or hard $\vec{\xi}$-contraction
$C_{g^n}({\psi}_1,\dots ,{\psi}_Z,\vec{\xi})$, of
$\vec{\xi}$-length $L$. We then have that
$PO[C_{g^n}({\psi}_1,\dots ,{\psi}_Z,\vec{\xi})]$ will consist
of complete contractions of length greater than
 or equal to $L$, or
$PO[C_{g^n}({\psi}_1,\dots ,{\psi}_Z,\vec{\xi})]=0$.

\par We also consider the hard or the stigmatized
$\vec{\xi}$-contractions that arise along the iterative
 integrations by parts. We claim that any such
$\vec{\xi}$-contraction has $\vec{\xi}$-length $\ge L$.
\end{lemma}

{\it Proof:} The proof is by induction. Initially, to make
things easier, consider the case where the are no factors
$|\vec{\xi}|^2$ in
$C_{g^n}({\psi}_1,\dots ,{\psi}_Z,\vec{\xi})$. So, think of
$C_{g^n}({\psi}_1,\dots ,{\psi}_Z,\vec{\xi})$ as being in the
 form (\ref{linisymxi2}) with $\vec{\xi}$-length $M$ and
with $E$ factors $\vec{\xi}$ and $C$
 factors $S{\nabla}^m \vec{\xi}$. We will perform induction
 on $C+E$.

\par Initially suppose $C+E=1$. Then if $C=1,E=0$, our
$\vec{\xi}$-contraction is hard, so
$PO[C_{g^n}({\psi}_1,\dots ,{\psi}_Z,\vec{\xi})]=0$.
 If $E=1,C=0$, the proof is the same as for the inductive step:

\par Suppose we know the claim is true for $E+C=p$ and we
 want to prove it for $E+C=p+1$. So, pick out a factor
$\vec{\xi}_i$ and do an integration by parts with respect to
 it. If ${\nabla}_i$ hits a factor ${\nabla}^mR_{ijkl}$
or ${\nabla}^pRic_{ij}$ or $\nabla^p\psi_k$, we get a
$\vec{\xi}$-contraction in the form
(\ref{linisymxi1}) or (\ref{linisymxi2}) with
 $E+C=p$ and $\vec{\xi}$-length $M$.
If ${\nabla}_i$ hits a factor $\vec{\xi}$, we get a
$\vec{\xi}$-contraction in the form
(\ref{linisymxi1}) or (\ref{linisymxi2}) with
 $E+C=p$ and $\vec{\xi}$-length $M+1$.
If it hits a factor $S{\nabla}^m\vec{\xi}$ ($m\ge 1$), then
 after applying identity
(\ref{symunsymgenxi}), we obtain a linear combination
of complete contractions in the form (\ref{linisymxi1}) or
(\ref{linisymxi2}) with $C+E=p$ and $\vec{\xi}$ length
$\ge M$.

\par Now, suppose we do allow factors
$|\vec{\xi}|^2$ in
$C_{g^n}({\psi}_1,\dots ,{\psi}_Z,\vec{\xi})$.
We again proceed by induction on the number $C+E$.
 If all the $\vec{\xi}$-factors in
$C_{g^n}({\psi}_1,\dots ,{\psi}_Z,\vec{\xi})$ are in the form
  $|\vec{\xi}|^2$ or $S{\nabla}^m\vec{\xi}$,
 we already have a stigmatized
$\vec{\xi}$-contraction. Hence, $PO[C_{g^n}({\psi}_1,\dots
,{\psi}_Z,\vec{\xi})]=0$ in that case. Otherwise, there is at
least one factor $\vec{\xi}_i$ that does not contract against
another factor $\vec{\xi}$. We integrate by parts with respect to
it. If ${\nabla}_i$ hits a factor ${\nabla}^mR_{ijkl}$ or
${\nabla}^aRic_{ij}$ or $\vec{\xi}$ or $|\vec{\xi}|^2$, we fall
under our induction hypothesis with $\vec{\xi}$-length $M$ or
$M+1$. If it hits a factor $S{\nabla}^m\vec{\xi}$,
 we apply (\ref{symunsymgenxi})  and obtain a
 a linear combination of $\vec{\xi}$-contractions that fall
 under our induction hypothesis, by the same reasoning as above.
 That completes the proof of the Lemma.$\Box$

\subsection{Conclusion: The Algorithm for the super
 divergence formula.}

\par We want to apply Proposition \ref{bigcanc}
and see how it can provide us with a divergence formula for
 $I^Z_{g^n}({\psi}_1,\dots ,{\psi}_Z)$.
The Proposition gives us an algorithm:

\par Write $I^Z_{g^n}({\psi}_1,\dots ,
{\psi}_Z)=
{\Sigma}_{r\in R} a_r C^r_{g^n}({\psi}_1,\dots ,{\psi}_Z)$,
where each complete contraction $C^r_{g^n}({\psi}_1,\dots ,{\psi}_Z)$
 is in the form (\ref{linisym}).

\par For each
 complete contraction $C^r_{g^n}({\psi}_1,\dots ,{\psi}_Z)$
 we consider the set of its good or undecided
descendants, along with their coefficients (see definition
\ref{descendants}),
 say $a_b C^{r,b}_{g^n}({\psi}_1,\dots ,{\psi}_Z,
\vec{\xi}), b\in B^r$. So each
$C^{r,b}_{g^n}({\psi}_1,\dots ,{\psi}_Z, \vec{\xi})$ is in
 the form (\ref{linisymxi1}) or (\ref{linisymxi2}) and
has $S_b$ $\vec{\xi}$-factors (see definition \ref{xi}).

\par We then begin to integrate by parts each
$\vec{\xi}$-contraction $C^{r,b}_{g^n}({\psi}_1,\dots ,
{\psi}_Z,\vec{\xi})$, and make the following convention:
Whenever along this iterative integration by parts we obtain
a hard or a stigmatized $\vec{\xi}$-contraction (see definition \ref{defstigma}),
 we discard it. For each $\vec{\xi}$-contraction
$C^{r,b}_{g^n}({\psi}_1,\dots ,{\psi}_Z,\vec{\xi})$,
 consider the $\vec{\xi}$-contractions
$a_x C^{r,b,x}_{g^n}({\psi}_1,\dots ,{\psi}_Z,\vec{\xi}),
x\in X^b$ we are
 left with after $S_b-1$ integrations by parts
(along with their coefficients). They are in
 the form (\ref{linisymxi1}) with one factor $\vec{\xi}$.

\par We then construct vector fields
$(C^{r,b,x}_{g^n})^j({\psi}_1,\dots ,{\psi}_Z)$ out of each
\\ $C^{r,b,x}_{g^n}({\psi}_1,\dots ,{\psi}_Z,\vec{\xi})$
by crossing out the factor $\vec{\xi}_j$ and making the index
 that contracted against $j$ into a free index. By virtue of
 Proposition \ref{bigcanc}, we have:

\begin{equation}
\label{supdv1} I^Z_{g^n}({\psi}_1,\dots ,{\psi}_Z)= {\Sigma}_{r\in
R} a_r {\Sigma}_{b\in B^r} a_b{\Sigma}_{x\in X^b}div_j a_x
(C^{r,b,x}_{g^n})^j(\psi_1,\dots ,\psi_Z)
\end{equation}

\par We will refer to this equation as the
{\it super divergence formula} and denote it by
$supdiv[I^Z_{g^n}(\psi_1,\dots ,\psi_Z)]=0$. We note that there
are many such formulas, since at each stage we pick a factor
$\vec{\xi}$ to
integrate by parts (subject to the restrictions that we have imposed
because of Remark \ref{toremark}).
\newline

\par Now, a notational convention and two observations:
Firstly, for any complete contraction
$C^r_{g^n}(\psi_1,\dots ,\psi_Z)$, define:

\begin{equation}
\label{deftail} Tail[C^r_{g^n}({\psi}_1,\dots ,{\psi}_Z)]=
C^{r}_{g^n} ({\psi}_1,\dots ,{\psi}_Z) +{\Sigma}_{b\in B_r} a_b
PO[C^{r,b}_{g^n} ({\psi}_1,\dots ,{\psi}_Z,\vec{\xi})]
\end{equation}

\par Then, notice that if the complete contraction
$C^r_{g^n}({\psi}_1,\dots ,{\psi}_Z)$ has length $L$, then
each complete contraction in its tail will have length
$\ge L$. This follows from Lemmas \ref{lxitol} and
\ref{Itosupdiv}.

\par Furthermore, we see that the super divergence formula
holds for any
\\ $I^Z_{g^n}({\psi}_1,\dots ,{\psi}_Z)=
{\Sigma}_{r\in R} a_r C^r_{g^n}({\psi}_1,\dots {\psi}_Z)$
where each complete contraction is in the form (\ref{linisym})
 with weight $-n$, for which
${\int}_{M^n}I^Z_{g^n}({\psi}_1,\dots ,{\psi}_Z)dV_{g^n}=0$
for every compact Riemannian $(M^n,g^n)$ and any
${\psi}_1,\dots ,{\psi}_Z\in C^{\infty}(M^n)$.
In other words, the super divergence formula does not depend on the fact that
$I^Z_{g^n}({\psi}_1,\dots ,{\psi}_Z)$
arises from a polarization of the transformation law of
$P(g^n)$.

\subsection{The shadow formula.}

\par We will draw another conclusion from the Lemma
\ref{finecanc} and Proposition \ref{bigcanc}.

\par As before, write $I^Z_{g^n}({\psi}_1,\dots ,
{\psi}_Z)=
{\Sigma}_{r\in R} a_r C^r_{g^n}({\psi}_1,\dots ,{\psi}_Z)$,
where each complete contraction $C^r_{g^n}({\psi}_1,\dots ,{\psi}_Z)$
 is in the form (\ref{linisym}).

\par For each complete contraction $C^r_{g^n}(\psi_1,\dots ,\psi_Z)$
 we consider the set of its good or undecided or hard
descendants, along with their coefficients (see definition
\ref{descendants}), say $a_b C^{r,b}_{g^n}(\psi_1,\dots ,\psi_Z,
\vec{\xi}), b\in B^r$. So each
$C^{r,b}_{g^n}({\psi}_1,\dots ,{\psi}_Z, \vec{\xi})$ is in
 the form (\ref{linisymxi1}) or (\ref{linisymxi2}) and
has $S_b$ $\vec{\xi}$-factors.

\par We then begin to integrate by parts each
$\vec{\xi}$-contraction $C^{r,b}_{g^n}({\psi}_1,\dots ,
{\psi}_Z,\vec{\xi})$, in the order explained in definition
\ref{po}. We make the following convention:

Whenever we encounter a hard or a stigmatized
$\vec{\xi}$-contraction, we put it aside. Whenever we
 encounter a good $\vec{\xi}$-contraction with $k=1$
(and $l=0$), we discard it.

\par We then consider the set of the hard or stigmatized
$\vec{\xi}$-contractions, along with their coefficients, that
 we are left with after this procedure. Suppose that set is
$\{a_t C^t_{g^n}({\psi}_1,\dots ,{\psi}_Z,
\vec{\xi})\}_{t\in T}$.
 We then have the {\it shadow formula} for
$I^Z_{g^n}({\psi}_1,\dots ,{\psi}_Z)$:

\begin{equation}
\label{shadow}
{\Sigma}_{t\in T} a_t C^t_{g^n}({\psi}_1,\dots ,
{\psi}_Z,\vec{\xi})=0
\end{equation}

for every $(M^n,g^n)$, every ${\psi}_1,\dots ,{\psi}_Z\in
C^{\infty}(M^n)$, any coordinate system and any $\vec{\xi}\in
\mathbb{R}^n$. We will denote this equation by
$Shad[I^Z_{g^n}({\psi}_1,\dots ,{\psi}_Z,\vec{\xi})]=0$. It
 follows, as for the super divergence formula, that
the shadow equation holds for any
$I^Z_{g^n}({\psi}_1,\dots ,{\psi}_Z)$
that integrates to zero on any $(M^n,g^n)$, for any
$\psi_1,\dots ,\psi_Z$ and for any coordinate system and
 any $\vec{\xi}\in \mathbb{R}^n$.
 It does not depend on the fact that
$I^Z_{g^n}(\psi_1,\dots ,\psi_Z)$ is the
 polarized transformation law for some $P(g^n)$.
\newline

\par Recalling the notation of Definition \ref{sieve}, we
additionally define:

\begin{equation}
\label{irish} \begin{split} &Tail^{Shad}[C^r_{g^n}({\psi}_1,\dots
,{\psi}_Z)]= \Sigma_{b\in B^r} a_b \{H[C^b_{g^n} (\psi_1,\dots
,\psi_Z,\vec{\xi})]
\\& +Stig^1 [C^b_{g^n}(\psi_1,\dots
,\psi_Z,\vec{\xi})] + Stig^2[C^b_{g^n}(\psi_1,\dots
,\psi_Z,\vec{\xi})]\}
\end{split}
\end{equation}

We may then re-express the shadow formula as:

\begin{equation}
\label{shadow3ana}
\Sigma_{r\in R} a_rTail^{Shad}[C^r_{g^n}(\psi_1,\dots ,\psi_Z)]=0
\end{equation}

\par The above equation follows straightforwardly from
Lemma \ref{finecanc} and also from equation (\ref{Laplbreak2}).

\par Moreover, for future reference we define:

\begin{equation}
\label{landofmordor}
\begin{split}
&O^{Shad}[C^b_{g^n}(\psi_1,\dots ,\psi_Z,\vec{\xi})]= H[C^b_{g^n}
(\psi_1,\dots ,\psi_Z,\vec{\xi})] +Stig^1
[C^b_{g^n}(\psi_1,\dots ,\psi_Z,\vec{\xi})] +
\\& Stig^2[C^b_{g^n}(\psi_1,\dots ,\psi_Z,\vec{\xi})]
\end{split}
\end{equation}

\par We furthermore show the following: For any
$m\ge 0$, let $T^m$ stand for the sublinear combination in
(\ref{shadow}) with $m$ factors $|\vec{\xi}|^2$. We then also
 have:

\begin{equation}
\label{shadowm}
{\Sigma}_{t\in T^m} a_t C^t_{g^n}({\psi}_1,\dots ,
{\psi}_Z,\vec{\xi})=0
\end{equation}
for every $(M^n,g^n)$, every ${\psi}_1,\dots ,{\psi}_Z\in
C^{\infty}(M^n)$, any coordinate system and any $\vec{\xi}\in
\mathbb{R}^n$.

\par This follows since (\ref{shadow}) must hold formally and
 the number of factors $\vec{\xi}$ that contract against
 another factor $\vec{\xi}$ is invariant under the
 permutations of definition \ref{realform}.

\par Furthermore it follows, from Lemma \ref{Itosupdiv} and
also from Lemma \ref{lxitol}, that
if a complete contraction $C^l_{g^n}(\psi_1,\dots ,
\psi_Z)$ of length $L$
in $I^Z_{g^n}(\psi_1,\dots ,\psi_Z)$
gives rise to a hard or stigmatized
$\vec{\xi}$-contraction $C^{l,z}_{g^n}(\psi_1,\dots ,
\psi_Z,\vec{\xi})$ in (\ref{shadow}), by the procedure outlined
 above, then
$C^{l,z}_{g^n}({\psi}_1,\dots ,{\psi}_Z,\vec{\xi})$
 will have $\vec{\xi}$-length $\ge L$.
\newline

\noindent{\bf Acknowledgement.} The work presented here is
extracted from the author's Ph.D.~dissertation at Princeton
University. I am greatly indebted to my thesis advisor Charles
Fefferman for suggesting this problem to me, for his constant
encouragement and support and for his endless patience.


\begin{thebibliography}{12}



\bibitem{a:rcipem} P. Albin \emph{Renormalizing Curvature
 Integrals on Poincare-Einstein manifolds.} arxiv math.DG/0504161

\bibitem{a:dgciII} S. Alexakis \emph{The decomposition of Global
Conformal Invariants II} to appear in Adv. in Math.

\bibitem{beg:itccg} T. N. Bailey, M. G. Eastwood, C. R. Graham
\emph{Invariant Theory for Conformal and CR Geometry} Ann. of Math (2),
 {\bf 139} (1994), 491-552.

\bibitem{bfg:spdcn} M. Beals, C. Fefferman, R. Grossman
\emph{Strictly Pseudoconvex Domains in $\mathbb{C}^n$}
 Bull. Amer. Math. Soc. (N.S.), {\bf 8} (1983), no.3,
125-322.

\bibitem{be:clwid8} N. Boulager, J. Erdmenger
\emph{A Classification of Local Weyl Invariants in D=8},
Class. Quantum Gravity {\bf 21} (2004), 4305-4316.

\bibitem{b:fd} T. Branson \emph{The functional determinant}, Global
Analysis Research Center Lecture Note Series, no. 4, Seoul
National University (1993).


\bibitem{bgp:ilcfm} T. Branson, P. Gilkey, J. Pohjanpelto
\emph{Invariants of locally conformally flat manifolds}
 Trans. Amer. Math. Soc.  347  (1995),  no. 3, 939--953.

\bibitem{b:gechofcg} S. Brendle \emph{Global existence and
convergence for a higher order flow in conformal geometry},
Ann. of Math.(2) 158 (2003), no.1, 323-343.

\bibitem{cgy:ematcg} S.Y.A. Chang, M. Gursky, P.C. Yang
\emph{An equation of Monge-Ampere type in conformal geometry,
and four-manifolds of positive Ricci curvature} Ann. of Math.(2) 155
 (2002) 709-787.

\bibitem{cqy:tccem} S.Y.A. Chang, J. Qing, P.C. Yang
\emph{On the topology of conformally compact Einstein 4-manifolds},
 Noncompact problems at the intersection of geometry, analysis,
 and topology,  49--61, Contemp. Math., 350.

\bibitem{cqy:pc}S.Y.A. Chang, J. Qing, P.C. Yang
\emph{On the renormalized volumes of conformally compact Einstein
manifolds}, private communication.








\bibitem{ds:gccaad} S. Deser, A. Schwimmer
\emph{Geometric classification of
 conformal anomalies in arbitrary dimensions}, Phys. Lett.
 B309 (1993) 279-284.




\bibitem{e:ncg} M. G. Eastwood \emph{Notes on Conformal Geometry}
Rend. Circ. Mat. Palermo      (2)     Suppl. No.43 (1996), 57-76



\bibitem{e:rg} L.P. Eisenhart \emph{Riemannian Geometry} Princeton University Press (1925)

\bibitem{e:ntrm} D.B.A. Epstein \emph{Natural Tensors on Riemannian
Manifolds}, Journal of Diff. Geom. {\bf 10} (1975), 631-645


\bibitem{fg:ci} C. Fefferman, C. R. Graham \emph{Conformal Invariants}
\'Elie Cartan et les mathematiques d'aujourd'hui,      Ast\'erisque,
numero hors serie, 1985, 95-116.

\bibitem{fh:amcqcccg} C. Fefferman, K. Hirachi \emph{Ambient Metric
Construction of Q-Curvature in Conformal and CR Geometries}, Math. Res.
Lett. {\bf 10} (2003), 819-831.




\bibitem{g:lierm} P. Gilkey \emph{Local Invariants of an Embedded
 Riemannian Manifold.} Ann. of Math. (2) {\bf 102} (1975), no.2.
187-203.



\bibitem{g:varccem} C.R. Graham \emph{Volume and area renormalizations for conformally compact Einstein metrics}
 Rend. Circ. Math. Palermo II. Ser. Suppl. 63, 31-42 (2000).

\bibitem{gjms:cipl} C. R. Graham, R. Jenne, L. J. Mason, G. Sparling
\emph{Conformally invariant powers of the Laplacian: existence},
J. London Math. Soc. (2) {\bf 46} (1992), 557-565.


\bibitem{gw:casoacc} C.R. Graham E. Witten \emph{Conformal
Anomaly of Submanifold Observables in AdS/CFT Correspondence},
Nucl. Phys. B 546 (1999), 52-64.


\bibitem{gz:smcg} C.R. Graham, M. Zworski
\emph{Scattering Matrix in Conformal
Geometry}, Invent. Math. 152 (2003) 89-118.




\bibitem{g:pecidowasep} Gursky, M. \emph{The Principal eigenvalue of
a conformally invariant operator, with an application to semilinear
elliptic PDE}, Comm. Math. Phys. 207
(1999), no.1 131-143.






\bibitem{hs:hwa} M. Henningson, K. Skenderis \emph{The holographic
Weyl anomaly}, J. High Energy Phys. 07 (1998)












\bibitem{q:rccem} J. Qing \emph{On the rigidity for conformally compact
Einstein manifolds}, arXiv:math.DG/0305084.


\bibitem{w:cg} H. Weyl \emph{The Classical Groups}, Princeton University Press 1946.

\bibitem{w:adssh} E. Witten \emph{Anti de Sitter space and
holography}, Adv. Theor. Math. Phys. 2 (1998), 253-291.


\bibitem{wy:cbacc} E. Witten, S.T. Yau \emph{Connectedness of
the boundary in the AdS/CFT correspondence}, Adv. Theor. Math.
Phys. 3 (1999), no.6, 1635-1655.

\end{thebibliography}
\end{document}